\documentclass{amsart}

\usepackage{amsmath,amssymb,thmtools}
\usepackage{graphicx}

\newcommand{\group}[1]{\mathbf{#1}}
\newcommand{\ring}[1]{\mathbb{#1}}
\newcommand{\field}[1]{\mathbb{#1}}
\newcommand{\vecspace}[1]{\mathsf{#1}}
	\renewcommand{\vector}[1]{\mathsf{#1}}
\newcommand{\algebra}[1]{\mathcal{#1}}
\newcommand{\lattice}[1]{\mathsf{#1}}
\newcommand{\domain}[1]{\mathsf{#1}}

\newcommand{\Cat}{\operatorname{Cat}}
\newcommand{\E}{\mathbb{E}}
\newcommand{\id}{\operatorname{id}}
\newcommand{\Hol}{\operatorname{Hol}}
\newcommand{\lis}{\operatorname{lis}}

\newcommand{\End}{\operatorname{End}}
\newcommand{\F}{\mathfrak{F}}
\newcommand{\vac}{\vector{v}_{\emptyset}}
\newcommand{\Tr}{\operatorname{Tr}}
\newcommand{\tr}{\operatorname{tr}}

\title{Three Lectures on Free Probability}

\author{Jonathan Novak\\ \textit{with illustrations by Michael LaCroix}}
\address{Department of Mathematics, Massachusetts Institute of Technology, 77 Massachusetts Avenue, Cambridge, MA 02139-4307}
\email{jnovak@math.mit.edu}

\begin{document}


\maketitle

\setcounter{section}{-1}

\tableofcontents
\pagebreak

\section{Introduction}
These are notes from a three-lecture mini-course on free probability given
at MSRI in the Fall of 2010 and repeated a year later at Harvard.  The lectures were
aimed at mathematicians and mathematical physicists working in combinatorics, probability, and 
random matrix theory.  
The first lecture was a staged rediscovery of free independence from first principles,
the second dealt with the additive calculus of free random variables, 
and the third focused on random matrix models.

Most of my knowledge of free probability was acquired through informal conversations
with my thesis supervisor, Roland Speicher, and while he is an expert in the field the same
cannot be said for me.  
These notes reflect my own limited understanding
and are no substitute for complete and rigorous treatments, such as
Voiculescu, Dykema and Nica \cite{VDN}, Hiai and Petz \cite{HP}, and Nica and Speicher \cite{NS}.
In addition to these sources, the expository articles of Biane \cite{Biane}, Shlyakhtenko \cite{Shly}
and Tao \cite{Tao} are very informative.

I would like to thank the organizers of the MSRI semester ``Random Matrix Theory, Interacting Particle Systems and
Integrable Systems'' for the opportunity to participate as a postdoctoral fellow.  Special thanks are owed to Peter Forrester
for coordinating the corresponding MSRI book series volume in which these notes appear.
I am also grateful to the participants of the Harvard random
matrices seminar for their insightful comments and questions.  

I am indebted to Michael LaCroix for making the illustrations which accompany these notes.

\section{Lecture One: Discovering the Free World}
	
	\subsection{Counting connected graphs}
	Let $m_n$ denote the number of simple, undirected graphs on the vertex set $[n]=\{1,\dots,n\}$.  We have $m_n=2^{n \choose 2}$,
	since each pair of vertices is either connected by an edge or not. 
	A more subtle quantity is the number $c_n$ of connected graphs on $[n]$.  The sequence $(c_n)_{n \geq 1}$ is listed
	as $\mathsf{A01187}$ in Sloane's Online Encyclopedia of Integer Sequences; its first few terms are
	
		\begin{equation*}
			1,\ 1,\ 4,\ 38,\ 728,\ 26\, 704,\ 1\, 866\, 256,\ \dots.
		\end{equation*}
		
	\noindent
	Perhaps surprisingly, there is no closed formula for $c_n$.  However, $c_n$ may be understood in terms of the transparent
	sequence $m_n$ in several ways, each of which corresponds to a combinatorial decomposition. 
	
		\begin{figure}
			\includegraphics{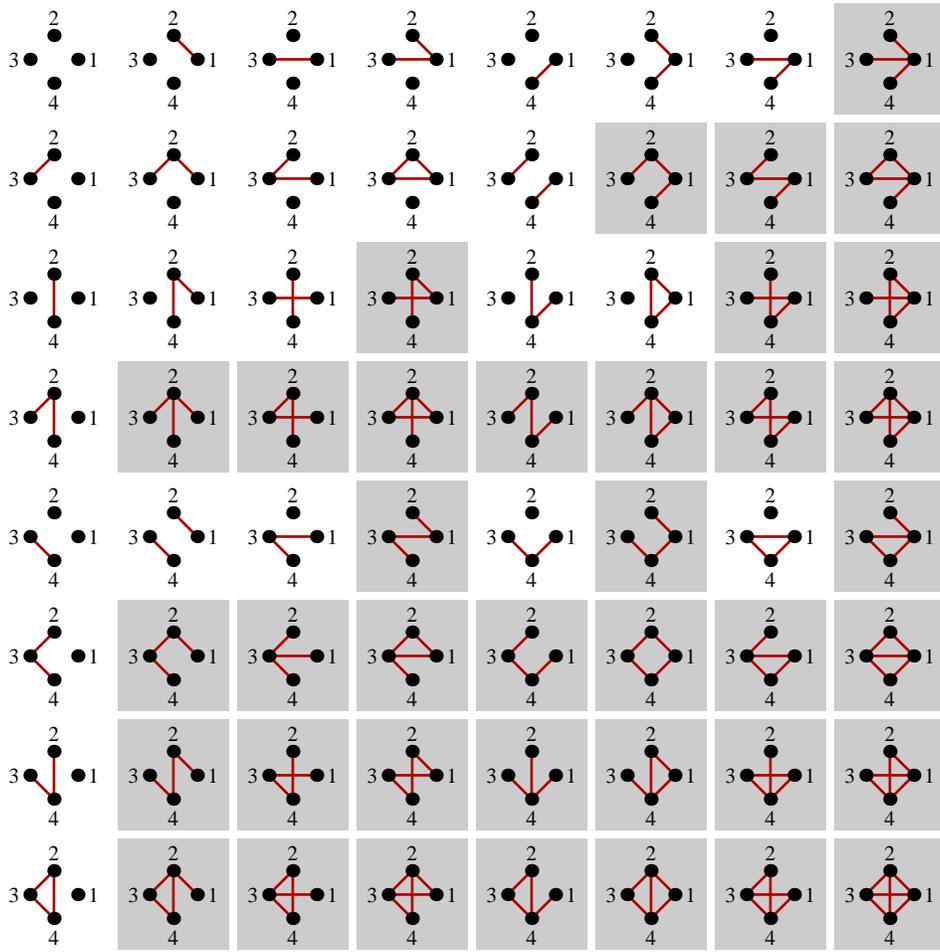}
			\caption{\label{fig:4-vertex}Thirty eight of sixty four graphs on four vertices are connected.}
		\end{figure}
		
	First, we may decompose a graph into two disjoint subgraphs: the connected component of a distinguished vertex, say $n$, and everything else, i.e.
	the induced subgraph on the remaining vertices. 
	Looking at this the other way around, we may build a graph as follows.  
	From the vertices $1,\dots,n-1$ we can choose $k$ of these in ${n-1 \choose k}$ ways, and then build an arbitrary graph on these vertices in 
	$m_k$ ways.  On the remaining $n-1-k$ vertices together with $n$, we may build a connected graph in $c_{n-k}$ ways. 
	This construction produces different graphs for different values of $k$, since the size of the connected component containing 
	the pivot vertex $n$ will be different.  Moreover, as $k$ ranges from one to $n-1$ we obtain all graphs in this fashion.  Thus we have
			
		\begin{equation*}
			\label{eqn:ThieleRecurrence}
			m_n = \sum_{k=0}^{n-1} {n-1 \choose k} m_k c_{n-k},
		\end{equation*}
		
	\noindent
	or equivalently
	
		\begin{equation*}
			c_n = m_n - \sum_{k=1}^{n-1} {n-1 \choose k} m_k c_{n-k}.
		\end{equation*}
		
	\noindent
	While this is not a closed formula, it allows the efficient computation of $c_n$ given $c_1,\dots,c_{n-1}$.
	
	A less efficient but ultimately more useful recursion can be obtained by viewing a graph
	as the disjoint union of its connected components. 
	We construct a graph by first choosing a partition of the underlying vertex set into disjoint non-empty subsets $B_1,\dots,B_k$,
	and then building a connected graph on each of these, which can be done in 
	$c_{|B_1|} \dots c_{|B_k|}$ ways.  This leads to the formula
			
		\begin{equation*}
			\label{eqn:MomentCumulant}
			m_n = \sum_{\pi \in \lattice{P}(n)} \prod_{B \in \pi} c_{|B|},
		\end{equation*}
		
	\noindent
	where the summation is over the set of all partitions of $[n]$.  We can split off the
	term of the sum corresponding to the partition $[n]=[n]$ to obtain the recursion
	
		\begin{equation*}
			c_n = m_n - \sum_{\substack{\pi \in \lattice{P}(n)\\ b(\pi) \geq 2}}  \prod_{B \in \pi} c_{|B|},
		\end{equation*}
		
	\noindent
	in which we sum over partitions with at least two blocks.
			
	The above reasoning is applicable much more generally.  Suppose that $m_n$ is the number of 
	``structures'' which can be built on a set of $n$ labelled points, and that $c_n$ is the number of ``connected structures'' on
	these points of the same type.
	Then the quantities $m_n$ and $c_n$ will satisfy the above (equivalent) relations.	
	This fundamental enumerative link between connected and disconnected structures is ubiquitous in mathematics and the sciences, 
	see \cite[Chapter 5]{Stanley:EC2}.  
	Prominent examples come from enumerative
	algebraic geometry \cite{Roth}, where connected covers of curves are counted in terms of all covers, 
	and quantum field theory \cite{Etingof}, where Feynman diagram sums are reduced to summation over
	connected terms.
		
	\subsection{Cumulants and connectedness}
	The relationship between connected and disconnected structures is well-known to probabilists, albeit from a different
	point of view.  
	In stochastic applications, $m_n=m_n(X)=\E[X^n]$ is the 
	moment
	sequence of a random variable $X$, and the quantities $c_n(X)$ defined by either of the equivalent
	recurrences	
	
		\begin{equation*}
			\begin{split}
				m_n(X) &= \sum_{k=0}^{n-1} {n-1 \choose k} m_k(X) c_{n-k}(X) \\
				m_n(X) &= \sum_{\pi \in \lattice{P}(n)} \prod_{B \in \pi} c_{|B|}(X)
			\end{split}
		\end{equation*}
		
	\noindent
	are called the cumulants of $X$.  
	This term was suggested by Harold Hotelling and subsequently popularized by Ronald Fisher and John Wishart in an influential 1932
	article \cite{FW}.
	Cumulants were, however, investigated as early as 1889 by the Danish mathematician and astronomer Thorvald Nicolai Thiele, 
	who called them half-invariants.  
	Thiele introduced the cumulant sequence as a transform of the moment sequence defined via the first of the above recurrences, 
	and some years later arrived at the equivalent formulation using the second recurrence.  The latter is now called the 
	moment-cumulant formula.  
	Thiele's contributions to statistics and the 
	early theory of cumulants have been detailed by Anders Hald \cite{Hald1,Hald2}.
	
	Cumulants are now well-established and frequently encountered in probability and statistics, sufficiently so that the first four have been given 
	names: mean, variance, skewness, and kurtosis\footnote{In practice, statisticians often define skewness and kurtosis to be 
	the third and fourth cumulants scaled by a power of the variance.}.  
	The formulas for mean and variance in terms of moments are simple and familiar,
	
		\begin{equation*}
			\begin{split}
			c_1(X) &=m_1(X) \\
			c_2(X) &=m_2(X)-m_1(X)^2,
			\end{split}
		\end{equation*}
		
	\noindent
	whereas the third and fourth cumulants are more involved,
	
		\begin{equation*}
			\begin{split}
				c_3(X) &= m_3(X) - 3m_2(X)m_1(X) + 2m_1(X)^3 \\
				c_4(X) &= m_4(X) - 4m_3(X)m_1(X) -3m_2(X)^2 + 12m_2(X)m_1(X)^2 -6m_1(X)^4.
			\end{split}
		\end{equation*}

	It is not immediately clear why the cumulants of a random variable are of interest.  
	If a random variable $X$ is uniquely determined by its moments, then we may think of the moment sequence
	
		\begin{equation*}
			(m_1(X),m_2(X),\dots,m_n(X),\dots)
		\end{equation*}
		
	\noindent
	as coordinatizing $X$.  Passing from moments to cumulants then amounts to a (polynomial) change of coordinates.
	Why is this advantageous?
	
\begin{figure}
			\includegraphics{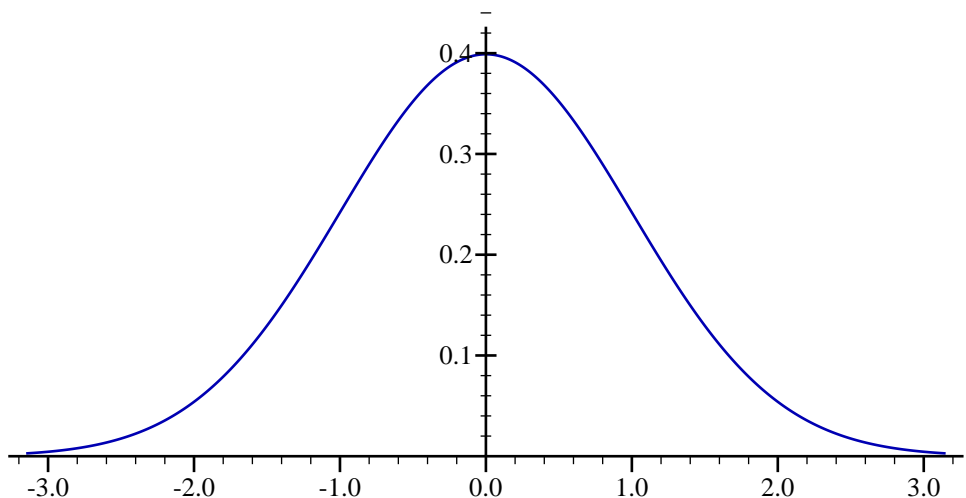}
			\caption{\label{fig:Gaussian}The Gaussian density}
		\end{figure}

	As a motivating example,
	let us compute the cumulant sequence of the most important random variable, 
	the standard Gaussian $X$.
	The distribution of $X$ has density given by the bell curve
	
		\begin{equation*}
			\mu_X(\mathrm{d}t) = \frac{1}{\sqrt{2\pi}} e^{-\frac{t^2}{2}} \mathrm{d}t
		\end{equation*}
		
	\noindent
	depicted in Figure \ref{fig:Gaussian}
	
	\noindent
	We will now determine the moments of $X$.
	Let $z$ be a complex variable, and define
	
		\begin{equation*}
			M_X(z) := \int\limits_{\field{R}} e^{tz} \mu_X(\mathrm{d}t).
		\end{equation*}
		
	\noindent
	Since $e^{-\frac{t^2}{2}}$ decays rapidly as $|t| \rightarrow \infty$, $M_X(z)$ is a well-defined 
	entire function of $z$ 
	whose derivatives can be computed by differentiation under the integral sign,
	
		\begin{equation*}
			M_X'(z) = \int\limits_{\field{R}} te^{tz} \mu_X(\mathrm{d}t),\ 
			M_X''(z) = \int\limits_{\field{R}} t^2e^{tz} \mu_X(\mathrm{d}t),\ \dots.
		\end{equation*}
		
	\noindent
	In particular, the $n^{\text{th}}$ derivative of $M_X(z)$ at $z=0$ is
	
		\begin{equation*}
			M_X^{(n)}(0) = \int\limits_{\field{R}} t^n \mu_X(\mathrm{d}t) = m_n(X),
		\end{equation*}
		
	\noindent
	so we have the Maclaurin series expansion
	
		\begin{equation*}
			M_X(z) = \sum_{n=0}^{\infty} m_n(X) \frac{z^n}{n!}.
		\end{equation*}
		
	\noindent
	Thus, the integral $M_X(z)$ acts as a generating function for the moments of $X$.  On the other hand, this integral
	may be explicitly evaluated.  Completing the square in the
	exponent of the integrand we find that
	
		\begin{equation*}
			M_X(z) = e^{\frac{z^2}{2}} \int\limits_{\field{R}} e^{-\frac{1}{2}(t-z)^2} \frac{\mathrm{d}t}{\sqrt{2\pi}},
		\end{equation*}
		
	\noindent
	whence
	
		\begin{equation*}
			M_X(z) = e^{\frac{z^2}{2}} = \sum_{k=0}^{\infty} \frac{z^{2k}}{2^k k!}
		\end{equation*}
	
	\noindent	
	by translation invariance of Lebesgue measure.
	We conclude that the odd moments of $X$ vanish while the even ones are given by the formula
	
		\begin{equation*}
			m_{2k}(X) = \frac{(2k)!}{2^k k!} = (2k-1) \cdot (2k-3) \cdot \dots \cdot 5 \cdot 3 \cdot 1.
		\end{equation*}
		
	\noindent
	This is the number of partitions of the set $[2k]$ into blocks of size two, also called ``pairings'': we have $2k-1$ choices for the element to be paired with $1$,
	then $2k-3$ choices for the element to be paired with the smallest remaining unpaired element, etc.
	Alternatively, we may say that $m_n(X)$ is equal to the number of $1$-regular graphs on $n$ labelled vertices.  It now follows from the 
	fundamental link between connected and disconnected structures that the
	cumulant $c_n(X)$ is equal to the number of connected $1$-regular graphs.  Consequently, the cumulant sequence
	of a standard Gaussian random variable is simply
	
		\begin{equation*}
			(0,1,0,0,0,\dots)
		\end{equation*}
		
	The fact that the universality of the Gaussian distribution is reflected in the simplicity of its cumulant sequence signals cumulants
	as a key concept in probability theory.  In Thiele's own words \cite{Hald2}, 
	
		\begin{center}
			\begin{quote}
			This remarkable proposition has originally led me to prefer the 
			half-invariants over every other system of symmetrical functions.
			\end{quote}
		\end{center}
			
	\noindent
	This sentiment persists amongst modern-day probabilists.  
	To quote Terry Speed \cite{Speed}, 
	
		\begin{center}
			\begin{quote}
			In a sense which it is hard to make precise, all of the important aspects of distributions seem to be simpler functions
			of cumulants than of anything else, and they are also the natural tools with which transformations of systems of random
			variables can be studied when exact distribution theory is out of the question.
			\end{quote}
		\end{center}
		
	\subsection{Cumulants and independence}
	The importance of cumulants stems, ultimately, from their relationship with stochastic independence.
	Suppose that $X$ and $Y$ are a pair of random variables whose moment sequences have been given to us
	by an oracle, and our task is to compute the moments of $X+Y$.  Since 
	$\E[X^aY^b]=\E[X^a]\E[Y^b]$, this can be done using the formula
	
		\begin{equation*}
			m_n(X+Y) = \sum_{k=0}^n {n \choose k} m_k(X)m_{n-k}(Y),
		\end{equation*}
		
	\noindent
	which is conceptually clear but computationally inefficient because of its dependence on $n$.
	For example, if we want to compute $m_{100}(X+Y)$ we must evaluate a sum with $101$ terms, each of 
	which is a product of three factors.  Computations with independent random variables 
	simplify dramatically if one works with cumulants rather than moments.
	Indeed, Thiele called cumulants ``half-invariants'' because
	
		\begin{equation*}
			\label{eqn:halfinvariant}
			X,Y \text{ independent } \implies c_n(X+Y) = c_n(X) + c_n(Y)\ \forall n \geq 1.
		\end{equation*} 
		
	\noindent
	Thanks to this formula,
	if the cumulant sequences of $X$ and $Y$ are given, then each cumulant of $X+Y$ can be computed simply by adding two
	numbers.  The mantra to be remembered is:
	
		\begin{center}
			\boxed{
			\emph{cumulants linearize addition of independent random variables}.}
		\end{center}
		
	\noindent
	For example, this fact together with the computation we did above yields 
	that the sum of two iid standard Gaussians is a Gaussian of variance two.
	
	In order to precisely understand the relationship between cumulants and independence, we need to 
	extend the relationship between moments and cumulants to a relationship between mixed moments and mixed
	cumulants.  Mixed moments are easy to define: given a set of 
	(not necessarily distinct) random variables $X_1,\dots,X_n$, 
	
		\begin{equation*}
			m_n(X_1,\dots,X_n) := \E[X_1 \dots X_n].
		\end{equation*}
		
	\noindent
	It is clear that $m_n(X_1,\dots,X_n)$ is a symmetric, multilinear function of its arguments.
	The new notation for mixed moments is related to our old notation for pure moments
	by
	
		\begin{equation*}
			m_n(X) = m_n(X,\dots,X),
		\end{equation*}
		
	\noindent
	which we may keep as a useful shorthand.
	
	We now define mixed cumulants recursively in terms of mixed moments using the natural
	extension of the moment-cumulant formula:
	
		\begin{equation*}
			m_n(X_1,\dots,X_n) = \sum_{\pi \in \lattice{P}(n)} \prod_{B \in \pi} c_{|B|}(X_i : i \in B).
		\end{equation*}
		
	\noindent
	For example, we have
	
		\begin{equation*}
			m_2(X_1,X_2) = c_2(X_1,X_2) + c_1(X_1)c_1(X_2),
		\end{equation*}
		
	\noindent
	from which we find that the second mixed cumulant of $X_1$ and $X_2$ is
	their covariance,
	
		\begin{equation*}
			c_2(X_1,X_2) = m_2(X_1,X_2)-m_1(X_1)m_2(X_2).
		\end{equation*}
		
	\noindent
	More generally, the recurrence
	
		\begin{equation*}
			c_n(X_1,\dots,X_n) = m_n(X_1,\dots,X_n) - \sum_{\substack{\pi \in \lattice{P}(n)\\ b(\pi) \geq 2}} 
			\prod_{B \in \pi} c_{|B|}(X_i : i \in B)
		\end{equation*}
		
	\noindent
	facilitates a straightforward inductive proof that $c_n(X_1,\dots,X_n)$ is a symmetric, $n$-linear function of
	its arguments, which explains Thiele's reference to cumulants as his preferred system of symmetric functions. 
	
	The fundamental relationship between cumulants and stochastic independence is the following:	
	$X$ and $Y$ are independent if and only if all their mixed cumulants vanish,
		
			\begin{equation*}
				\begin{split}
					&c_2(X,Y) = 0 \\
					&c_3(X,X,Y) = c_3(X,Y,Y) = 0 \\
					&c_4(X,X,X,Y) =c_4(X,X,Y,Y)=c_4(X,Y,Y,Y) = 0\\
					&\vdots
				\end{split}
			\end{equation*}	
	
	The forward direction of this theorem,
	
		\begin{equation*}	
				X,Y \text{ independent} \implies \text{ mixed cumulants vanish, }
		\end{equation*}
	
	\noindent
	immediately yields Thiele's linearization property, since by multilinearity we have
	
		\begin{equation*}
			\begin{split}
				c_n(X+Y) &= c_n(X+Y,\dots,X+Y) \\
				&= c_n(X,\dots,X) + \text{ mixed cumulants } + c_n(Y,\dots,Y) \\
				&= c_n(X)+c_n(Y).
			\end{split}
		\end{equation*}
		
\begin{figure}
			\includegraphics{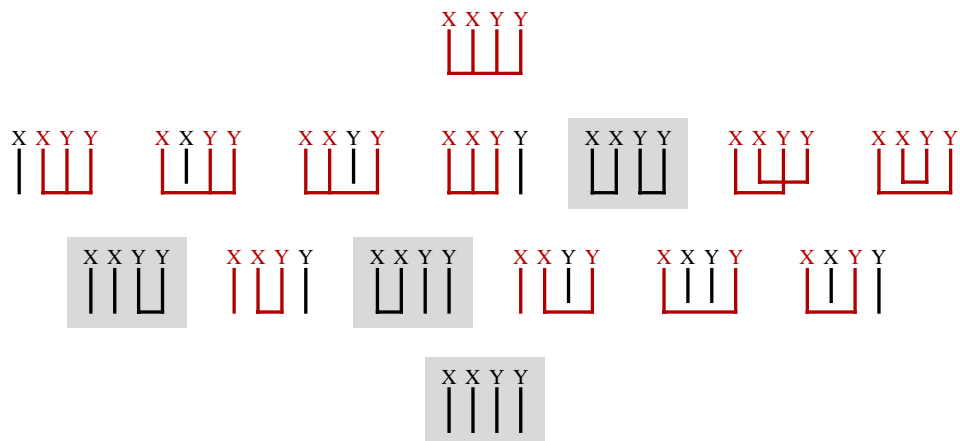}
			\caption{\label{fig:AllPartXXYY}Graphical evaluation of $m_4(X,X,Y,Y)$.}
		\end{figure}

	Conversely, let $X,Y$ be a pair of random variables whose mixed cumulants vanish.  
	Let us check in a couple of concrete cases that this condition forces $X$ and $Y$ to 
	obey the algebraic identities associated with independent random variables.
	In the first non-trivial case, $n=2$, vanishing of mixed cumulants reduces the 
	extended moment-cumulant formula to
	
		\begin{equation*}
			m_2(X,Y) = c_1(X)c_1(Y) = m_1(X)m_1(Y),
		\end{equation*}
		
	\noindent
	which is consistent with the factorization rule $\E[XY]=\E[X]\E[Y]$ for independent random variables.
	Now let us try an $n=4$ example.
	We compute $m_4(X,X,Y,Y)$ directly from the extended moment cumulant formula.  
	Referring to Figure \ref{fig:AllPartXXYY}, we find that vanishing of mixed cumulants implies
	
		\begin{equation*}
			\begin{split}
			m_4(X,X,Y,Y) &= c_2(X,X)c_2(Y,Y) + c_2(X,X)c_1(Y)c_1(Y) + c_2(Y,Y)c_1(X)c_1(X)\\
			&+ c_1(X)c_1(X)c_1(Y)c_1(Y),
			\end{split}
		\end{equation*}
		
	\noindent
	which reduces to the factorization identity $\E[X^2Y^2] = \E[X^2]\E[Y^2]$.  
	
\begin{figure}
			\includegraphics{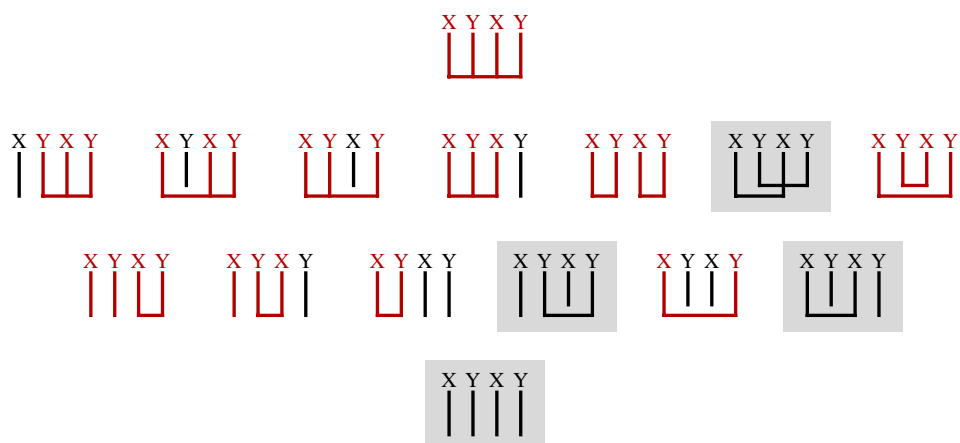}
			\caption{\label{fig:AllPartXYXY}Graphical evaluation of $m_4(X,Y,X,Y)$.}
		\end{figure}

	\noindent
	Of course, if we compute $m_4(X,Y,X,Y)$ using the extended moment-cumulant formula we should
	get the same answer, and indeed this is the case, but it is important to note that the contributions
	to the sum come from different partitions, as indicated in Figure \ref{fig:AllPartXYXY}.
	
	\subsection{Central Limit Theorem by cumulants}
	We can use the theory of cumulants presented thus far to prove an elementary 
	version of the Central Limit Theorem.		
	Let $X_1,X_2,X_3 \dots$ be a sequence of iid random variables, and let $X$ be a standard Gaussian.  
	Suppose that the common distribution of the variables $X_i$ has mean zero, variance one, and finite moments of all orders.  
	Put
	
		\begin{equation*}
			S_N := \frac{X_1 + \dots + X_N}{\sqrt{N}}.
		\end{equation*}
			
	\noindent
	Then, for each positive integer $n$,
		
		\begin{equation*}
			\lim_{N \rightarrow \infty} m_n(S_N) = m_n(X).
		\end{equation*}

	Since moments and cumulants mutually determine one another, in order to prove this CLT it suffices to prove that
		
		\begin{equation*}
			\lim_{N \rightarrow \infty} c_n(S_N) = c_n(X)
		\end{equation*}
			
	\noindent 
	for each $n \geq 1$.  Now, by multilinearity of $c_n$ and independence of the $X_i$'s, we have
		
		\begin{align*}
			c_n(S_N) &= c_n(N^{-\frac{1}{2}}(X_1+\dots+X_N)) \\
			&= N^{-\frac{n}{2}}(c_n(X_1) + \dots + c_n(X_N)) \\
			&= N^{1-\frac{n}{2}} c_n(X_1),
		\end{align*}
			
	\noindent
	where the last line follows from the fact that the $X_i$'s are equidistributed.
	Thus: if $n=1$,
			
		\begin{equation*}
			c_1(S_N) = N^{\frac{1}{2}} c_1(X_1) =0;
		\end{equation*}
			
	\noindent
	if $n=2$,
			
		\begin{equation*}
			c_2(S_N) = c_2(X_1) = 1;
		\end{equation*}
			
	\noindent
	if $n>2$,
		
		\begin{equation*}
			c_n(S_N) = N^{\text{negative number}}c_n(X_1).
		\end{equation*}
	
	\noindent		
	We conclude that
		
		\begin{equation*}
			\lim_{N \rightarrow \infty} c_n(S_N) = \delta_{n2},
		\end{equation*}
			
	\noindent
	which we have already identified as the cumulant sequence of a 
	standard Gaussian random variable.
			
	\subsection{Geometrically connected graphs}
	Let us now consider a variation on our original graph-counting question.  Given a 
	graph $G$ on the vertex set $[n]$, we may represent its vertices by $n$ distinct points
	on the unit circle (say, the $n^{\text{th}}$ roots of unity) and its edges by straight line segments joining these points.
	This is how we represented the set of four-vertex graphs in Figure \ref{fig:4-vertex}.
	We will denote this geometric realization of $G$ by $|G|$.  The geometric realization of a graph carries extra structure
	which we may wish to consider.  For example, it may happen that $|G|$ is a connected set of points in the plane 
	even if the graph $G$ is not connected in the usual sense of graph theory.  Let $\kappa_n$ denote the number of 
	geometrically connected graphs on $[n]$.  This is sequence $\mathsf{A136653}$ in Sloane's database; its first
	few terms are
	
		\begin{equation*}
			1,\ 1,\ 4,\ 39,\ 748,\ 27\, 162,\ 1\, 880\, 872,\ \dots.
		\end{equation*}
		
\begin{figure}
			\includegraphics{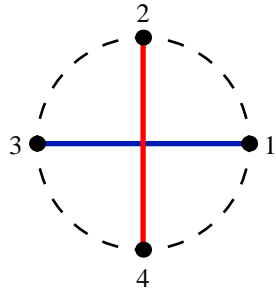}
			\caption{\label{fig:Crosshairs}The crosshairs graph.}
		\end{figure}

	\noindent
	Since geometric connectivity is a weaker condition than set-theoretic connectivity, $\kappa_n$ grows faster
	than $c_n$; these sequences diverge from one another at $n=4$, where the unique disconnected but geometrically
	connected graph is the ``crosshairs'' graph shown in Figure \ref{fig:Crosshairs}.
	
	Consider now the problem of computing $\kappa_n$.  As with $c_n$, we can address this problem by 
	means of a combinatorial decomposition of the set of graphs with $n$ vertices.  However, this decomposition
	must take into account the planar nature of geometric connectivity, 
	which our previous set-theoretic decompositions do not.  Consequently, we must
	formulate a new decomposition.

\begin{figure}
			\includegraphics{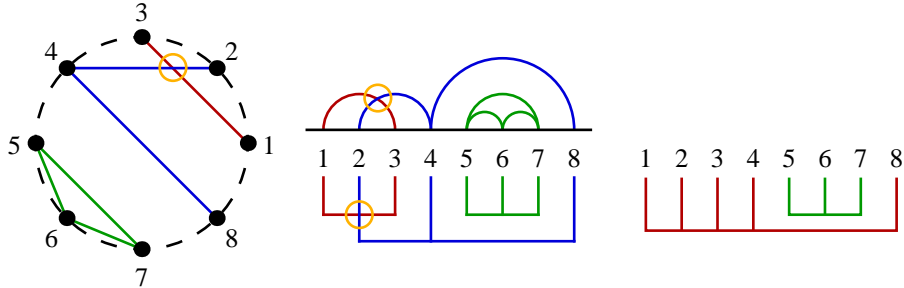}
			\caption{\label{fig:PartitionFusion}Partition fusion accounts for geometric connectedness.}
		\end{figure}

	Given a graph $G$ on $[n]$, let $\pi(G)$ denote the partition of $[n]$ induced by the connected components
	of $G$ ($i$ and $j$ are in the same block of $\pi(G)$ if and only if they are in the same connected component of $G$), 
	and let $\pi(|G|)$ denote the partition of $[n]$ induced by the geometrically connected components
	of $|G|$ ($i$ and $j$ are in the same block of $\pi(|G|)$ if and only if they are in the same geometrically connected 
	component of $|G|$).  How are $\pi(G)$ and $\pi(|G|)$ related?  To understand this, let us view our geometric graph
	realizations as living in the hyperbolic plane rather than the Euclidean plane.  Thus
	Figure \ref{fig:4-vertex} depicts line systems in the Klein model, in which the plane is an open disc and straight lines are chords of the boundary circle.
	We could alternatively represent a graph in the Poincar\'e disc model, where straight lines are arcs of circles orthogonal
	to the boundary circle, or in the Poincar\'e half-plane model, where space is an open-half plane and straight lines are
	arcs of circles orthogonal to the boundary line.  The notion of geometric connectedness does not depend on the particular
	realization chosen.  The half-plane model has the useful feature that the geometric realization $|G|$ essentially coincides
	with the pictorial representation of $\pi(G)$, and we can see clearly that crossings in $|G|$ correspond exactly
	to crossings in $\pi(G)$.  Thus, $\pi(|G|)$ is obtained by fusing together crossing blocks of $\pi(G)$.  
	The resulting partition $\pi(|G|)$ no longer has any crossings --- by construction, it is a non-crossing partition, see
	figure \ref{fig:PartitionFusion}.
	
	We can now obtain a recurrence for $\kappa_n$.  We construct a graph by 
	first choosing a non-crossing partition of the underlying vertex set into blocks $B_1,\dots,B_k$ and
	then building a geometrically connected graph on each block, which can be done in 
	$\kappa_{|B_1|} \dots \kappa_{|B_k|}$ ways.  This leads to the formula

		\begin{equation*}
			m_n = \sum_{\pi \in \lattice{NC}(n)} \prod_{B \in \pi} \kappa_{|B|},
		\end{equation*}
		
	\noindent
	where the summation is over non-crossing partitions of $[n]$.
	Just as before, we can split off the term of the sum corresponding to the partition with only one
	block to obtain the recursion
	
		\begin{equation*}
			\kappa_n =m_n- \sum_{\substack{\pi \in \lattice{NC}(n)\\ b(\pi) \geq 2}} \prod_{B \in \pi} \kappa_{|B|},
		\end{equation*}
		
	\noindent
	in which we sum over non-crossing partitions with at least two blocks.

	\subsection{Non-crossing cumulants}
	We have seen above that the usual set-theoretic notion of connectedness manifests itself probabilistically as the 
	cumulant concept.  We have also seen that set-theoretic connectedness has an interesting geometric variation,
	which we called geometric connectedness.  This begs the question: 
	
		\begin{center}
			Is there a probabilistic interpretation of geometric connectedness?
		\end{center}
				
	Let $X$ be a random variable, with moments $m_n(X)$.  Just as the classical cumulants $c_n(X)$ were
	defined recursively using the relation between all structures and connected structures, we define the 
	non-crossing cumulants of $X$ recursively using the relation between all structures 
	and geometrically connected structures:
	
		\begin{equation*}
			m_n(X) = \sum_{\lattice{NC}(n)} \prod_{B \in \pi} \kappa_{|B|}(X).
		\end{equation*}
	
	\noindent
	We will call this the non-crossing moment-cumulant formula.
	Since connectedness and geometric connectedness coincide for structures of size $n=1,2,3$, the first
	three non-crossing cumulants of $X$ are identical to its first three classical cumulants.  However, for $n \geq 4$,
	the non-crossing cumulants become genuinely new statistics of $X$.  
	
	Our first step in investigating these new statistics is to look for a non-crossing analogue of the most
	important random variable, the standard Gaussian.  This should be a random variable whose
	non-crossing cumulant sequence is 
	
		\begin{equation*}
			0,\ 1,\ 0,\ 0,\ \dots.
		\end{equation*}
		
	\noindent
	If this search leads to something interesting, we may be motivated to further investigate non-crossing probability theory.
	If not, we will reject the idea as a will-o'-the-wisp.
	
	From the non-crossing moment-cumulant formula, we find that the moments of the non-crossing Gaussian $X$ are
	given by
	
		\begin{equation*}
			m_n(X) = \sum_{\pi \in \lattice{NC}(n)} \prod_{B \in \pi} \delta_{|B|,2} = \sum_{\pi \in \lattice{NC}_2(n)} 1.
		\end{equation*}
		
	\noindent
	That is, $m_n(X)$ is equal to the number of partitions in $\lattice{NC}(n)$ all of whose blocks have size
	$2$, i.e. non-crossing pairings of $n$ points.  We know that there are no pairings at all on an odd number of 
	points, so the odd moments of $X$ must be zero, which indicates that $X$ likely has a symmetric distribution.  
	The number of pairings on $n=2k$ points is given by a factorial going down in steps of two, $(2k-1)!! = (2k-1) \cdot (2k-3) 
	\dots \cdot 5 \cdot 3 \cdot 1$, so the number of non-crossing pairings must be smaller than this double factorial.
	
\begin{figure}
			\includegraphics{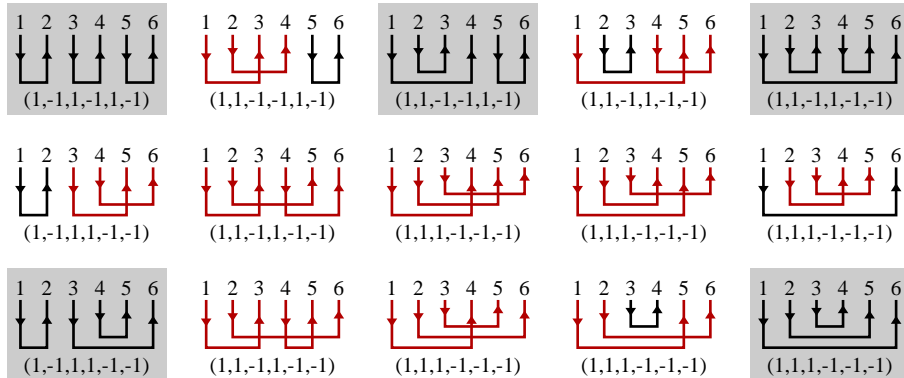}
			\caption{\label{fig:Bitstrings}Construction of the function $f$ from pairings to bitstrings.}
		\end{figure}

	In order to count non-crossing pairings on $2k$ points, we construct a function $f$ from the set 
	of all pairings on $2k$ points to length $2k$ sequences of $\pm1$'s.  This function is easy to describe:
	if $i<j$ constitute a block of $\pi$, then the $i^{\text{th}}$ element of $f(\pi)$ is $+1$ and the $j^{\text{th}}$
	element of $f(\pi)$ is $-1$.  See Figure \ref{fig:Bitstrings} for an illustration of this function in the case 
	$k=3$.  By construction, $f$ is a surjection from the set of pairings on $2k$ points onto
	the set of length $2k$ sequences of $\pm1$'s all of whose partial sums are non-negative and whose total
	sum is zero.  We leave it to the reader to show that the fibre of $f$ over any such sequence contains exactly
	one non-crossing pairing, so that $f$ restricts to a bijection from non-crossing pairings onto its image.  
	The image sequences can be neatly enumerated using the Dvoretzky-Motzkin-Raney cyclic shift lemma, as
	in \cite[\S7.5]{GKP}.  They are counted by the Catalan numbers

		\begin{equation*}
			\Cat_k = \frac{1}{k+1} {2k \choose k},
		\end{equation*}

	\noindent
	which are smaller than the double factorials by a factor of $2^k/(k+1)!$.  This indicates
	that the distribution of $X$ decays even more rapidly than the Gaussian distribution and might
	even be compactly supported.
	
	We have discovered that
	
		\begin{equation*}
			m_n(X) = \begin{cases}
					0, \text{ if $n$ odd} \\
					\Cat_{\frac{n}{2}}, \text{ if $n$ even}.
				\end{cases}
		\end{equation*}
		
	\noindent
	The Catalan numbers are ubiquitous in enumerative combinatorics, see
	\cite[Exercise 6.19]{Stanley:EC2} as well as \cite{Stanley:CA}, 
	and their appearance in this context is the first sign that we are onto something interesting.
	We are now faced with an inverse problem: we are not trying to calculate the moments of a random variable 
	given its distribution, rather we know that the moment sequence of $X$ is
	
		\begin{equation*}
			0,\ \Cat_1,\ 0,\ \Cat_2,\ 0,\ \Cat_3,\ 0,\ \dots.
		\end{equation*}
		
	\noindent
	and we would like to write down its distribution $\mu_X$.
	Equivalently, we are looking for an integral representation of the entire function
	
		\begin{equation*}
			M_X(z) = \sum_{n=0}^{\infty} \Cat_n \frac{z^{2n}}{(2n)!}= \sum_{n=0}^{\infty} \frac{z^{2n}}{n!(n+1)!}
		\end{equation*}
		
	\noindent
	which has the form
	
		\begin{equation*}
			M_X(z) = \int\limits_{\field{R}} e^{tz} \mu_X(\mathrm{d}t),
		\end{equation*}
		
	\noindent
	with $\mu_X$ a probability measure on the real line.  
	The solution to this problem can 
	be extracted from the classical theory of Bessel functions.  
	
	The modified Bessel function $I_\alpha(z)$
	of order $\alpha$ is one of two linearly independent solutions to the modified Bessel equation
	
		\begin{equation*}
			\bigg{(}z^2\frac{d^2}{dz^2} + z\frac{d}{dz} - (z^2+\alpha^2)\bigg{)}F =0,
		\end{equation*}
		
	\noindent
	the other being the Macdonald function 
	
		\begin{equation*}
			K_\alpha(z) = \frac{\pi}{2} \frac{I_{-\alpha}(z)-I_{\alpha}(z)}{\sin(\alpha\pi)}.
		\end{equation*}
		
	\noindent
	The modified Bessel equation (and hence the functions $I_\alpha,K_\alpha$) appears in many problems of physics and engineering
	since it is related to solutions of Laplace's equation with cylindrical symmetry.  
	An excellent reference on this topic
	is \cite[Chapter 4]{AAR}.
	
	Interestingly, Bessel functions also occur in the combinatorics of permutations: a remarkable
	identity due to Ira Gessel asserts that 
	
		\begin{equation*}
			\det [I_{i-j}(2z)]_{i,j=1}^k = \sum_{n=0}^{\infty} \lis_k(n) \frac{z^{2n}}{(n!)^2},
		\end{equation*}
		
	\noindent
	where $\lis_k(n)$ is the number of permutations in the symmetric group $\group{S}(n)$ with
	no increasing subsequence of length $k+1$.  Gessel's identity was the point of departure in the 
	work of Jinho Baik, Percy Deift and Kurt Johansson who, answering a question posed by
	Stanislaw Ulam, proved that the limit distribution of the length of the longest increasing subsequence
	in a uniformly distributed random permutation is given by the ($\beta=2$) Tracy-Widom distribution.
	This non-classical distribution was isolated and studied by Craig Tracy and Harold Widom in a series
	of works on random matrix theory in the early 1990's where it emerged as the limiting distribution
	of the top eigenvalue of large random Hermitian matrices.  
	It has a density which may also be described in terms
	of Bessel functions, albeit indirectly.  Consider the ordinary differential equation
	
		\begin{equation*}
			\frac{d^2}{dx^2}u = 2u^3+xu
		\end{equation*}
		
	\noindent
	for a real function $u=u(x)$, which is known as the Painlev\'e II equation after the French mathematician
	(and two-time Prime Minister of France) Paul Painlev\'e.
	It is known that this equation has a unique solution, called the 
	Hastings-McLeod solution, with the asymptotics $u(x) \sim -\operatorname{Ai}(x)$ as
	$x \rightarrow \infty$, where 
	
		\begin{equation*}
			\operatorname{Ai}(x) = \frac{1}{\pi} \sqrt{\frac{x}{3}}K_{\frac{1}{3}}(\frac{2}{3}x^{\frac{3}{2}})
		\end{equation*}
		
	\noindent
	is a scaled specialization of the Macdonald function known as the Airy function.  Define the Tracy-Widom
	distribution function by
	
		\begin{equation*}
			F(t) = e^{-\int_t^{\infty} (x-t)u(x)^2 \mathrm{d}x},
		\end{equation*}
		
	\noindent
	where $u$ is the Hastings-McLeod solution to Painlev\'e II.  The theorem of Baik, Deift and Johansson asserts
	that 
	
		\begin{equation*}
			\lim_{n \rightarrow \infty} \frac{1}{n!}\lis_{2\sqrt{n}+tn^{1/6}}(n) = F(t)
		\end{equation*}
	
	\noindent
	for any $t \in \field{R}$.  From this one may conclude, for example, that the probability a permutation 
	drawn uniformly at random from the symmetric group $\group{S}(n^2)$ avoids the pattern 
	$1\ 2\ \dots\ 2n+1$ converges to $F(0) = 0.9694\dots$.
	We refer the interested reader to Richard Stanley's survey \cite{Stanley:ICM}
	for more information on this topic.
		
	Nineteenth century mathematicians knew how to describe the modified Bessel function both as a series,
	
		\begin{equation*}
			I_\alpha(z) = \sum_{n=0}^{\infty} \frac{(\frac{z}{2})^{2n+\alpha}}{n!\Gamma(n+1+\alpha)},
		\end{equation*}
		
	\noindent
	and as an integral,
	
		\begin{equation*}
			I_\alpha(z) = \frac{(\frac{z}{2})^{\alpha}}{\sqrt{\pi}\Gamma(\alpha+\frac{1}{2})} \int\limits_0^{\pi} 
			e^{(\cos \theta)z} (\sin \theta)^{2\alpha} \mathrm{d}\theta.
		\end{equation*}
		
	\noindent
	From the series representation we find that
	
		\begin{equation*}
			M_X(z) = \frac{I_1(2z)}{z},
		\end{equation*}
		
	\noindent
	and consequently we have the integral representation
	
		\begin{equation*}
			M_X(z) = \frac{2}{\pi} \int\limits_0^{\pi} e^{2(\cos\theta)z} \sin^2 \theta \mathrm{d}\theta.
		\end{equation*}
		
\begin{figure}
			\includegraphics{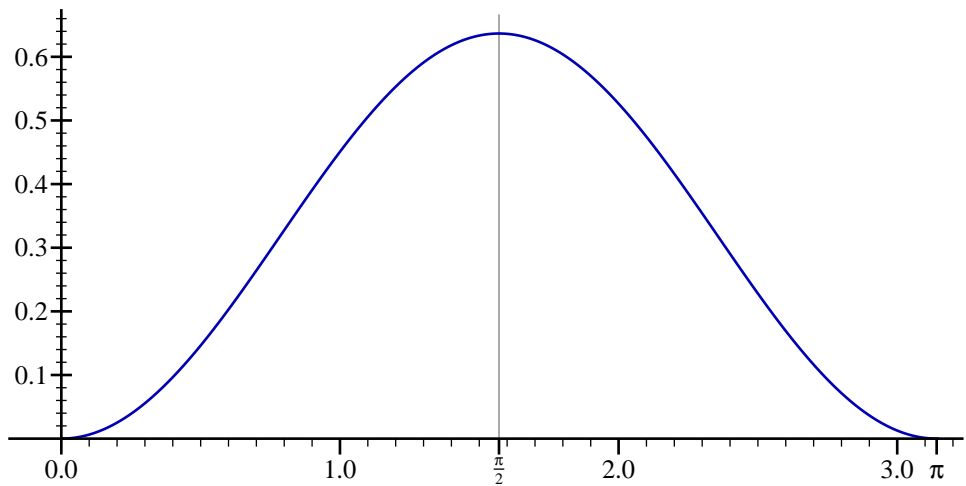}
			\caption{\label{fig:SatoTate}The Sato-Tate density}
		\end{figure}

	\noindent
	This is one step removed from what we want: it tells us that the Catalan numbers are the even 
	moments of the random variable $X=2\cos(Y)$, where $Y$ is a random variable with distribution 
	
		\begin{equation*}
			\mu_Y(\mathrm{d}\theta) = \frac{2}{\pi} \sin^2\theta \mathrm{d}\theta
		\end{equation*}
		
	\noindent
	supported on the interval $[0,\pi]$.  However, this is a rather interesting intermediate step
	since the above measure appears in number theory, where it is called the Sato-Tate distribution, 
	see Figure \ref{fig:SatoTate}.  
	
\begin{figure}
			\includegraphics{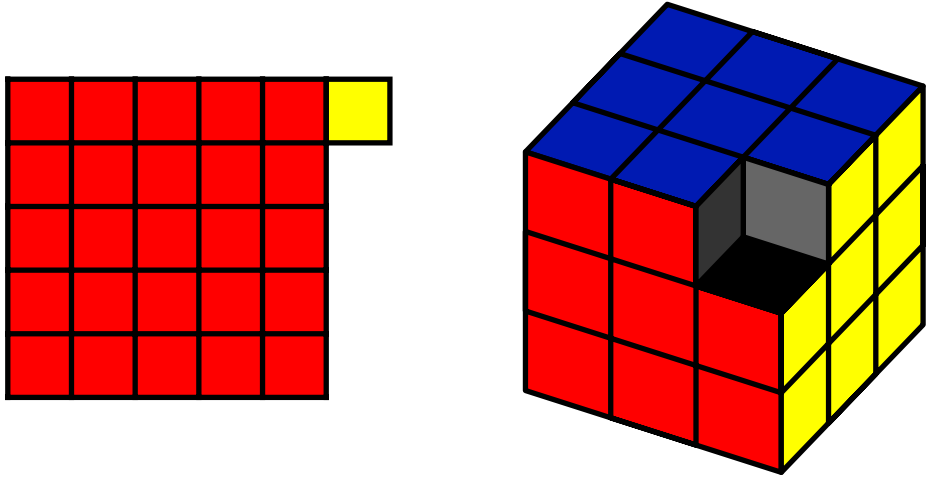}
			\caption{\label{fig:Fermat}Diophantine perspectives on twenty-six}
		\end{figure}

	The Sato-Tate distribution arises in the arithmetic statistics of elliptic curves.  
	The location of integer points on elliptic curves	
	is a classical topic in number theory.  For example, Diophantus of Alexandria 
	wrote that the equation
	
		\begin{equation*}
			y^2=x^3-2
		\end{equation*}
		
	\noindent
	has the solution $x=3,y=5$, and in the 1650's Pierre de Fermat claimed that there are no other positive
	integer solutions.  This is the striking assertion that $26$ is the only number 
	one greater than a perfect square and one less than a perfect cube, see Figure \ref{fig:Fermat}.  
	That this is indeed the case
	was proved by Leonhard Euler in 1770, although according to some sources Euler's proof was incomplete and the
	solution to this problem should be attributed to Axel Thue in 1908.  	
	
	Modern number theorists study solutions to elliptic Diophantine 
	equations by reducing modulo primes.
	Given an elliptic curve
			
		\begin{equation*}
			y^2=x^3+ax+b, \quad a,b \in \ring{Z},
		\end{equation*}
		
	\noindent
	let $\Delta=-16(4a^3+27b^2)$ be sixteen times the 
	discriminant of $x^3+ax+b$, and let $S_p$ be the number of solutions of the 
	congruence
	
		\begin{equation*}
			y^2 \equiv x^3+ax+b \quad \mod p
		\end{equation*}
		
	\noindent
	where $p$ is a prime which does not divide $\Delta$.
	In his 1924 
	doctoral thesis, Emil Artin conjectured that 
	
		\begin{equation*}
			|S_p - p| \leq 2\sqrt{p}
		\end{equation*}
		
	\noindent
	for all such good reduction primes.  This remarkable inequality states that the number of solutions 
	modulo $p$ is roughly $p$ itself, up to an error of order $\sqrt{p}$.
	Artin's conjecture was proved by Helmut Hasse in 1933.  Around 1960, Mikio Sato and 
	John Tate became interested in the finer question of the distribution of the centred and scaled
	solution count $(S_p-p)/\sqrt{p}$ for typical elliptic 
	curves $E$ (meaning those without complex multiplication) as $p$ ranges over the infinitely many 
	primes not dividing the discriminant of $E$.  Because of Hasse's theorem, this 
	amounts to studying the distribution of the angle $\theta_p$ defined by
		
		\begin{equation*}
			\frac{S_p-p}{\sqrt{p}} = 2\cos\theta_p
		\end{equation*}
		
	\noindent
	in the the interval $[0,\pi]$.  Define a sequence $\mu_N^E$ of empirical probability measures 
	associated to $E$ by
	
		\begin{equation*}
			\mu_N^E = \frac{1}{\pi(N)} \sum_{p \leq N} \delta_{\theta_p},
		\end{equation*}
		
	\noindent
	where $\pi(N)$ is the number of prime numbers less than or equal to $N$.
	Sato and Tate conjectured that, for any elliptic curve $E$ without complex multiplication,
	$\mu_N^E$ converges weakly to the Sato-Tate distribution as $N \rightarrow \infty$.  
	This is a universality conjecture: it posits that certain limiting behaviour
	is common to a large class of elliptic curves irrespective of their fine structural details.
	Major progress on the Sato-Tate
	conjecture has been made within the last decade; we refer the reader to the surveys of Barry Mazur \cite{Mazur} and
	Ram Murty and Kumar Murty \cite{MM} for further information.
	
	The random variable we seek is not the Sato-Tate variable $Y$, but twice its cosine, $X=2\cos(Y)$.  Making the substitution 
	$s=\arccos(\theta)$ in the integral representation of $M_X(z)$ obtained above, we obtain
	
		\begin{equation*}
			M_X(z) = \frac{2}{\pi}\int\limits_{-1}^1 e^{2sz} \sqrt{1-s^2} \mathrm{d}s,
		\end{equation*}
		
	\noindent
	and further substituting $t=2s$ this becomes
	
		\begin{equation*}
			M_X(z) = \frac{1}{2\pi} \int\limits_{-2}^2 e^{tz} \sqrt{4-t^2}\mathrm{d}t.
		\end{equation*}
		
\begin{figure}
			\includegraphics{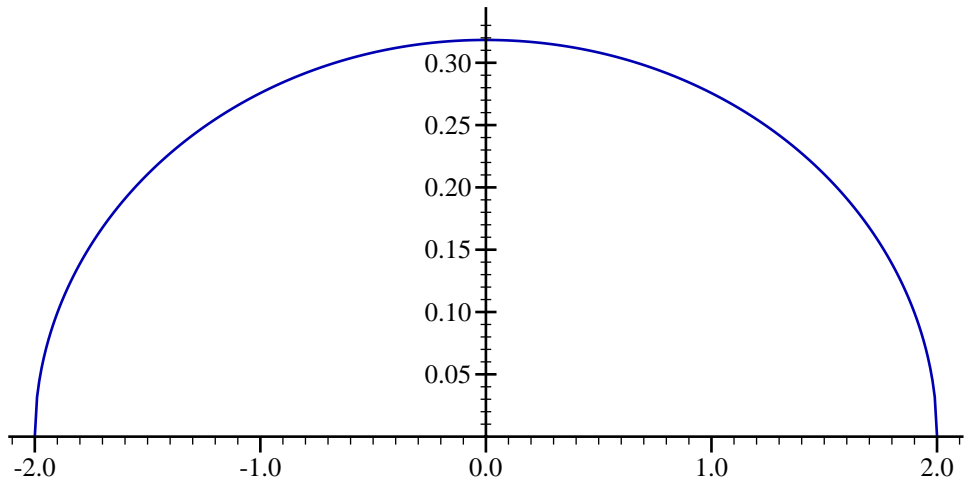}
			\caption{\label{fig:Wigner}The Wigner semicircle density}
		\end{figure}

	\noindent
	Thus the random variable $X$ with even moments the Catalan numbers and vanishing odd moments
	is distributed in the interval $[-2,2]$ with density
	
		\begin{equation*}
			\mu_X(\mathrm{d}t) = \frac{1}{2\pi} \sqrt{4-t^2}\mathrm{d}t,
		\end{equation*}
		
	\noindent
	which is both symmetric and compactly supported.
	This is another famous distribution: it is called the Wigner semicircle distribution
	after the physicist Eugene Wigner, who considered it in the 1950's in a context ostensibly
	unrelated to elliptic curves.  The density 
	of $\mu_X$ is shown in Figure \ref{fig:Wigner} --- note that it is not a semicircle, but rather half an ellipse of 
	semi-major axis two and semi-minor axis $1/\pi$.
	
	Wigner was interested in constructing models for the energy levels of complex systems, and hit on the 
	idea that the eigenvalues of large symmetric random matrices provide a good approximation.  Wigner considered
	$N \times N$ symmetric matrices $X_N$ whose entries $X_N(ij)$ are independent random variables, up to 
	the symmetry constraint $X_N(ij)=X_N(ji)$.  
	Random matrices of this form are now known as Wigner matrices, and their study remains a topic 
	of major interest today.  Wigner studied the empirical spectral distribution of the eigenvalues of $X_N$,
	i.e. the probability measure
		
		\begin{equation*}
			\mu_N = \frac{1}{N} \sum_{k=1}^N \delta_{\lambda_k(N)}
		\end{equation*}
		
	\noindent
	which places mass $1/N$ at each eigenvalue of $X_N$.  Note that, unlike in the setting above where we considered
	the sequence of empirical measures associated to a fixed elliptic curve $E$, the measure $\mu_N$ is a random measure
	since $X_N$ is a random matrix.
	Wigner showed that the limiting behaviour of $\mu_N$ does not depend on the details of the random variables
	which make up $X_N$.  In \cite{Wigner}, he made the following hypotheses:
	
		\begin{enumerate}
			
			\item
			Each $X_N(ij)$ has a symmetric distribution;
			
			\item
			Each $X_N(ij)$ has finite moments of all orders, each of which is 
			bounded by a constant independent of $N,i,j$;
			
			\item
			The variance of $X_N(ij)$ is $1/N$.
						
		\end{enumerate}
		
	\noindent 
	Wigner prove that, under these hypotheses, $\mu_N$ converges weakly to the semicircle law which now bears 
	his name.  We will see a proof of Wigner's theorem for random matrices with (complex) Gaussian entries in 
	Lecture Three.  The universality of the spectral structure of real and complex Wigner matrices holds at a much finer level, 
	and under much weaker hypotheses,
	both at the edges of the semicircle \cite{Soshnikov} and in the bulk \cite{ESY,TV}.
				
	\subsection{Non-crossing independence}
	Our quest for the non-crossing Gaussian has brought us into contact with interesting objects (random permutations,
	elliptic curves, random matrices)
	and the limit laws which govern them (Tracy-Widom distribution, Sato-Tate distribution, Wigner semicircle distribution).  This motivates us 
	to continue developing the rudiments of non-crossing probability theory --- perhaps we have hit on a framework
	within which these objects may be studied.
	
	Our next step is to introduce a notion of non-crossing independence. 
	We know that classical independence is characterized by the vanishing of mixed cumulants.  
	Imitating this, we will define non-crossing independence via the vanishing of mixed non-crossing cumulants.
	Like classical mixed cumulants,
	the non-crossing mixed cumulant functionals are defined recursively via the 
	multilinear extension of the non-crossing moment-cumulant formula,
			
		\begin{equation*}
			m_n(X_1,\dots,X_n) = \sum_{\pi \in \lattice{NC}(n)} \prod_{B \in \pi} \kappa_{|B|}(X_i : i \in B).
		\end{equation*}
				
	\noindent
	The recurrence 
	
		\begin{equation*}
			\kappa_n(X_1,\dots,X_n) = m_n(X_1,\dots,X_n)-\sum_{\pi \in \lattice{NC}(n)} \prod_{B \in \pi} \kappa_{|B|}(X_i : i \in B)
		\end{equation*}
		
	\noindent
	and induction establish that $\kappa_n(X_1,\dots,X_n)$ is a symmetric multilinear function of its arguments.
	Two random variables $X,Y$ are said to be non-crossing independent if their mixed 
	non-crossing cumulants vanish:
	
		\begin{equation*}
				\begin{split}
					&\kappa_2(X,Y) = 0 \\
					&\kappa_3(X,X,Y) = \kappa_3(X,Y,Y) = 0 \\
					&\kappa_4(X,X,X,Y) =\kappa_4(X,X,Y,Y)=\kappa_4(X,Y,Y,Y) = 0\\
					&\vdots
				\end{split}
			\end{equation*}	
		
	An almost tautological consequence of this definition is: 
	
		\begin{equation*}
			X,Y \text{ non-crossing independent } \implies \kappa_n(X+Y) = \kappa_n(X) + \kappa_n(Y)\ \forall n \geq 1.
		\end{equation*} 
		
	\noindent
	Thus, just as classical cumulants linearize the addition of classically independent random variables,
	
		\begin{center}
			\boxed{
			\emph{non-crossing cumulants linearize addition of non-crossing independent random variables}.}
		\end{center}
		
	We can also note that the semicircular random variable $X$, whose non-crossing cumulant sequence is 
	$0,1,0,0,\dots$, plays the role of the standard Gaussian with respect to this new notion of independence.
	For example, since non-crossing cumulants linearize non-crossing independence, the sum of two 
	non-crossing independent semicircular random variables is a semicircular random variable of variance two.
	The non-crossing analogue of the Central Limit Theorem asserts that, if $X_1,X_2,\dots$ is a sequence of non-crossing
	independent and identically distributed random variables with mean zero and variance one, then the moments of
	
		\begin{equation*}
			S_N = \frac{X_1+\dots+X_N}{\sqrt{N}}
		\end{equation*}
		
	\noindent
	converge to the moments of the standard semicircular $X$ as $N \rightarrow \infty$.  The proof of this fact is identical
	to the proof of the classical Central Limit Theorem given above, except that classical cumulants are replaced
	by non-crossing cumulants.
	
	Of course, we don't really know what non-crossing independence means.  For example, if $X$ and $Y$ are non-crossing 
	independent, is it true that $\E[XY]=\E[X]\E[Y]$?  The answer is yes, since classical and non-crossing mixed cumulants
	agree up to and including order three,  
	
		\begin{equation*}
			c_1(X)=\kappa_1(X), \quad c_2(X,Y)=\kappa_2(X,Y), \quad c_3(X,Y,Z)=\kappa_3(X,Y,Z).
		\end{equation*}
	
	\noindent
	But what about higher order mixed moments?
		
\begin{figure}
			\includegraphics{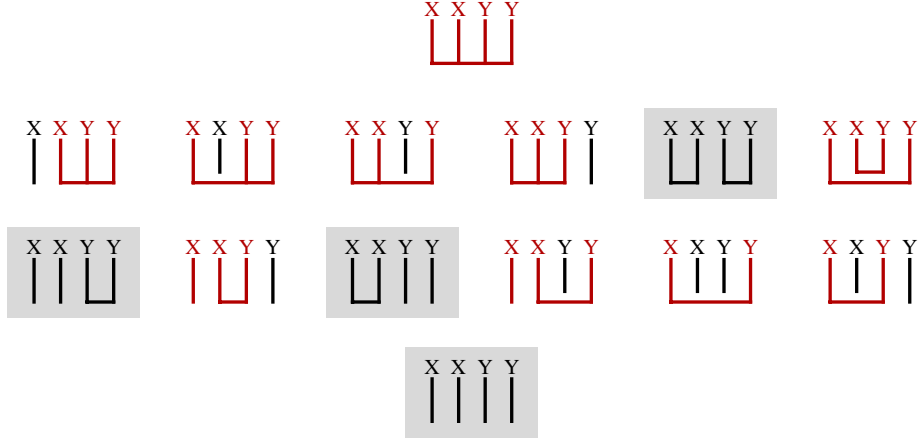}
			\caption{\label{fig:NoncrossingXXYY}Graphical evaluation of $m_4(X,X,Y,Y)$ using non-crossing cumulants.}
		\end{figure}

	We observed above that, in the classical case, vanishing of mixed cumulants allows us to recover the familiar
	algebraic identities governing the expectation of independent random variables.
	We do not have a priori knowledge of the algebraic identities governing the expectation of non-crossing independent
	random variables, so we must discover them using the 
	vanishing of mixed non-crossing cumulants.  Let us see what this implies
	for the mixed moment $m_4(X,X,Y,Y)=\E[X^2Y^2]$.  Referring to Figure \ref{fig:NoncrossingXXYY} we see that in this case
	the non-crossing moment-cumulant formula reduces to 
	
		\begin{equation*}
			\begin{split}
			m_4(X,X,Y,Y) &= \kappa_2(X,X)\kappa_2(Y,Y) + \kappa_2(X,X)\kappa_1(Y)\kappa_1(Y) + \kappa_2(Y,Y)\kappa_1(X)\kappa_1(X) \\
				&+\kappa_1(X)\kappa_1(X)\kappa_1(Y)\kappa_1(Y)
			\end{split}
		\end{equation*}
		
	\noindent
	which is exactly the formula we obtained for classically independent random variablesusing the classical moment-cumulant formula.
		
\begin{figure}
			\includegraphics{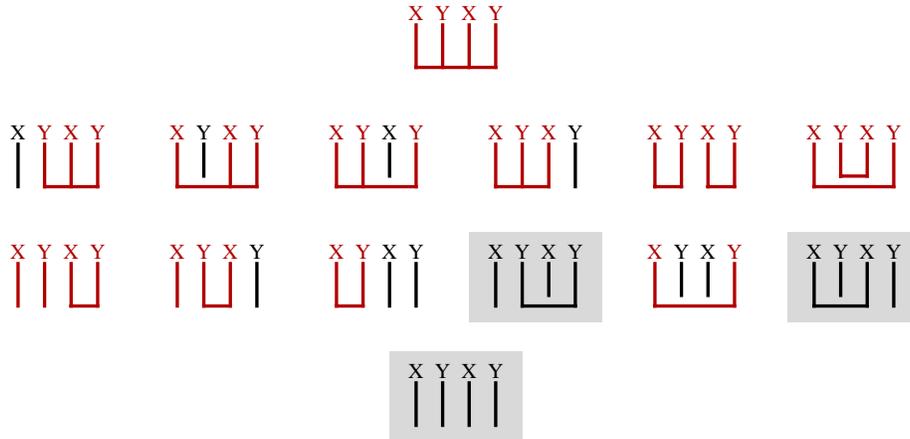}
			\caption{\label{fig:NoncrossingXYXY}Graphical evaluation of $m_4(X,Y,X,Y)$ using non-crossing cumulants.}
		\end{figure}

	\noindent
	However, when we use the non-crossing moment-cumulant formula to evaluate the same mixed moment
	with its arguments permuted, we instead get
	
		\begin{equation*}
			m_4(X,Y,X,Y) = \kappa_2(X,X)\kappa_1(Y)\kappa_1(Y) + \kappa_2(Y,Y)\kappa_1(X)\kappa_1(X) 
			+\kappa_1(X)\kappa_1(X)\kappa_1(Y)\kappa_1(Y),
		\end{equation*}
		
	\noindent
	see Figure \ref{fig:NoncrossingXYXY}.
	
	\noindent
	Since $m_4(X,X,Y,Y)=m_4(X,Y,X,Y)$, we are forced to conclude that the two expressions
	obtained are equal, which in turn forces
	
		\begin{equation*}
			\kappa_2(X,X)\kappa_2(Y,Y)=0.
		\end{equation*}
		
	\noindent
	Thus, if $X,Y$ are non-crossing independent random variables, at least one of them must have vanishing variance,
	and consequently must be almost surely constant.
	The converse is also true --- one can show that a (classical or non-crossing) mixed cumulant
	vanishes if any of its entries are constant random variables.  So we have classified pairs of non-crossing independent
	random variables: they look like $\{X,Y\}= \{\text{arbitrary},\text{constant}\}$.  Such pairs of random variables are
	of no interest from a probabilistic perspective.  It would seem that non-crossing probability is a dead end.
		
	\subsection{The medium is the message}
	If $\Omega$ is a compact Hausdorff space then the algebra $\algebra{A}(\Omega)$ of 
	continuous functions $X:\Omega \rightarrow \field{C}$
	is a commutative $C^*$-algebra.  This means that in addition to its standard algebraic structure (pointwise addition, multiplication
	and scalar multiplication of functions) $\algebra{A}(\Omega)$ is equipped with a norm satisfying the Banach algebra
	axioms and an antilinear involution which is compatible with the norm, $\|X^*X\| = \|X\|^2$.		
	The norm comes from the topology of the source, $\|X\| = \sup_\omega |X(\omega)|$, and the
	involution comes from the conjugation automorphism of the target, $X^*(\omega) = \overline{X(\omega)}$.
	Conversely, a famous theorem of Israel Gelfand asserts that any unital commutative $C^*$-algebra $\algebra{A}$
	can be realized as the algebra of continuous functions on a compact Hausdorff space $\Omega(\algebra{A})$ 
	in an essentially unique way.  In fact, $\Omega(\algebra{A})$ may be constructed as the set of maximal ideals
	of $\algebra{A}$ equipped with a suitable topology.  The associations $\Omega \mapsto \algebra{A}(\Omega)$
	and $\algebra{A} \mapsto \Omega(\algebra{A})$ are contravariantly functorial and set up a dual equivalence
	between the category of compact Hausdorff spaces and the category of 
	unital commutative $C^*$-algebras.  
	
	There are many situations in which one encounters a category of spaces dually equivalent to a category
	of algebras.  In a wonderful book  \cite{Nestruev}, 
	the mathematicians collectively known as Jet Nestruev develop the theory of smooth real manifolds entirely upside-down: the theory is built 
	in the dual algebraic category, whose objects Nestruev terms smooth complete geometric $\field{R}$-algebras,  
	and then exported to the geometric one by a contravariant functor.  
	In many situations, given a category of spaces dually equivalent to a category of 
	algebras it pays to shift our stance and view the the algebraic category as primary.  In particular, the algebraic
	point of view is typically easier to generalize.  This
	is the paradigm shift driving Alain Connes' non-commutative geometry programme, and the reader
	is referred to \cite{Connes} for much more information.
		
	This paradigm shift is precisely what is needed in order to salvage non-crossing probability theory.
	In probability theory, the notion of space is that of a Kolmogorov triple $(\Omega,\algebra{F},P)$ which
	models the probability to observe a stochastic system in a given state or collection of states.  
	The dual algebraic object associated to a Kolmogorov triple is $L^{\infty}(\Omega,\algebra{F},P)$, 
	the algebra of essentially bounded complex random variables $X:\Omega \rightarrow \field{C}$.  Just like in 
	the case of continuous functions on a compact Hausdorff space, this algebra has a very special
	structure: it is a commutative von Neumann algebra equipped with a unital faithful tracial state,
	$\tau[X]=\int_\Omega X \mathrm{d}P$.  Moreover, there is an analogue of Gelfand's
	theorem in this setting which says that any commutative von Neumann algebra can be realized
	as the algebra of bounded complex random variables on a Kolmogorov triple
	in an essentially unique way.  This is the statement that the categories of Kolmogorov triples and 
	commutative von Neumann algebras are dual equivalent.  
	
	Non-crossing independence was rendered trivial by the commutativity of random variables.
	We can rescue it from the abyss by following the lead of non-commutative geometry and dropping
	commutativity in the dual category: we shift our stance and define a non-commutative 
	probability space to be a pair $(\algebra{A},\tau)$ consisting of a possibly non-commutative 
	complex associative unital algebra $\algebra{A}$ together with a unital linear functional 
	$\tau:\algebra{A} \rightarrow \field{C}$.  If we reinstate commutativity and insist that $\algebra{A}$ is a 
	von Neumann algebra and $\tau$ a faithful tracial state,
	we are looking at essentially bounded random variables on a Kolmogorov triple, but a general
	non-commutative probability space need not be an avatar of any classical probabilistic entity.  
	
	As a nod to the origins of this definition, and in order to foster
	analogies with classical probability, we refer to the elements of $\algebra{A}$ as
	random variables and call $\tau$ the expectation functional.  This prompts some natural questions.
	Before this subsection we only discussed real random variables --- complex numbers 
	crept in with the abstract nonsense.  What is the analogue of the notion of real random variable in a non-commutative
	probability space?  Probabilists characterize random variables in terms of their distributions.  
	Can we assign distributions to random variables living in a non-commutative 
	probability space?  Is it possible to give meaning to the phrase ``the distribution of a bounded real random variable
	living in a non-commutative probability space is a compactly supported probability measure on the line''?
	We will deal with some of these questions at the end of Lecture Two.
	For now, however, we remain in the 
	purely algebraic framework, where the closest thing to the distribution of a random variable $X \in \algebra{A}$ is its
	moment sequence $m_n(X)=\tau[X^n]$.  As in \cite[Page 12]{VDN},
	
		\begin{quote}
		 	The algebraic context is not used in the pursuit of generality, but rather of transparence.
		\end{quote}
															
	\subsection{A brief history of the free world}
	Having cast off the yoke of commutativity, we are free --- free to explore non-crossing
	probability in the new framework provided by the non-commutative probability
	space concept.  Non-crossing probability has become {\sc Free Probability}, and will henceforth be 
	referred to as such.  Accordingly, non-crossing cumulants will now be referred to as 
	free cumulants, and non-crossing independence will be termed free independence.
	
	The reader is likely aware that free probability is a flourishing area of contemporary
	mathematics.  This first lecture has been historical fiction, and is essentially an extended version of  \cite{NovakSniady}.  
	Free probability was not discovered
	in the context of graph enumeration problems, or by tampering with the cumulant concept, although in retrospect
	it might have been.
	Rather, free probability theory was invented by Dan-Virgil Voiculescu in the 1980's in order to address
	a famous open problem in the theory of von Neumann algebras, the free group factors isomorphism 
	problem.  The problem is to determine when the von Neumann algebra of the free group on $a$ generators
	is isomorphic to the von Neumann algebra of the free group on $b$ generators.  It is generally believed that
	these are isomorphic von Neumann algebras if and only if $a=b$, but this remains an open problem. 
	Free probability theory (and its name) originated in this operator-algebraic context.
	
	Voiculescu's definition of free independence, which was modelled on the free product of groups,
	is the following: random variables $X,Y$ in a non-commutative probability space $(\algebra{A},\tau)$
	are said to be freely independent if 
	
		\begin{equation*}
			\tau[f_1(X)g_1(Y) \dots f_k(X)g_k(Y)]=0 
		\end{equation*}
		
	\noindent
	whenever $f_1,g_1,\dots,f_k,g_k$ are polynomials such that 
	
		\begin{equation*}
			\tau[f_1(X)]=\tau[g_1(X)] = \dots=\tau[f_k(X)]=\tau[g_k(Y)]=0.  
		\end{equation*}
			
	\noindent
	This should be compared with the definition of classical independence:
	random variables $X,Y$ in a non-commutative probability space $(\algebra{A},\tau)$ are said to 
	be classically independent if they commute, $XY=YX$, and if 
	
		\begin{equation*}
			\tau[f(X)g(Y)]=0
		\end{equation*}
		
	\noindent
	whenever $f$ and
	$g$ are polynomials such that $\tau[f(X)]=\tau[g(Y)]=0$.  These two definitions are 
	antithetical: classical independence has commutativity built into it, while free independence
	becomes trivial if commutativity is imposed.  Nevertheless, both notions are
	accommodated within the non-commutative probability space framework. 
	
	The precise statement of equivalence between classical independence and vanishing of mixed cumulants is due to 
	Gian-Carlo Rota \cite{Rota}.  In the 1990's, knowing both of Voiculescu's new free probability 
	Theory and Rota's approach to classical probability theory, Roland Speicher made the beautiful
	discovery that by excising the lattice of set partitions from Rota's foundations and replacing it with 
	the lattice of non-crossing partitions, much of Voiculescu's theory could be recovered and extended
	by elementary combinatorial methods.  In particular, Speicher showed that free independence is equivalent
	to the vanishing of mixed free cumulants.  The combinatorial approach to free probability is exhaustively applied
	in \cite{NS}, while the original analytic approach of Voiculescu is detailed in \cite{VDN}.  
			
	\section{Lecture Two: Exploring the Free World}
	Lecture One culminated in the notion of a non-commutative probability space
	and the realization that this framework supports two types of independence: classical independence
	and free independence.  From here we can proceed in several ways.
	One option is to prove an abstract result essentially stating that  these are the only notions of independence which can occur.  
	This result, due to Speicher, places classical and free independence on equal footing. 
	Another possibility is to present concrete problems of intrinsic interest where free independence naturally appears.
	We will pursue the second route, and examine problems emerging from the theory of 
	random walks on groups which can be recast as questions about free random variables. 	
	In the course of solving these problems we will develop the calculus of free random variables
	and explore the terrain of the free world. 
			
	\subsection{Random walk on the integers}
	The prototypical example of a random walk on a group is the simple random walk on $\group{Z}$:
	a walker initially positioned at zero tosses a fair coin at each tick of the clock --- if it lands heads he takes 
	a step of $+1$, if it lands tails he takes a step of $-1$.  
	A random walk is said to be recurrent if it returns to its initial position with probability one, 
	and transient if not.  Is the simple random walk on $\group{Z}$ recurrent or transient?
	
	Let $\alpha(n)$ denote the number of walks which return to zero for the first time after $n$ steps, and 
	let $\phi(n)=2^{-n}\alpha(n)$ denote the corresponding probability that the first return occurs at time $n$.
	Note that $\alpha(0)=\phi(0)=0$, and define
	
		\begin{equation*}
			F(z) = \sum_{n=0}^\infty \phi(n)z^n.
		\end{equation*}
		
	\noindent
	Then 
	
		\begin{equation*}
			F(1) = \sum_{n=0}^{\infty} \phi(n) \leq 1
		\end{equation*}
		
	\noindent
	is the probability we seek.  The radius of convergence of $F(z)$ is at least one, and
	by Abel's theorem 
	
		\begin{equation*}
			F(1) = \lim_{x \rightarrow 1} F(x)
		\end{equation*}
		
	\noindent
	as $x$ approaches $1$ in the interval $[0,1)$.
	
	Let $\lambda(n)$ denote the number of length $n$ loops on $\group{Z}$ based at $0$, 
	and let $\rho(n)=2^{-n}\lambda(n)$
	be the corresponding probability of return at time $n$ (regardless of whether this is the first return
	or not).  Note that $\lambda(0)=\rho(0)=1$.
	We have
	
		\begin{equation*}
			\lambda(n) = \begin{cases}
					0, \text{ if $n$ odd}\\
					{n \choose \frac{n}{2}}, \text{ if $n$ even}
				\end{cases}.
		\end{equation*}
		
	\noindent
	From Stirling's formula, we see that 
	
		\begin{equation*}
			\rho(2k) \sim \frac{1}{\sqrt{\pi k}}
		\end{equation*}
		
	\noindent
	as $k \rightarrow \infty$.  Thus the radius of convergence of
	
		\begin{equation*}
			R(z) = \sum_{n=0}^{\infty} \rho(n)z^n
		\end{equation*}
		
	\noindent
	is one.
	
	We can decompose the set of loops of given length according to the number of
	steps taken to the first return.  This produces the equation
	
		\begin{equation*}
			\lambda(n) = \sum_{k=0}^n \alpha(k)\lambda(n-k).
		\end{equation*}
		
	\noindent
	Equivalently, since all probabilities are uniform,
			
		\begin{equation*}
			\rho(n) = \sum_{k=0}^n \phi(k)\rho(n-k).
		\end{equation*}
		
	\noindent
	Summing on $z$, this becomes the identity
	
		\begin{equation*}
			R(z)-1=F(z)R(z)
		\end{equation*}
		
	\noindent
	in the algebra of holomorphic functions on the open unit disc in $\field{C}$.
	Since $R(z)$ has non-negative coefficients, it is non-vanishing for $x \in [0,1)$ and we can write

		\begin{equation*}
			F(x) = 1-\frac{1}{R(x)}, \quad 0 \leq x < 1.
		\end{equation*}
	
	\noindent
	Thus
	
		\begin{equation*}
			F(1) = \lim_{x \rightarrow 1} F(x) = 1-\frac{1}{\lim_{x \rightarrow 1} R(x)}.
		\end{equation*}
		
	\noindent
	If $R(1)<\infty$, then by Abel's theorem $\lim_{x \rightarrow 1} R(x)=R(1)$ and we 
	obtain $F(1)<1$.  On the other hand, if 
	$R(1)=\infty$, then $\lim_{x \rightarrow 1} R(x) = \infty$ and we get $F(1)=1$.  Thus
	the simple random walk is transient or recurrent according to 
	the convergence or divergence of the series $\sum \rho(n)$.		
	From the Stirling estimate above we find that this sum diverges, so the simple random 
	walk on $\group{Z}$ is recurrent.
					
	\subsection{P\'olya's theorem}
	In the category of abelian groups, coproduct is direct sum:
	
		\begin{equation*}
			\coprod_{i \in I} \group{G}_i = \bigoplus_{i \in I} \group{G}_i.
		\end{equation*}
	
	\noindent	
	In 1921, George P\'olya  \cite{Polya} proved that the simple random 
	walk on 
		
		\begin{equation*}
			\group{Z}^d = \underbrace{\group{Z} \oplus \dots \oplus \group{Z}}_d
		\end{equation*}

	\noindent
	is recurrent for $d=1,2$ and transient for $d>2$.  This striking result can 
	be deduced solely from an understanding of the simple random walk on $\group{Z}$.
		
	Let us give a proof of P\'olya's theorem.  
	Let $\lambda_d(n)$ denote the number of length $n$ loops on $\group{Z}^d$ based at
	$\vector{0}^d$.  Let $\rho_d(n)$ denote the probability of return to $\vector{0}^d$ 
	after $n$ steps,
	
		\begin{equation*}
			\rho_d(n) = \frac{1}{(2d)^n} \lambda_d(n).		
		\end{equation*}
		
	\noindent
	As above, the simple random walk on $\group{Z}^d$ is recurrent 
	if the sum $\sum \rho_d(n)$ diverges, and transient otherwise.
	Form the loop generating function
	
		\begin{equation*}
			L_d(z) = \sum_{n=0}^{\infty} \lambda_d(n)z^n.
		\end{equation*}	
		
	\noindent
	We aim to prove that 
	
		\begin{equation*}
			L_d\bigg{(}\frac{1}{2d}\bigg{)}=\sum_{n=0}^{\infty} \rho_d(n)
		\end{equation*} 
	
	\noindent	
	diverges for $d=1,2$ and converges for $d>2$.
	
	While the ordinary loop generating function is hard to analyze directly, the exponential 
	loop generating function
	
		\begin{equation*}
			E_d(z) = \sum_{n=0}^{\infty} \lambda_d(n)\frac{z^n}{n!}
		\end{equation*}
		
	\noindent
	is quite accessible.  Indeed, as in the last subsection we have 
	
		\begin{equation*}
			\lambda_1(n) = \begin{cases}
				0, \text{ if $n$ odd}\\
				{n \choose \frac{n}{2}}, \text{ if $n$ even}
				\end{cases},
		\end{equation*}
		
	\noindent
	so that
	
		\begin{equation*}
			E_1(z) = \sum_{k=0}^{\infty} \frac{z^{2k}}{k!k!}=I_0(2z)
		\end{equation*}
		
	\noindent
	is precisely the modified Bessel function of order zero.
	Since a loop on $\group{Z}^d$ is just a shuffle of loops on $\group{Z}$, 
	the product formula for exponential generating functions \cite{Stanley:EC2} yields
	
		\begin{equation*}
			E_d(z) = E_1(z)^d = I_0(2z)^d.
		\end{equation*}
	
	What we have is the exponential generating function for the loop counts $\lambda_d(n)$,
	and what we want is the ordinary generating function of this sequence.	
	The integral transform 
	
		\begin{equation*}
			L_f(z) = \int\limits_0^{\infty} f(tz)e^{-t} \mathrm{d}t,
		\end{equation*}
		
	\noindent
	which looks like the Laplace transform of $f$ but with the $z$-parameter in the wrong place,
	converts exponential generating functions into ordinary generating functions.
	This can be seen by differentiating under the integral sign and using the fact that
	the moments of the exponential distribution are the factorials,
	
		\begin{equation*}
			\int\limits_0^{\infty} t^n e^{-t} \mathrm{d}t = n!.
		\end{equation*}
		
	\noindent
	This trick is constantly used in quantum field theory in connection with Borel summation of divergent
	series \cite{Etingof}.
	In particular, we have
	
		\begin{equation*}
			L_d(z) = \int\limits_0^\infty E_d(tz)e^{-t} \mathrm{d}t = \int\limits_0^\infty I_0(2tz)^de^{-t} \mathrm{d}t.
		\end{equation*}
		
	\noindent
	Thus it remains only to show that the integral
	
		\begin{equation*}
			L_d\bigg{(}\frac{1}{2d}\bigg{)} = \int\limits_0^\infty I_0 \bigg{(}\frac{t}{d}\bigg{)}^de^{-t} \mathrm{d}t
		\end{equation*}
		
	\noindent
	is divergent for $d=1,2$ and convergent for $d >2$.  This in turn amounts to understanding the 
	asymptotics of $I_0(t/d)$ as $t \rightarrow \infty$ along the real line.  
	
	We already encountered
	Bessel functions in Lecture One, and we know that
	
		\begin{equation*}
			I_0(t/d) = \frac{1}{\pi} \int\limits_0^\pi e^{t(\frac{\cos \theta}{d})} \mathrm{d}\theta.
		\end{equation*}
		
	\noindent
	This is an integral of Laplace type, 
	
		\begin{equation*}
			\int\limits_a^b e^{tf(\theta)} \mathrm{d}\theta,
		\end{equation*}
	
	\noindent
	and Laplace integrals localize as $t \rightarrow \infty$ with asymptotics given
	by the classical steepest descent formula (maximum at an endpoint case),
	
		\begin{equation*}
			\int\limits_a^b e^{tf(\theta)} \mathrm{d}\theta \sim \sqrt{\frac{\pi}{2t|f''(a)|}}e^{tf(a)}.
		\end{equation*} 	
		
	\noindent
	For our integral, this specializes to 
	
		\begin{equation*}
			I_0(t/d) \sim \sqrt{\frac{1}{2\pi^{3/2} t}} e^{t/d}, \quad t \rightarrow \infty,
		\end{equation*}				
		
	\noindent
	from which it follows that $L_d((2d)^{-1})$ diverges or converges according to the divergence 
	or convergence of the integral
	
		\begin{equation*}
			\int\limits_1^\infty t^{-d/2}\mathrm{d}t.
		\end{equation*}		
		
	\noindent
	This integral diverges for $d=1,2$ and converges for $d \geq 3$, which proves P\'olya's result. 
	In fact, the probability that the simple random walk on $\group{Z}^3$ returns to its initial position
	is already less than thirty five percent.
	
	\subsection{Kesten's problem}
	The category of abelian groups is a full subcategory of the category of groups.
	In the category of groups, coproduct is free product:
	
		\begin{equation*}
			\coprod_{i \in I} \group{G}_i = *_{i \in I} \group{G}_i.
		\end{equation*}

\begin{figure}
			\includegraphics{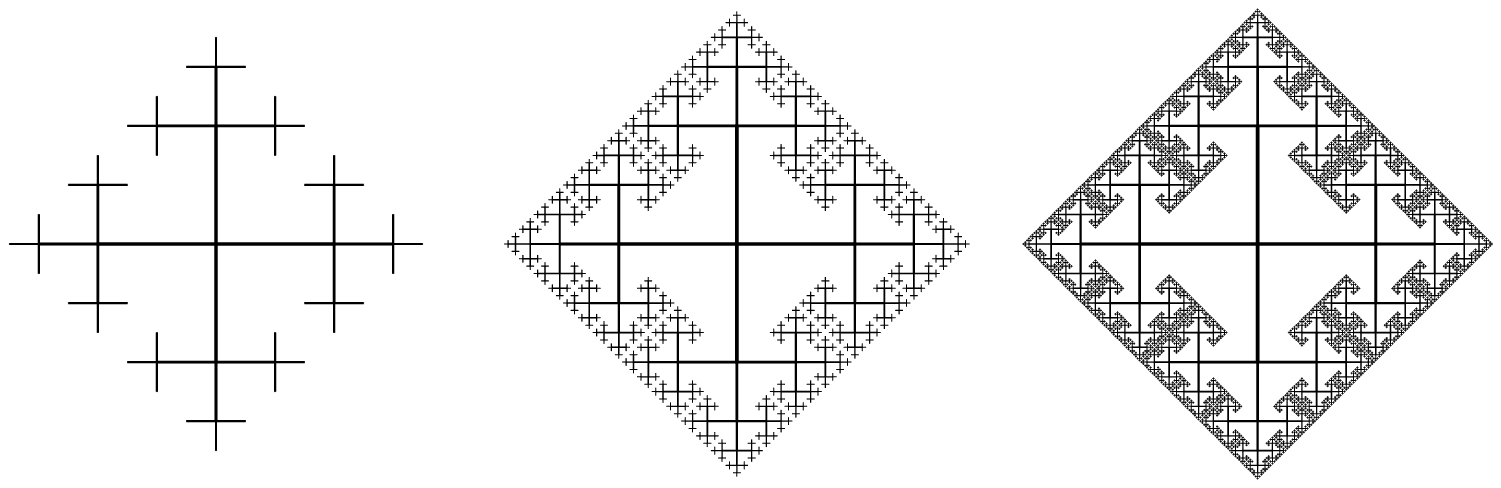}
			\caption{\label{fig:Kesten}Balls of increasing radius in $\group{F}_2$.}
	\end{figure}

	\noindent
	Thus one could equally well ask about the recurrence or transience of the simple
	random walk on 
	
		\begin{equation*}
			\group{F}_d = \underbrace{\group{Z} * \dots * \group{Z}}_d,
		\end{equation*}
		
	\noindent
	the free group on $d$ generators.  Whereas the Cayley graph of the abelian group
	$\group{Z}^d$ is the $(2d)$-regular hypercubic lattice, the Cayley graph of the free group $\group{F}_d$
	is the $(2d)$-regular tree, see Figure \ref{fig:Kesten}.  What is the free analogue of P\'olya's
	theorem?  We will see that the random walk on $\group{F}_d$ 
	can be understood entirely in terms of the random walk on $\group{F}_1 = \group{Z}$,
	just like in the abelian category.  However, the tools we will use are quite different,
	and the concept of free random variables plays the central role.
	
	The study of random walks on groups was initiated by Harry Kesten in his 1958 Ph.D. thesis,
	with published results appearing in \cite{Kesten}.  A good source of information on this topic,
	with many pointers to the literature, is Laurent Saloff-Coste's survey article \cite{SC}.
	Kesten related the behaviour of the simple random walk on a finitely-generated group
	$\group{G}$ to other properties of $\group{G}$, such as amenability.  A countable group
	is said to be amenable if it admits a finitely additive $\group{G}$-invariant 
	probability measure.  The notion of amenability was introduced by John von Neumann 
	in 1929.  Finite groups are amenable since they can be equipped with the 
	uniform measure $P(g)=|\group{G}|^{-1}$.  For infinite groups the situation is not so clear,
	and many different characterizations of amenability have been derived.
	For example, Alain Connes showed that a group is amenable if and only if its von Neumann 
	algebra is hyperfinite.
	Kesten proved that $\group{G}$ is non-amenable if and only if the probability $\rho_\group{G}(n)$
	that the simple random walk on $\group{G}$ returns to its starting point at time $n$ 
	decays exponentially in $n$.  We saw above that for $\group{G}=\group{Z}$ the return 
	probability has square root decay, so $\group{Z}$ is amenable.  In fact, amenability is preserved by
	direct sum so all abelian groups are amenable.  Is the free group $\group{F}_d$ amenable?  Let $\lambda_d(n)$ denote
	the number of length $n$ loops on $\group{F}_d$ based at $\id$.  We will refer to the 
	problem of finding an explicit expression for the loop generating function
		
		\begin{equation*}
			L_d(z) = 1+\sum_{n=1}^{\infty} \lambda_d(n)z^n
		\end{equation*}
		
	\noindent
	as Kesten's problem.  Presumably, if we can obtain an explicit expression for this function
	then we can read off the asymptotics of $\rho_d(n)$, which is the coefficient of $z^n$ in 
	$L_d(z/2d)$, via the usual methods of singularity analysis of generating functions.
	
	We begin at the beginning: $d=2$.
	Let $A$ and $B$ denote the generators of $\group{F}_2$, and let $\algebra{A}=\algebra{A}[\group{F}_2]$
	be the group algebra consisting of formal $\field{C}$-linear combinations of words in these
	generators and their inverses, $A^{-1}$ and $B^{-1}$.  The identity element of $\algebra{A}$ is the 
	empty word, which is identified with $\id \in \group{F}_2$.  Introduce the expectation functional
	
		\begin{equation*}
			\tau[X] = \text{coefficient of $\id$ in $X$}
		\end{equation*}
		
	\noindent
	for each $X \in \algebra{A}$.  Then $(\algebra{A},\tau)$ is a non-commutative probability space.
	A loop $\id \rightarrow \id$ in $\group{F}_2$
	is simply a word in $A,A^{-1},B,B^{-1}$ which reduces to $\id$.  Thus the number of
	length $n$ loops in $\group{F}_2$  is 
	
		\begin{equation*}
			\lambda_2(n)=m_n(X+Y) = \tau[(X+Y)^n],
		\end{equation*}
		
	\noindent
	where $X,Y \in \algebra{A}$ are the random variables
	
		\begin{equation*}
			X=A+A^{-1}, \quad Y=B+B^{-1}.
		\end{equation*}
		
	\noindent
	We see that the loop generating function for $\group{F_2}$ is precisely the 
	moment generating function for the random variable $X+Y$ in the non-commutative
	probability space $(\algebra{A},\tau)$,
	
		\begin{equation*}
			L_2(z) = 1 + \sum_{n=1}^{\infty} m_n(X+Y) z^n.
		\end{equation*}
				
	We want to compute the moments of the sum $X+Y$ of two non-commutative random variables, 
	and what we know are the moments of its summands:
	
		\begin{equation*}
			m_n(X)=m_n(Y) = \begin{cases}
				0, \text{ if $n$ odd} \\
				{n \choose \frac{n}{2}}, \text{ if $n$ even}
			\end{cases}.
		\end{equation*}
		
	\noindent
	Now we make the key observation: the random variables $X,Y$ are freely independent.
	Indeed, suppose that $f_1,g_1,\dots,f_k,g_k$ are polynomials such that 
	
		\begin{equation*}
			\tau[f_1(X)]=\tau[g_1(Y)] = \dots = \tau[f_k(X)]=\tau[g_k(Y)] =0.
		\end{equation*}
		
	\noindent
	This means that $f_i(X)=f_i(A+A^{-1})$ is a Laurent polynomial in $A$ with zero constant 
	term, and $g_j(Y)=g_j(B+B^{-1})$ is a Laurent polynomial in $B$ with zero constant term.  Since there are no relations
	between $A$ and $B$, an alternating product of polynomials of this form cannot produce any occurrences
	of the empty word, and we have
	
		\begin{equation*}
			\tau[f_1(X)g_1(Y) \dots f_k(X)g_k(Y)]=0.
		\end{equation*}
		
	\noindent
	This is precisely Voiculescu's definition of free independence.  
	
	We conclude that the problem of computing $\lambda_2(n)$ 
	is a particular case of the problem of computing the moments
	$m_n(X+Y)$ of the sum of two free random variables given their individual moments, $m_n(X)$ and $m_n(Y)$. 
	This motivates us to solve a fundamental problem in free probability 
	theory:
	
		\begin{center}
			\emph{
			Given a pair of free random variables $X$ and $Y$, compute the moments of 
			$X+Y$ in terms of the moments of $X$ and the moments of $Y$}.
		\end{center}
		
	We can, in principle, solve this problem using the fact that free cumulants linearize the addition of free random 
	variables, $\kappa_n(X+Y)=\kappa_n(X)+\kappa_n(Y)$.  This solution is implemented as the following recursive algorithm.
	
		\begin{description}
			
			\item[Input] 
				$\kappa_1(X),\dots,\kappa_{n-1}(X),\kappa_1(Y),\dots,\kappa_{n-1}(Y)$.
				
					\bigskip
				
				\item[\emph{Step 1}]
				Compute $m_n(X),m_n(Y)$.
				
					\medskip
				
				\item[\emph{Step 2}]
				Compute $\kappa_n(X),\kappa_n(Y)$ using
					
					\begin{equation*}
						\begin{split}
						\kappa_n(X) &= m_n(X) - \sum_{\substack{\pi \in \lattice{NC}(n)\\ b(\pi) \geq 2}} \prod_{B \in \pi} \kappa_{|\beta|}(X) \\
						\kappa_n(Y) &= m_n(Y) - \sum_{\substack{\pi \in \lattice{NC}(n)\\ b(\pi) \geq 2}} \prod_{B \in \pi} \kappa_{|\beta|}(Y).
						\end{split}
					\end{equation*}
					
					\medskip
					
				\item[\emph{Step 3}]
				Add,
				
					\begin{equation*}
						\kappa_n(X+Y) = \kappa_n(X) + \kappa_n(Y).
					\end{equation*}
				
					\medskip
				
				\item[\emph{Step 4}]
				Compute $m_n(X+Y)$ using
				
					\begin{equation*}
						m_n(X+Y) = \kappa_n(X+Y) + \sum_{\substack{\pi \in \lattice{NC}(n)\\ b(\pi \geq 2}} \prod_{B \in \pi} \kappa_{|B|}(X+Y).
					\end{equation*}
					
					\bigskip
			
			\item[Output] 
			$m_n(X+Y)$.
						
		\end{description}
		
		This recursive algorithm is conceptually simple but virtually useless as is.  In particular, it is not clear how 
		to coax it into computing the loop generating function $L_2(z)$.  We need to develop an additive calculus of 
		free random variables which parallels the additive calculus of classically independent random variables.
				
		\subsection{The classical algorithm}
		If $X,Y$ are classically independent random variables, we can compute the moments of their 
		sum $X+Y$ using the recursive algorithm above, replacing free cumulants with classical cumulants.
		But this is not what probabilists do in their daily lives.  They have a much better algorithm which uses
		analytic function theory to efficiently handle the recursive nature of the naive algorithm.  The classical 
		algorithm associates to $X$ and $Y$ analytic functions $M_X(z)$ and $M_Y(z)$ which have the 
		property that $M_{X+Y}(z):=M_X(z)M_Y(z)$ encodes the moments of $X+Y$ as its derivatives 
		at $z=0$.  We will give a somewhat roundabout derivation of this algorithm, which is presented in this way
		specifically to highlight the analogy with Voiculescu's algorithm presented in the next section.
		
		The classical algorithm for summing two random variables is developed in two stages.  In the first stage,
		the relation between the moments and classical cumulants of a random variable is packaged as an identity
		in the ring of formal power series $\field{C}[[z]]$.  Suppose that $(m_n)_{n=1}^{\infty}$ and $(c_n)_{n=1}^{\infty}$
		are two numerical sequences related by the chain of identities
		
			\begin{equation*}
				m_n = \sum_{\pi \in \lattice{P}(n)} \prod_{B \in \pi} c_{|B|}, \quad n \geq 1.
			\end{equation*}
			
		\noindent
		The $\pi^{\text{th}}$ term of the sum on the right only depends on the ``spectrum''
		of $\pi$, i.e. the integer vector $\Lambda(\pi) = (1^{b_1(\pi)},2^{b_2(\pi)},\dots,n^{b_n(\pi)}),$ 
		where $b_i(\pi)$ is the
		number of blocks of size $i$ in $\pi$.  We may view $\Lambda(\pi)$ as the Young diagram 
		with $b_i$ rows of length $i$.  Consequently, we can perform a change of variables to push the summation
		forward onto a sum over Young diagrams with $n$ boxes provided we can compute the ``Jacobian''
		of the map $\Lambda:\lattice{P}(n) \rightarrow \lattice{Y}(n)$ sending $\pi$ on its spectrum:
						
			\begin{equation*}
				m_n = \sum_{b_1+2b_2+ \dots + nb_n=n} c_1^{b_1}c_2^{b_2} \dots c_n^{b_n}
				|\Lambda^{-1}(1^{b_1},2^{b_2},\dots,n^{b_n})|.
			\end{equation*}
			
		\noindent
		The volume of the fibre of $\Lambda$ over any given Young diagram can be explicitly computed,
		
			\begin{equation*}
				|\Lambda^{-1}(1^{b_1},2^{b_2},\dots,n^{b_n})| = \frac{n!}{(1!)^{b_1}(2!)^{b_2} \dots (n!)^{b_n}b_1!b_2!\dots b_n!},
			\end{equation*}
			
		\noindent
		so that we have the chain of identities
		
			\begin{equation*}
				\frac{m_n}{n!} = \sum_{b_1 + 2b_2 + \dots +nb_n=n}\frac{(c_1/1!)^{b_1} (c_2/2!)^{b_2} \dots 
				(c_n/n!)^{b_n}}{b_1!b_2! \dots b_n!}, \quad n \geq 1.
			\end{equation*}
			
		\noindent
		We can bundle these identities as a single relation between power series.  Summing on $z$ we obtain
		
			\begin{equation*}
				\begin{split}
				1+\sum_{n=1}^{\infty} m_n \frac{z^n}{n!} &= 1+\sum_{n=1}^{\infty} 
				\bigg{(} \sum_{b_1 + 2b_2 + \dots +nb_n=n}\frac{(c_1/1!)^{b_1} (c_2/2!)^{b_2} \dots 
				(c_n/n!)^{b_n}}{b_1!b_2! \dots b_n!} \bigg{)} z^n \\
				&=1 + \frac{1}{1!}\bigg{(} \sum_{n=1}^{\infty} c_n \frac{z^n}{n!} \bigg{)}^1+
				\frac{1}{2!}\bigg{(}\sum_{n=1}^{\infty} c_n\frac{z^n}{n!} \bigg{)}^2+
				\dots \\
				&= e^{\sum_{n=1}^{\infty} c_n\frac{z^n}{n!}}.
				\end{split}
			\end{equation*}
			
		\noindent
		We conclude that the chain of moment-cumulant formulas is equivalent to the 
		single identity $M(z)=e^{C(z)}$ in $\field{C}[[z]]$, where
			
			\begin{equation*}
				M(z)=1+\sum_{n=1}^{\infty} m_n\frac{z^n}{n!}, \quad C(z) = \sum_{n=1}^{\infty} c_n \frac{z^n}{n!}
			\end{equation*}
			
		\noindent
		This fact is known in enumerative 
		combinatorics as the exponential formula.  In other branches of science it goes
		by other names, such as the the polymer expansion formula or the linked cluster theorem.
		In the physics literature, the exponential formula is often invoked using colourful phrases
		such as ``connected vacuum bubbles exponentiate'' \cite{Samuel}.
		The exponential 
		formula seems to have been first written down precisely by Adolf Hurwitz in 1891 \cite{Hurwitz}.

		The exponential formula becomes particularly powerful when combined with complex analysis.	
		Suppose that $X,Y$ are classically independent random variables living in a non-commutative
		probability space $(\algebra{A},\tau)$.  Suppose moreover that an oracle has given us probability
		measures $\mu_X,\mu_Y$ on the real line which behave like distributions for $X,Y$ insofar as
		
			\begin{equation*}
				\tau[X^n] = \int\limits_{\field{R}} t^n \mu_X(\mathrm{d}t), \quad 
				\tau[Y^n] = \int\limits_{\field{R}} t^n \mu_Y(\mathrm{d}t), \quad n \geq 1.
			\end{equation*}
		
		\noindent
		Let us ask for even more, and insist that $\mu_X,\mu_Y$ are compactly supported. 
		Then the functions\footnote{The restriction of $M_X$ to the real axis, $M_X(-x)$, is the 
		two-sided Laplace transform, while the restriction of $M_X$ to the imaginary axis,
		$M_X(-iy)$, is the Fourier transform.}
		
			\begin{equation*}
				M_X(z) = \int\limits_{\field{R}} e^{tz}\mu_X(\mathrm{d}t), \quad
				M_Y(z) = \int\limits_{\field{R}} e^{tz}\mu_Y(\mathrm{d}t)
			\end{equation*}
			
		\noindent
		are entire, and their derivatives can be 
		computed by differentiation under the integral sign.
		Consequently, we
		have the globally convergent power series expansions
		
			\begin{equation*}
				M_X(z) = 1+\sum_{n=1}^{\infty} m_n(X) \frac{z^n}{n!}, \quad 
				M_Y(z) = 1+\sum_{n=1}^{\infty} m_n(Y) \frac{z^n}{n!}.
			\end{equation*}
		
		Since $M_X(0)=M_Y(0)=1$ and the zeros of holomorphic functions are discrete,
		we can restrict to a complex domain $\domain{D}$ containing the origin on which $M_X(z),M_Y(z)$ 
		are non-vanishing.  Let $\Hol(\domain{D})$ denote the algebra of holomorphic functions on 
		$\domain{D}$.  The following algorithm produces a function $M_{X+Y}(z) \in \Hol(\domain{D})$
		whose derivatives at $z=0$ are the moments of $X+Y$.
											
			\begin{description}
			
				\item[Input] 
					$\mu_X$ and $\mu_Y$.
				
						\bigskip
						
					\item[\emph{Step 1}]
					Compute
					
						\begin{equation*}
							M_X(z) =\int\limits_{\field{R}} e^{tz} \mu_X(\mathrm{d}t),\ 
							M_Y(z) =\int\limits_{\field{R}} e^{tz} \mu_Y(\mathrm{d}t).
						\end{equation*}
				
					\item[\emph{Step 2}]
					Solve
					
						\begin{equation*}
							M_X(z) = e^{C_X(z)}, \quad M_Y(z) = e^{C_Y(z)}
						\end{equation*}
						
					\noindent
					in $\Hol(\domain{D})$ subject to $C_X(0)=C_Y(0)=0$.
								
						\medskip
							
					\item[\emph{Step 3}]
					Add, 
					
						$$C_{X+Y}(z) := C_X(z) + C_Y(z).$$
								
						\medskip
				
					\item[\emph{Step 4}]
					Exponentiate, 
					
						$$M_{X+Y}(z):=e^{C_{ X+Y}(z)}.$$
									
						\bigskip
			
				\item[Output] 
				$M_{X+Y}(z)$.
				
			\end{description}
		
		In Step One, we try to compute the integral transforms $M_X(z),M_Y(z)$ in terms of elementary functions,
		like $e^z,\log(z), \sin(z),\cos(z),\sinh(z),\cosh(z),\dots$ etc, or other classical functions like Bessel functions,
		Whittaker functions, or anything else that can be looked up in \cite{AAR}.  
		This is often feasible if the distributions $\mu_X,\mu_Y$ have known densities, and
		we saw some examples in Lecture One.
				
		The equations in Step Two have unique solutions.  The required functions
		$C_X(z),C_Y(z) \in \Hol(\domain{D})$ are the principal branches of the logarithms of 
		$M_X(z),M_Y(z)$ on $\domain{D}$, and can be represented as contour integrals: 
		
			\begin{equation*}
				C_X(z) = \log M_X(z) = \oint_0^z \frac{M_X'(\zeta)}{M_X(\zeta)}\mathrm{d}\zeta, \quad
				C_Y(z) = \log M_Y(z) = \oint_0^z \frac{M_Y'(\zeta)}{M_Y(\zeta)}\mathrm{d}\zeta
			\end{equation*}
			
		\noindent
		for $z \in \domain{D}$.  Since $\log$ has the usual formal properties associated with the 
		logarithm, if Step One outputs a reasonably explicit expression then so will Step Two.
		
		Step Two is the crux of the algorithm.  It is performed precisely to change gears from a moment 
		computation to a cumulant computation.  Appealing to the exponential formula, we 
		conclude that the holomorphic functions $C_X(z),C_Y(z)$ passed to 
		Step Three by Step Two have Maclaurin series
		
			\begin{equation*}
				C_X(z)= \sum_{n=1}^{\infty} c_n(X) \frac{z^n}{n!}, \quad C_Y(z)= \sum_{n=1}^{\infty} c_n(Y) \frac{z^n}{n!},
			\end{equation*}
			
		\noindent
		where $c_n(X),c_n(Y)$ are the cumulants of $X$ and $Y$.  Since cumulants linearize the addition of 
		independent random variables, the new function $C_{X+Y}(z) :=C_X(z) + C_Y(z)$ defined in Step Three
		encodes the cumulants of $X+Y$ as its derivatives at $z=0$.
					
		In Step Four we define a new function $M_{X+Y}(z) \in \Hol(\domain{D})$ by 
		$M_{X+Y}(z) := e^{C_{X+Y}(z)}$.  The exponential formula and the moment-cumulant formula now combine
		in the reverse direction to tell us 
		that the Maclaurin series of $M_{X+Y}(z)$ is
		
			\begin{equation*}
				M_{X+Y}(z) = 1+\sum_{n=1}^{\infty} m_n(X+Y) \frac{z^n}{n!}.
			\end{equation*}
		
		In summary, assuming that $X,Y$ are classically independent random variables living in a non-commutative
		probability space $(\algebra{A},\tau)$ with affiliated distributions $\mu_X,\mu_Y$ having nice properties,
		the classical algorithm takes these distributions as input and outputs a function $M_{X+Y}(z)$ 
		analytic at $z=0$ whose derivatives are the moments of $X+Y$.  It works by combining the exponential 
		formula and the moment-cumulant formula to convert the moment problem into the (linear) cumulant problem, adding,
		and then converting back to moments.  An optional Fifth Step is 
		to extract the distribution $\mu_{X+Y}$ from $M_{X+Y}(z)$ using the Fourier inversion formula:
		
			\begin{equation*}
				\mu_{X+Y}([a,b]) = \lim_{T \rightarrow \infty} \frac{1}{2\pi}
				\int\limits_{-T}^T \frac{e^{-iat}-e^{-ibt}}{it}M_{X+Y}(it)\mathrm{d}t.
			\end{equation*}
																			
		\subsection{Voiculescu's algorithm}
		We wish to develop a free analogue of the classical algorithm.  
		Suppose that $X,Y$ are freely independent random variables living in a non-commutative
		probability space $(\algebra{A},\tau)$ possessing compactly supported real 
		distributions $\mu_X,\mu_Y$.
		The free algorithm should take these distributions as input,
		build a pair of analytic functions which encode the moments of $X$ and $Y$ respectively,
		and then convolve these somehow to produce a new analytic function which encodes
		the moments of $X+Y$.  A basic hurdle to be overcome is that, even assuming we know 
		how to construct $\mu_X$ and $\mu_Y$, we don't know what to do with them.  We could
		repeat Step One of the classical algorithm to obtain analytic functions $M_X(z),M_Y(z)$ whose 
		derivatives at $z=0$ are the moments of $X$ and $Y$.  If we then perform Step Two we 
		obtain analytic functions $C_X(z),C_Y(z)$ whose derivatives encode the classical cumulants
		of $X$ and $Y$.  But classical cumulants do not linearize addition of free random variables.
		
		The classical algorithm is predicated on the existence of a formal power series identity equivalent
		to the chain of classical moment-cumulant identities.  
		We need a free analogue of this, namely a power series identity equivalent to the chain of 
		numerical identities
		
			\begin{equation*}
				m_n = \sum_{\pi \in \lattice{NC}(n)} \prod_{B \in \pi} \kappa_{|B|}, \quad n \geq 1.
			\end{equation*}	
			
		\noindent
		Proceeding as in the classical case, rewrite this in the form
		
			\begin{equation*}
				m_n = \sum_{b_1+2b_2+ \dots + nb_n=n} \kappa_1^{b_1}\kappa_2^{b_2} \dots \kappa_n^{b_n}
				|\Lambda^{-1}(1^{b_1},2^{b_2},\dots,n^{b_n}) \cap \lattice{NC}(n)|,
			\end{equation*}
			
		\noindent
		where as above $\Lambda:\lattice{P}(n) \rightarrow \lattice{Y}(n)$ is the surjection which sends a partition $\pi$
		with $b_i$ blocks of size $i$ to the Young diagram with $b_i$ rows of length $i$.  Now we 
		have to compute the volume of the fibres of $\Lambda$ intersected with the non-crossing partition lattice.  
		The solution to this
		enumeration problem is again known in explicit form,
				
			\begin{equation*}
				|\Lambda^{-1}(1^{m_1},2^{m_2},\dots,n^{m_n}) \cap \lattice{NC}(n)| =
				\frac{n!}{(n+1-(b_1+b_2+\dots+b_n))!b_1!b_2! \dots b_n!}.
			\end{equation*}
			
		\noindent
		This formula allows us to obtain the desired power series identity, 
		though the manipulations required are
		quite involved and require either the use of Lagrange inversion or an understanding of
		the poset structure of $\lattice{NC}(n)$.  In any event,
		what ultimately comes out of the computation is the fact that two numerical 
		sequences satisfy the chain of free moment-cumulant identities if and only if the 
		ordinary (not exponential) generating functions
		
			\begin{equation*}
				L(z) = 1+\sum_{n=1}^{\infty} m_nz^n, \quad K(z) = 1 + \sum_{n=1}^{\infty}
				\kappa_n z^n
			\end{equation*}
			
		\noindent
		solve the equation
		
			\begin{equation*}
				L(z)=K(zL(z))
			\end{equation*}	
			
		\noindent
		in the formal power series ring $\field{C}[[z]]$.  This is the free analogue of the exponential formula.
			
		As in the classical case, we wish to turn this formal power series encoding into an analytic encoding.
		Suppose that $X,Y$ admit distributions $\mu_X,\mu_Y$ supported in the real interval $[-r,r]$.  We then 
		have $|m_n(X)|,|m_n(Y)| \leq r^n$, so the moment generating functions 
		
			\begin{equation*}
				L_X(z) = 1 + \sum_{n=1}^{\infty} m_n(X)z^n, \quad L_Y(z) = 1 + \sum_{n=1}^{\infty} m_n(Y)z^n,
			\end{equation*}
			
		\noindent
		are absolutely convergent in the open disc $\domain{D}(0,\frac{1}{r})$.  One can use the relation between
		moments and free cumulants to show that the free cumulant generating functions
		
			\begin{equation*}
				K_X(z) = 1 + \sum_{n=1}^{\infty} \kappa_n(X)z^n, \quad K_Y(z) = 1 + \sum_{n=1}^{\infty} \kappa_n(Y)z^n
			\end{equation*}
			
		\noindent
		are also absolutely convergent on a (possibly smaller) neighbourhood of $z=0$.
		However, it turns out that the correct environment for the free algorithm is a neighbourhood 
		of infinity rather than a neighbourhood of zero.
		This is because what we really want
		is an integral transform which realizes ordinary
		generating functions in the same way as the Fourier (or Laplace) transform realizes
		exponential generating functions.  Access to such a transform will allow us to obtain closed forms for 
		generating functions by evaluating integrals, just like in classical probability.
		Such an object is well-known in analysis,
		where it goes by the name of the Cauchy (or Stieltjes) transform.  The Cauchy transform of a
		random variable $X$ with real distribution $\mu_X$ is
		
			\begin{equation*}
				G_X(z) = \int\limits_{\field{R}} \frac{1}{z-t} \mu_X(\mathrm{d}t).
			\end{equation*}
			
		\noindent
		The Cauchy transform is well-defined on the complement of the support of $\mu_X$, and differentiating
		under the integral sign shows that $G_X(z)$ is holomorphic on its domain of definition.  
		In particular, if $\mu_X$ is supported in $[-r,r]$ then
		$G_X(z)$ admits the convergent Laurent expansion 
		
			\begin{equation*}
				G_X(z) = \frac{1}{z}\sum_{n=0}^{\infty} \frac{\int t^n \mu_X(\mathrm{d}t)}{z^n} =
				\sum_{n=0}^{\infty} \frac{m_n(X)}{z^{n+1}}
			\end{equation*}
			
		\noindent
		on $|z|>r$.  This is an ordinary generating function for the moments of $X$ with $z^{-1}$ 
		playing the role of the formal variable.  
		
		To create an interface between the free moment-cumulant formula and the Cauchy 
		transform, we must  re-write the formal power series identity $L(z)=K(zL(z))$ as an identity in 
		$\field{C}((z))=\operatorname{Quot}\field{C}[[z]]$, the field of formal 
		Laurent series.  Introduce the formal Laurent series
		
			\begin{equation*}
				G(z) =\frac{1}{z}L(\frac{1}{z})= \sum_{n=0}^{\infty} \frac{m_n}{z^{n+1}}.
			\end{equation*}
		
		\noindent
		The automorphism $z \mapsto \frac{1}{z}$ transforms
		the non-crossing exponential formula into the identity 
		
			\begin{equation*}
				\frac{K(G(z))}{G(z)} = z.
			\end{equation*}
			
		\noindent
		Setting 
				
			\begin{equation*}
				V(z) = \frac{K(z)}{z} = \frac{1}{z} + \sum_{n=0}^{\infty} \kappa_{n+1}z^n,
			\end{equation*}
			
		\noindent
		this becomes the identity
		
			\begin{equation*}
				V(G(z))=z
			\end{equation*}

		\noindent
		in $\field{C}((z))$.
				
		We have now associated two analytic functions to $X$.  The first is the Cauchy transform $G_X(z)$, which
		is defined as an integral transform and admits a convergent Laurent expansion in a neighbourhood
		of infinity in the $z$-plane.  The second is the Voiculescu transform $V_X(w)$, which is defined by the 
		convergent Laurent series 
		
			\begin{equation*}
				V_X(w) = \frac{1}{w} + \sum_{n=0}^{\infty} \kappa_{n+1}w^n
			\end{equation*}
			
		\noindent
		in a neighbourhood of zero in the $w$-plane.  The Voiculescu transform is a meromorphic function with a simple
		pole of residue one at $w=0$.  The Voiculescu transform less its principal part, $R_X(w)=V_X(w)-\frac{1}{w}$,
		is an analytic function known as the $R$-transform of $X$.  From the formal identities 
		$V(G(z))=z,\ G(V(w))=w$ and the asymptotics $G_X(z) \sim \frac{1}{z}$ as 
		$|z| \rightarrow \infty$ and $V_X(w) \sim \frac{1}{w}$ as $|w| \rightarrow 0$,
		we expect to find a neighbourhood $\domain{D}_{\infty}$ of infinity in the $z$-plane and a 
		neighbourhood $\domain{D}_0$ of zero in the $w$-plane such that $G_X:\domain{D}_\infty
		\rightarrow \domain{D}_0$ and $V_X:\domain{D}_0 \rightarrow \domain{D}_\infty$ are 
		mutually inverse holomorphic bijections.  The existence of the required domains hinges on identifying 
		regions where the Cauchy and Voiculescu transforms are injective,
		and this can be established through a complex-analytic argument, see \cite[Chapter 4]{MS}.
				
		With these pieces in place, we can state Voiculescu's algorithm for the addition of 
		free random variables. 
						
		\begin{description}
			
				\item[Input] 
					$\mu_X$ and $\mu_Y$.
				
						\bigskip
						
					\item[\emph{Step 1}]
					Compute
					
						\begin{equation*}
							G_X(z) =\int\limits_{\field{R}} \frac{1}{z-t} \mu_X(\mathrm{d}t),\ 
							G_Y(z) =\int\limits_{\field{R}} \frac{1}{z-t} \mu_Y(\mathrm{d}t)
						\end{equation*}
				
					\item[\emph{Step 2}]
					Solve the first Voiculescu functional equations,
					
						\begin{equation*}
							(G_X \circ V_X)(w)=w, \quad (G_Y \circ V_Y)(w)=w
						\end{equation*}
						
					\noindent
					subject to $V_X(w) \sim \frac{1}{w}$ near $w=0$.
								
						\medskip
							
					\item[\emph{Step 3}]
					Remove principal part,
					
						$$R_X(w)=V_X(w)-\frac{1}{w}, \quad R_Y(w) = V_Y(w)-\frac{1}{w},$$
						
					\noindent
					add,
					
						$$R_{X+Y}(w):=R_X(w)+R_Y(w),$$
					
					\noindent
					restore principal part,
					
						$$V_{X+Y}(w):=R_{X+Y}(w)+\frac{1}{w}.$$
						\medskip
				
					\item[\emph{Step 4}]
					Solve the second Voiculescu functional equation, 
					
						$$(V_{X+Y} \circ G_{X+Y})(z)=z,$$
						
					\noindent
					subject to $G_{X+Y}(z) \sim \frac{1}{z}$ near $z=\infty$.
									
						\bigskip
			
				\item[Output] 
				$G_{X+Y}(z)$.
				
			\end{description}
			
		Voiculescu's algorithm is directly analogous to the classical algorithm presented in the previous section.
		The analogy can be succinctly summarized as follows:
		
			\begin{center}
			\boxed{
			\emph{The $R$-transform is the free analogue of the logarithm of the Fourier transform}.}
			\end{center}
		
		In Step One, we try to compute the integral transforms $G_X(z),G_Y(z)$ in terms of elementary functions.
			
		Step Two changes gears from a moment 
		computation to a cumulant computation.  Since free cumulants linearize the addition of 
		free random variables, the new function $V_{X+Y}(w) :=R_X(w) + R_Y(w)+\frac{1}{w}$ defined in Step Three
		encodes the free cumulants of $\kappa_n(X+Y)$ as its Laurent coefficients of non-negative degree.
					
		In Step Four we define a new function $G_{X+Y}(z)$ by 
		solving the second Voiculescu functional equation.  The free exponential formula and the free
		moment-cumulant formula combine
		in the reverse direction to tell us 
		that the Laurent series of $G_{X+Y}(z)$ is
		
			\begin{equation*}
				G_{X+Y}(z) = \sum_{n=0}^{\infty} \frac{m_n(X+Y)}{z^{n+1}}.
			\end{equation*}
		
		An optional Fifth Step is 
		to extract the distribution $\mu_{X+Y}$ from $G_{X+Y}(z)$ using the Stieltjes inversion formula:
		
			\begin{equation*}
				\mu_{X+Y}(\mathrm{d}t) = -\frac{1}{\pi} \lim_{\varepsilon \rightarrow 0}
				\Im G_{X+Y}(t+i\varepsilon).
			\end{equation*}
				
		\subsection{Solution of Kesten's problem}
		Our motivation for building up the additive theory of free random variables came from Kesten's
		problem: explicitly determine the loop generating function of the free group $\group{F}_2$,
		and more generally of the free group $\group{F}_d,\ d \geq 2$. This amounts 
		to computing the moment generating function
		
			\begin{equation*}
				L_d(z) = 1+ \sum_{n=1}^{\infty} m_n(S_d)z^d
			\end{equation*}
			
		\noindent
		of the sum 
		
			\begin{equation*}
				S_d= X_1+ \dots + X_d
			\end{equation*}

		\noindent 
		of fid (free identically distributed) random variables with moments
				
			\begin{equation*}
				\tau[X_i^n]=\begin{cases}
						0, \text{ $n$ odd} \\
						{n \choose n/2}, \text{ $n$ even}.
					\end{cases}
			\end{equation*}
			
		\noindent
		Voiculescu's algorithm gives us the means to obtain this generating function
		provided we can feed it the required input, namely 
		a compactly supported probability measure on $\field{R}$ 
		with moment sequence
		
			\begin{equation*}
				0, {2 \choose 1}, 0, {4 \choose 2}, 0, {6 \choose 3},0,\dots.
			\end{equation*}
			
		\noindent
		As we saw above, the exponential generating function of this moment sequence,
		
			\begin{equation*}
				M(z) = \sum_{k=0}^{\infty} \frac{z^{2k}}{k!k!} = I_0(2z),
			\end{equation*}
			
		\noindent
		coincides with the modified Bessel function of order zero.  From the integral
		representation
		
			\begin{equation*}
				I_0(2z) = \frac{1}{\pi}\int\limits_0^\pi e^{2(\cos \theta)z} \mathrm{d}\theta
			\end{equation*}
			
		\noindent
		we conclude that a random variable $X$ with odd moments zero and even moments
		the central binomial coefficients is given by
		$X=2\cos(Y)$, where $Y$ has uniform distribution over $[0,\pi]$.  Making the 
		same change of variables that we did in Lecture One, we obtain
		
			\begin{equation*}
				M_X(z) = \frac{1}{\pi} \int\limits_{-2}^2 e^{tz} \frac{1}{\sqrt{4-t^2}}\mathrm{d}t,
			\end{equation*}
			
		\noindent
		so that $\mu_X$ is supported on $[-2,2]$ with density
		
			\begin{equation*}
				\mu_X(\mathrm{d}t) = \frac{1}{\pi\sqrt{4-t^2}}\mathrm{d}t.
			\end{equation*}
			
		\noindent 
		This measure is known as the arcsine distribution because its cumulative distribution function is
		
			\begin{equation*}
				\int_{-2}^x \mu_X(\mathrm{d}t) = \frac{1}{2}+\frac{\operatorname{arcsine}(\frac{x}{2})}{\pi}.
			\end{equation*}
			
			\begin{figure}
			\includegraphics{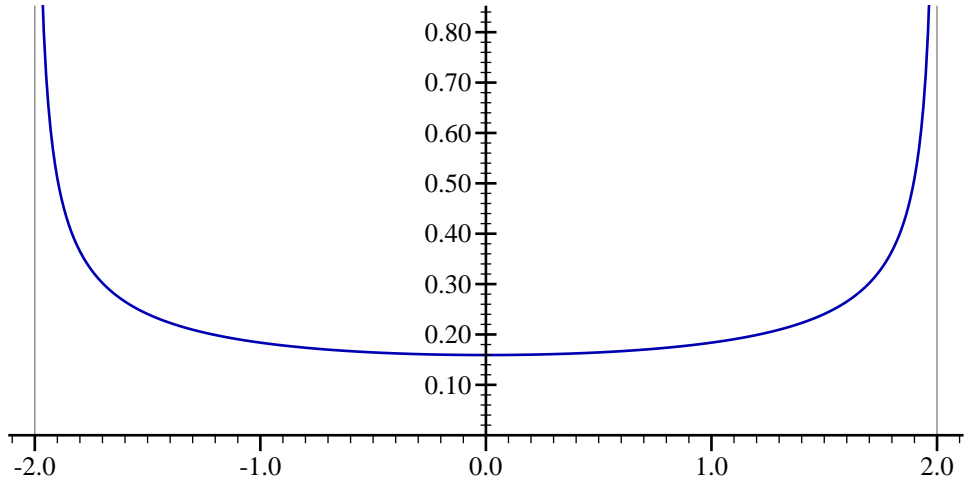}
			\caption{\label{fig:Arcsine}The arcsine density}
			\end{figure}
			
		\noindent
		So to obtain the loop generating function $L_2(z)$ for the simple random walk on $\group{F}_2$, we 
		should run Voiculescu's algorithm with input $\mu_X=\mu_Y=\text{ arcsine}$. 
				
		Let us warm up with an easier computation.  Suppose that $X,Y$ are not fid arcsine random variables,
		but rather fid $\pm1$-Bernoulli random variables:
		
			\begin{equation*}
				\mu_X = \mu_Y = \frac{1}{2}\delta_{-1} + \frac{1}{2}\delta_{+1}.
			\end{equation*}
			
		\noindent
		We will use Voiculescu's algorithm to obtain the distribution of $X+Y$.
		If $X,Y$ were classically iid Bernoullis, we would of course obtain the binomial distribution
		
			\begin{equation*}
				\mu_{X+Y} = \frac{1}{4}\delta_{-2} + \frac{1}{2} \delta_0 + \frac{1}{4}\delta_{+2}
			\end{equation*}
			
		\noindent
		giving the distribution of the simple random walk on $\group{Z}$ at time two.
		The result is quite different in the free case.
		
		Step One.  Obtain the Cauchy transform,
		
			\begin{equation*}
				G_X(z)=G_Y(z)=\frac{1}{2}\bigg{(} \frac{1}{z+1} + \frac{1}{z-1} \bigg{)} = \frac{z}{z^2-1}= \sum_{n=0}^{\infty} \frac{1}{z^{2n+1}}.
			\end{equation*}
						
		Step Two.  Solve the first Voiculescu functional equation.  From Step One, this is
		
			\begin{equation*}
				wV^2(w) -V(w)-w=0,
			\end{equation*}
			
		\noindent
		which has roots
		
			\begin{equation*}
				\frac{1+\sqrt{1+4w^2}}{2w} = \frac{1}{w}+w-w^3+2w^5-5w^7+\dots,\quad
				\frac{1-\sqrt{1+4w^2}}{2w} = -w+w^3-2w^5+\dots.
			\end{equation*}
			
		\noindent
		We identify the first of these as the Voiculescu transform $V_X(w)=V_Y(w)$.
		
		Step Three.  Compute the $R$-transform,
		
			\begin{equation*}
				R_X(w) = R_Y(w)=\frac{1+\sqrt{1+4w^2}}{2w}-\frac{1}{w} = \frac{\sqrt{1+4w^2}-1}{2w},
			\end{equation*}
			
		\noindent
		and sum to obtain
		
			\begin{equation*}
				R_{X+Y}(w) = R_X(w)+R_Y(w) = \frac{\sqrt{1+4w^2}-1}{w}.
			\end{equation*}	
			
		\noindent
		Now restore the principal part,
		
			\begin{equation*}
				V_{X+Y}(w) = R_{X+Y}(w)+\frac{1}{w} = \frac{\sqrt{1+4w^2}}{w}.
			\end{equation*}
			
		Step Four.  Solve the second Voiculescu functional
		equation.  From Step Three, this is the equation
		
			\begin{equation*}
				\frac{\sqrt{1+4G(z)^2}}{G(z)}=z,
			\end{equation*}
			
		\noindent
		which has roots
		
			\begin{equation*}
				\frac{\pm1}{\sqrt{z^2-4}} = \frac{\pm1}{z} + \frac{\pm2}{z^3}+\frac{\pm6}{z^5}+\frac{\pm20}{z^7}+
				\frac{\pm70}{z^9}+\frac{\pm252}{z^{11}}+\dots.
			\end{equation*}
			
		\noindent
		The positive root is identified as $G_{X+Y}(z)$.
		
		Finally, we perform the optional fifth step to recover the distribution $\mu_{X+Y}$
		whose Cauchy transform is $G_{X+Y}(z)$.  This can be done in two ways.  First, we could
		notice that the non-zero Laurent coefficients of $G_{X+Y}$ are the central binomial coefficients 
		${2k \choose k}$, and we just determined that these are the moments
		of the arcsine distribution.  Alternatively we could use Stieltjes inversion:
		
			\begin{equation*}
				\mu_{X+Y}(\mathrm{d}t) = -\frac{1}{\pi} \lim_{\varepsilon \rightarrow 0} \frac{1}{\sqrt{(t+i\varepsilon)^2-4}}
				=  -\frac{1}{\pi}\Im\frac{1}{\sqrt{t^2-4}} =  \frac{1}{\pi\sqrt{4-t^2}}\delta_{|t|\leq2}.
			\end{equation*}
	
		We conclude that the sum of two fid Bernoulli random variables has arcsine distribution.  
		Note the surprising feature that the outcome of a free coin toss
		has continuous distribution over $[-2,2]$.  More generally, we can say that the sum
		
			\begin{equation*}
				S_d = X_1+ \dots + X_{2d}
			\end{equation*}
			
		\noindent
		of $2d$ fid $\pm1$-Bernoulli random variables, i.e. the sum of $2d$ free coin tosses, 
		encodes all information about the simple random walk on $\group{F}_d$ in its moments.
		
		Let us move on to the solution of Kesten's problem for $\group{F}_2$.  Here $X,Y$ are fid 
		arcsine random variables.
				
		Step One.  The Cauchy transform $G_X(z)=G_Y(z)$ is the output of our last application
		of the algorithm, namely 
		
			\begin{equation*}
				G_X(z) = G_Y(z) = \frac{1}{\sqrt{z^2-4}}.
			\end{equation*}
			
		Step Two.  Solve the first Voiculescu functional equation to obtain
		
			\begin{equation*}
				V_X(w)=V_Y(w)=\frac{\sqrt{1+4w^2}}{w}=\frac{1}{w}+2w-2w^3+\dots.
			\end{equation*}
			
		Step Three.  Switch to $R$-transforms, add, switch back to get the Voiculescu
		transform of $X+Y$,
		
			\begin{equation*}
				V_{X+Y}(w) = \frac{2\sqrt{1+4w^2}-1}{z}=\frac{1}{w}+4w-4w^3+ \dots.
			\end{equation*}
			
		Step Four.  Solve the second Voiculescu functional equation to obtain
		
			\begin{equation*}
				G_{X+Y}(z) = \frac{-z+2\sqrt{z^2-12}}{z^2-16}=\frac{1}{z}+\frac{4}{z^3} + \frac{28}{z^5}+\frac{232}{z^7}+\frac{2092}{z^9}+\dots.
			\end{equation*}
			
		We can now calculate the loop generating function for $\group{F}_2$,
		
			\begin{equation*}
				L_2(z)=\frac{1}{z}G_{X+Y}(\frac{1}{z}) = \frac{-1+2\sqrt{1-12z^2}}{1-16z^2}=1+4z^2+28z^4+232z^6+2092z^8+\dots.
			\end{equation*}
			
		\noindent
		More generally, we can run through the above steps for general $d$ to obtain the loop generating function
		
			\begin{equation*}
				L_d(z) = \frac{-(d-1)+d\sqrt{1-4(2d-1)z^2}}{1-16z^2}
			\end{equation*}
			
		\noindent
		for the free group $\group{F}_d$, $d \geq 2$,
		which in turn leads to the probability generating function
		
			\begin{equation*}
				L_d(\frac{z}{2d}) = \frac{-(d-1)+d\sqrt{1-(2d-1)(\frac{z}{d})^2}}{1-4(\frac{z}{d})^2}.
			\end{equation*}
			
		\noindent
		Applying standard methods from analytic combinatorics \cite{FS}, this expression leads to 
		the asymptotics
		
			\begin{equation*}
				\rho_d(n) \sim \text{const}_d \cdot n^{-\frac{3}{2}} \bigg{(} \frac{2\sqrt{d}}{d+1}\bigg{)}^n
			\end{equation*}
			
		\noindent
		for the return probability of the simple random walk on $\group{F}_d,\ d \geq 2$.  From this
		we can conclude that the simple random walk on $\group{F}_d$ is transient for all $d \geq 2$,
		and indeed that $\group{F}_d$ is non-amenable for all $d \geq 2$.
		
		\subsection{Spectral measures and free convolution}
		Voiculescu's algorithm outputs a function $G_{X+Y}(z)$ which encodes the moments of 
		the sum of two freely independent random variables $X$ and $Y$.  As input, it requires a 
		pair of compactly supported real measures $\mu_X,\mu_Y$ which act as distributions
		for $X$ and $Y$ in the sense that
		
			\begin{equation*}
				\tau[X^n] = \int\limits_{\field{R}} t^n \mu_X(\mathrm{d}t), \quad
				\tau[Y^n] = \int\limits_{\field{R}} t^n \mu_Y(\mathrm{d}t).
			\end{equation*}
			
		\noindent
		In our applications of Voiculescu's algorithm we were able to find such measures
		by inspection.  Nevertheless, it is of theoretical and psychological importance to 
		determine sufficient conditions guaranteeing the existence of measures with the
		required properties.
		
		If $X:\Omega \rightarrow \field{C}$ is a random variable defined on a Kolmogorov
		triple $(\Omega, \algebra{F},P)$, its distribution $\mu_X$ is the pushforward 
		of $P$ by $X$,
		
			\begin{equation*}
				\mu_X(B) = (X_*P)(B) = P(X^{-1}(B))
			\end{equation*}
			
		\noindent
		for any Borel (or Lebsegue) set $B \subseteq \field{C}$.  One has the 
		general change of variables formula
		
			\begin{equation*}
				\E[f(X)] = \int\limits_{\field{C}} f(z) \mu_X(\mathrm{d}z)
			\end{equation*}
			
		\noindent
		for any reasonable $f:\field{C} \rightarrow \field{C}$.  If $X$ is essentially 
		bounded and real-valued, $\mu_X$ is compactly supported in $\field{R}$.
		As a random variable $X$ living in an abstract non-commutative probability 
		space $(\algebra{A},\tau)$ is not a function, one must obtain $\mu_X$
		by some other means.  
		
		The existence of distributions is too much to expect within the framework of
		a non-commtative probability space, which is a purely algebraic object.  We 
		need to inject some analytic structure into $(\algebra{A},\tau)$.  
		This is achieved by upgrading $\algebra{A}$ to a $*$-algebra, i.e. a complex
		algebra equipped with a map $*:\algebra{A} \rightarrow \algebra{A}$ satisfying
		
			\begin{equation*}
				(X^*)^* = X, \quad (\alpha X + \beta Y)^* = \overline{\alpha}X^* + 
				\overline{\beta}Y^*, \quad (XY)^*=Y^*X^*.
			\end{equation*}
			
		\noindent
		This map, which is an abstraction of complex conjugation, 
		is required to be compatible with the expectation $\tau$
		in the sense that
		
			\begin{equation*}
				\tau[X^*] = \overline{\tau[X]}.
			\end{equation*}
			
		\noindent
		A non-commutative probability space equipped with this extra structure is 
		called a non-commutative $*$-probability space.
		
		In the framework of a $*$-probability space we can single out a class of 
		random variables analogous to real random variables in classical probability.
		These are the fixed points of $*$, 
		$X^*=X$.  A random variable with this property is called self-adjoint.  
		Self-adjoint random variables have real expected values,
		$\tau[X] = \tau[X^*] = \overline{\tau[X]}$, and more generally 
		$\tau[f(X)] \in \field{R}$ for any polynomial $f$ with real coefficients.
		
		The identification of bounded random variables requires one more upgrade.
		Given a $*$-probability space $(\algebra{A},\tau)$, we can introduce a Hermitian
		form $B:\algebra{A} \times \algebra{A} \rightarrow \field{C}$ defined by
		
			\begin{equation*}
				B(X,Y) = \tau[XY^*].
			\end{equation*}
			
		\noindent
		If we require that $\tau$ has the positivity property $\tau[XX^*] \geq 0$ for 
		all $X \in \algebra{A}$, then we obtain a semi-norm 
		
			\begin{equation*}
				\|X\| = B(X,X)^{1/2}
			\end{equation*}
			
		\noindent
		on $\algebra{A}$, and we can access the Cauchy-Schwarz inequality 
		
			\begin{equation*}
				|B(X,Y)| \leq \|X\| \|Y\|.
			\end{equation*}
			
		\noindent
		Once we have Cauchy-Schwarz, we can prove the monotonicity inequalities
		
			\begin{equation*}
				\begin{split}
				&|\tau[X]| \leq |\tau[X^2]|^{1/2} \leq |\tau[X^4]|^{1/4} \\
				&|\tau[X^3]| \leq |\tau[X^4]|^{1/4} \leq |\tau[X^6]|^{1/6} \\
				&|\tau[X^5]| \leq |\tau[X^6]|^{1/6} \leq |\tau[X^8]|^{1/8} \\
				&\vdots
				\end{split}
			\end{equation*}
			
		\noindent
		from which the chain of inequalities 
		
			\begin{equation*}
				|\tau[X]| \leq |\tau[X^2]|^{1/2} \leq |\tau[X^4]|^{1/4} \leq |\tau[X^6]|^{1/6} \leq |\tau[X^8]|^{1/8} \leq \dots
			\end{equation*}
			
		\noindent
		can be extracted.  From this we conclude that the limit
		
			\begin{equation*}
				\rho(X) := \lim_{k \rightarrow \infty} |\tau[X^{2k}]|^{1/(2k)}
			\end{equation*}
			
		\noindent
		exists in $\field{R}_{\geq 0} \cup \{\infty\}$.  This limit is called the spectral radius of $X$.
		A random variable $X \in \algebra{A}$ is said to be bounded if its spectral radius is
		finite, $\rho(X) < \infty$.  
		
		In the framework of a non-commutative $*$-probability space $(\algebra{A},\tau)$
		with non-negative expectation, bounded self-adjoint random variables
		play the role 
		of essentially bounded real-valued random variables in classical probability theory. 
		With some work, one may deduce from the Riesz representation theorem that to each
		bounded self-adjoint $X$ corresponds a unique Borel measure 
		$\mu_X$ supported in $[-\rho(X),\rho(X)]$
		such that
		
			\begin{equation*}
				\tau[f(X)] = \int\limits_{\field{R}} f(t) \mu_X(\mathrm{d}t)
			\end{equation*}
			
		\noindent
		for all polynomial functions $f:\field{C} \rightarrow \field{C}$.  The details of this argument,
		in which a reverse-engineered Cauchy transform plays the key role,
		are given in Tao's notes \cite{Tao}.  The measure $\mu_X$ is often called the spectral measure
		of $X$, but we will refer to it as the distribution of $X$.  There is also a converse to 
		this result: given any compactly supported measure $\mu$ on $\field{R}$, there exists
		a bounded self-adjoint random variable $X$ living in some non-commutative 
		$*$-probability space $(\algebra{A},\tau)$ whose distribution is $\mu$.
		Consequently, given two compactly supported real probability measures $\mu,\nu$ we 
		may define a new measure $\mu \boxplus \nu$ as ``the distribution of the 
		random variable $X+Y$, where $X$ and $Y$ are freely independent
		bounded self-adjoint random variables with 
		distributions $\mu$ and $\nu$, respectively.''  Since the sum of two bounded self-adjoint
		random variables is again bounded self-adjoint, $\mu \boxplus \nu$ is another compactly
		supported real probability measure.  Moreover, $\mu \boxplus \nu$ does not depend 
		on the particular random variables chosen to realize $\mu$ and $\nu$.  Thus we get 
		a bona fide binary operation $\boxplus$ on the set of compactly supported real measures,
		which is known as the additive free convolution.  For example, we computed above that
		
			\begin{equation*}
				\text{Bernoulli} \boxplus \text{Bernoulli} = \text{Arcsine}.
			\end{equation*}
			
		The additive free convolution of measures is induced by the addition of free random variables.
		As such, it is the free analogue of the classical convolution of measures induced by the addition
		of classically independent random variables.  Like classical convolution, free convolution
		can be defined for unbounded measures, but this requires more work \cite{BV}.
																																
		\subsection{Free Poisson limit theorem}
		Select positive real numbers $\lambda$ and $\alpha$.
		Consider the measure 
		
			\begin{equation*}
				\mu_N = (1-\frac{\lambda}{N})\delta_0 + \frac{\lambda}{N}\delta_{\alpha}
			\end{equation*}
			
		\noindent
		which consists of an atom of mass $1-\frac{\lambda}{N}$ placed at zero and an atom
		of mass $\frac{\lambda}{N}$ placed at $\alpha$.  For $N$ sufficiently large,
		$\mu_N$ is a probability measure.  Its moment sequence is 
		
			\begin{equation*}
				m_n(\mu_N) = \frac{\lambda}{N}\alpha^n, \quad n \geq 1.
			\end{equation*}
		
		\noindent
		The $N$-fold classical convolution of $\mu_N$ with itself,
		
			\begin{equation*}
				\mu_N^{*N} = \underbrace{\mu_N * \dots * \mu_N}_N,
			\end{equation*}
			
		\noindent
		converges weakly to the Poisson measure of rate $\lambda$ and jump
		size $\alpha$ as $N \rightarrow \infty$.  This is a classical limit theorem 
		in probability known as the Poisson Limit Theorem, or the Law of Rare Events.
		
		Let us obtain a free analogue of the Poisson Limit Theorem.  This should be a limit
		law for the iterated free convolution 
		
			\begin{equation*}
				\mu_N^{\boxplus N} = \underbrace{\mu_N \boxplus \dots \boxplus \mu_N}_N.
			\end{equation*}
		
		From the free 
		moment-cumulant formula, we obtain the estimate
		
			\begin{equation*}
				\kappa_n(\mu_N) = m_n(\mu_N) + O\bigg{(} \frac{1}{N^2} \bigg{)} = \frac{\lambda}{N}\alpha^n
				+O\bigg{(} \frac{1}{N^2} \bigg{)}.
			\end{equation*}
			
		\noindent
		Since free cumulants linearize free convolution, we have
			
			\begin{equation*}
				\kappa_n(\mu_N^{\boxplus N}) =N \kappa_n(\mu_N) = \lambda\alpha^n +O\bigg{(} \frac{1}{N} \bigg{)} .
			\end{equation*}
			
		\noindent
		Thus 
		
			\begin{equation*}
				\lim_{N \rightarrow \infty} \kappa_n(\mu_N) = \lambda\alpha^n,
			\end{equation*}
			
		\noindent
		and it remains to determine the measure $\mu$ with this free cumulant sequence.
		The Voiculescu transform of $\mu$ is
		
			\begin{equation*}
				V_\mu(w) = \frac{1}{w} + \sum_{n=0}^{\infty} \lambda\alpha^{n+1}w^n = 
				\frac{1}{w} + \frac{\lambda\alpha}{1-\alpha w},
			\end{equation*}
			
		\noindent
		so the second Voiculescu functional equation $V_\mu(G_\mu(z))=z$ yields
		
			\begin{equation*}
				\frac{1}{G_\mu(z)} + \frac{\lambda\alpha}{1-\alpha G_\mu(z)}=z.
			\end{equation*}
			
		\noindent
		This equation has two solutions, and the one which behaves like $1/z$ for $|z| \rightarrow \infty$
		is the Cauchy transform of $\mu$.  We obtain
		
			\begin{equation*}
				G_\mu(z) = \frac{z+\alpha(1-\lambda) - \sqrt{(z-\alpha(1+\lambda))^2-4\lambda\alpha^2}}{2\alpha z}.
			\end{equation*}
			
		\noindent
		Applying Stieltjes inversion, we find that the density of $\mu$ is given by 
		
			\begin{equation*}
				\mu(\mathrm{d}t) = \begin{cases}
					(1-\lambda)\delta_0 + \lambda m(t)\mathrm{d}t,\ 0 \leq \lambda \leq 1\\
					m(t)\mathrm{d}t,\ \lambda > 1
				\end{cases}
			\end{equation*}
			
		\noindent
		where 
		
			\begin{equation*}
				m(t) = \frac{1}{2\pi\alpha t}\sqrt{4\lambda\alpha^2 - (t-\alpha(1+\lambda))^2}.
			\end{equation*}
			
		\noindent
		This measure is known as the Marchenko-Pastur distribution after the Ukrainian mathematical 
		physicists Vladimir Marchenko and Leonid Pastur, who discovered it in their
		study of the asymptotic eigenvalue distribution of a certain class of random matrices.
		
		\subsection{Semicircle flow}
		Given $r>0$, let $\mu_r$ be the semicircular measure of radius $r$,
		
			\begin{equation*}
				\mu_r(\mathrm{d}t) = \frac{2}{\pi r^2}\sqrt{r^2-t^2} \mathrm{d}t.
			\end{equation*}
			
		\noindent
		Taking $r=2$ yields the standard semicircular distribution.  
		Let $\mu$ be an arbitrary compactly supported probability measure on $\field{R}$.
		The function
		
			\begin{equation*}
				f_\mu: \{\text{positive real numbers} \} \rightarrow \{\text{compactly supported real 
				measures}\}
			\end{equation*}
			
		\noindent
		defined by
		
			\begin{equation*}
				f_\mu(r) = \mu \boxplus \mu_r
			\end{equation*}
			
		\noindent
		is called the semicircle flow.  The semicircle flow has very interesting dynamics:
		in one of his earliest articles on free random variables \cite{Voiculescu},
		Voiculescu showed that 
		
			\begin{equation*}
				\frac{\partial G(r,z)}{\partial r} + G(r,z) \frac{\partial G(r,z)}{\partial z} = 0,
			\end{equation*}
			
		\noindent
		where $G(r,z)$ is the Cauchy transform of $f_\mu(r)=\mu \boxplus \mu_r$.  Thus the free analogue 
		of the heat equation is the complex inviscid Burgers equation.  For a detailed 
		analysis of the semicircle flow, see \cite{Biane:semicircle}.
													
	\section{Lecture Three: Modelling the Free World}
		Free random variables are of interest for many reasons.  First and foremost,
		Voiculescu's free probability 
		theory is an intrinsically appealing subject
		worthy of study from a purely esthetic point of view.  Adding to 
		this are the many remarkable connections between free
		probability and other parts of mathematics, including operator
		algebras, representation theory, and random matrix theory.
		This lecture is an exposition of Voiculescu's 
		discovery that random matrices provide asymptotic models
		of free random variables.  We follow the treatment of 
		Nica and Speicher \cite{NS}.  
					
		\subsection{Algebraic model of a free arcsine pair}
		In Lecture Two we gave a group-theoretic construction of a pair of free random variables
		each of which has an arcsine distribution.  In this example, the algebra of random 
		variables is the group algebra $\algebra{A}=\algebra{A}[\group{F}_2]$ of the free group 
		on two generators $A,B$, and the expectation $\tau$ is the coefficient-of-$\id$ functional.  
		We saw that the random variables
		
			\begin{equation*}
				X=A+A^{-1}, \quad Y=B+B^{-1}
			\end{equation*}
			
		\noindent
		are freely independent, and each has an arcsine distribution:
		
			\begin{equation*}
				\tau[X^n]=\tau[Y^n] = \begin{cases}
					0, \text{ if $n$ odd}\\
					{n \choose \frac{n}{2}}, \text{ if $n$ even}
					\end{cases}.
			\end{equation*}
		
		\subsection{Algebraic model of a free semicircular pair}
		We can give a linear-algebraic model of a pair of free random variables each of which 
		has a semicircular distribution.  
		The ingredients in this construction
		are a complex vector space $\vecspace{V}$ and an inner product 
		$B:\vecspace{V} \times \vecspace{V} \rightarrow \field{C}$.  Our random variables will be 
		endomorphisms of the tensor algebra over $\vecspace{V}$,
				
					\begin{equation*}
				\mathfrak{F}(\vecspace{V}) = \bigoplus_{n=0}^{\infty} \vecspace{V}^{\otimes n},
			\end{equation*}
		
		\noindent			
		which physicists and operator algebraists call the full Fock space over $\vecspace{V}$ after the 
		Russian physicist Vladimir Fock.  We view
		the zeroth tensor power $V^{\otimes 0}$ as the line in $\vecspace{V}$ spanned by a distinguished unit 
		vector $\vector{v}_{\emptyset}$ called the vacuum vector.  Let $\algebra{A}=\End \F(V)$.
		This is a unital algebra, with unit the identity operator $I:\F(V) \rightarrow \F(V)$.  To make $\algebra{A}$
		into a non-commutative probability space we need an expectation. 
		We get an expectation by lifting the inner product on $\vecspace{V}$ to the inner product
		$\mathfrak{F}(B):\mathfrak{F}(\vecspace{V}) \times \mathfrak{F}(\vecspace{V}) \rightarrow \field{C}$ defined by
		
			\begin{equation*}
				\mathfrak{F}(B)(\vector{v}_1 \otimes \dots \otimes \vector{v}_m, \vector{w}_1 \otimes \dots \otimes \vector{w}_n)
				= \delta_{mn} B(\vector{v}_1,\vector{w}_1) \dots B(\vector{v}_n,\vector{w}_n).
			\end{equation*}
			
		\noindent 
		Note that this inner product makes $\algebra{A}=\End\F(B)$ into a $*$-algebra: for
		each $X \in \algebra{A}$, $X^*$ is that linear operator for which the equation
		
			\begin{equation*}
				\F(B)(X\vector{s},\vector{t})
				=\F(B)(\vector{s},X^*\vector{t})
			\end{equation*}
			
		\noindent 
		holds true for every pair of tensors $\vector{s},\vector{t} \in \F(\vecspace{V})$.  The expectation on $\algebra{A}$
		is the linear functional $\tau:\algebra{A} \rightarrow \field{C}$ defined by 
		
			\begin{equation*}
				\tau[X] = \F(B)(X\vac,\vac).
			\end{equation*}
			
		\noindent
		This functional is called vacuum expectation.  It is unital because
		
			\begin{equation*}
				\tau[I] = \F(B)(I\vac,\vac)=B(\vac,\vac)=1.
			\end{equation*}
			
		\noindent
		Thus $(\algebra{A},\tau)$ is a non-commutative $*$-probability space.
				
		To construct a semicircular element in $(\algebra{A},\tau)$, notice that to
		every non-zero vector $\vector{v} \in \vecspace{V}$ is naturally associated a pair of linear operators
		$R_\vector{v},L_{\vector{v}}:\F(V) \rightarrow \F(V)$ whose action on decomposable tensors
		is defined by tensoring,
		
			\begin{equation*}
				\begin{split}
				R_\vector{v}(\vector{v}_{\emptyset})&=\vector{v} \\
				R_\vector{v}(\vector{v}_1 \otimes \dots \otimes \vector{v}_n) &= \vector{v} \otimes \vector{v}_1 \otimes \dots \otimes \vector{v}_n,\quad
				\quad n \geq 1,
				\end{split}
			\end{equation*}
			
		\noindent
		and insertion-contraction
		
			\begin{equation*}
				\begin{split}
					L_\vector{v}(\vac) &= \vector{0} \\
					L_\vector{v}(\vector{v}_1) &= B(\vector{v}_1,\vector{v})\vac \\
					L_\vector{v}(\vector{v}_1 \otimes \vector{v}_2 \otimes \dots \otimes \vector{v}_n) 
					&= B(\vector{v}_1,\vector{v}) \vector{v}_2 \otimes \dots \otimes \vector{v}_n,
					\quad n \geq 2.
				\end{split}
			\end{equation*}
			
		\noindent
		Since $R_\vector{v}$ maps $\vecspace{V}^{\otimes n} \rightarrow \vecspace{V}^{\otimes n+1}$ for each $n \geq 0$,
		it is called the raising (or creation) operator associated to $\vector{v}$.  Since  
		$L_\vector{v}$ maps $\vecspace{V}^{\otimes n} \rightarrow \vecspace{V}^{\otimes n-1}$ for each $n \geq 1$
		and kills the vacuum, it 
		is called the lowering (or annihilation) operator associated to $\vector{v}$.
		We have $R_\vector{v}^*=L_\vector{v}$,
		and also
		
			\begin{equation*}
				L_\vector{v}R_\vector{w} = B(\vector{w},\vector{v})I
			\end{equation*}
			
		\noindent
		for any vectors $\vector{v},\vector{w} \in \vecspace{V}$.
		
		Let $\vector{v} \in \vecspace{V}$ be a unit vector, $B(\vector{v},\vector{v})=1$, and consider the self-adjoint random
		variable 
		
			\begin{equation*}
				X_\vector{v}=L_\vector{v} + R_\vector{v}.
			\end{equation*}		
			
		\noindent
		We claim that $X_\vector{v}$ has a semicircular distribution:
		
			\begin{equation*}
				m_n(X_\vector{v}) =\tau[X_\vector{v}^n]= \begin{cases}
					0, \text{ if $n$ odd}\\
					\Cat_{\frac{n}{2}}, \text{ if $n$ even}
				\end{cases}.
			\end{equation*}
			
		\noindent
		To see this, we expand
		
			\begin{equation*}
				\tau[X_\vector{v}^n] = \tau[(X_\vector{v}+Y_\vector{v})^n] = \sum_{W \in \{L_\vector{v},R_\vector{v}\}^n}\tau[W],
			\end{equation*}	
			
		\noindent
		where the summation is over all words of length $n$ in the operators $L_\vector{v},R_\vector{v}$.
		Only a very small fraction of these words have non-zero vacuum expectation.
		Using the relation $L_\vector{v}R_\vector{v}=I$ to remove occurrences of the substring $L_\vector{v}R_\vector{v}$,
		we see that any such word can be placed in normally ordered form
		
			\begin{equation*}
				W=\underbrace{R_\vector{v} \dots R_\vector{v}}_a\underbrace{L_\vector{v} \dots L_\vector{v}}_b
			\end{equation*}
			
		\noindent
		with $a+b \leq n$.  Since the lowering operator kills the vacuum vector, the vacuum expectation of $W$ 
		can only be non-zero if $b=0$.  On the other hand, since $\vecspace{V}^{\otimes a}$ is
		$\F(B)$-orthogonal to $\vecspace{V}^{\otimes 0}$ for $a>0$, we must also have $a=0$
		to obtain a non-zero contribution.  Thus the only words which contribute to the above sum are those
		whose normally ordered form is that of the identity operator.  If we replace each occurrence of 
		$L_\vector{v}$ in $W$ with a $+1$ and each occurrence of $R_\vector{v}$ in $W$ with a $-1$, the
		condition that $W$ reduces to $I$ becomes the condition that the corresponding bitstring
		has total sum zero and all partial sums non-negative.  There are no such bitstrings for $n$ odd,
		and as we saw in Lecture One when $n$ is even the required bitstrings are counted by the 
		Catalan number $\Cat_{n/2}$.
			
		Now let $\vecspace{V}_1$ and $\vecspace{V}_2$ be $B$-orthogonal vector subspaces of 
		$\vecspace{V}$, each of dimension at least one, and choose unit vectors $\vector{x} \in \vecspace{V}_1,
		\vector{y} \in \vecspace{V}_2$.  According to the above construction, the random variables
		
			\begin{equation*}
				X = L_{\vector{x}}+R_{\vector{x}}, \quad Y = L_{\vector{y}}+R_{\vector{y}}
			\end{equation*}
			
		\noindent
		are semicircular.  In fact, they are freely independent.  To prove this, we must demonstrate 
		that 
		
			\begin{equation*}
				\tau[f_1(X)g_1(Y) \dots f_k(X)g_k(Y)]=0
			\end{equation*}
			
		\noindent
		whenever $f_1,g_1,\dots,f_k,g_k$ are polynomials such that
		
			\begin{equation*}
				\tau[f_1(X)]=\tau[g_1(Y)] = \dots =\tau[f_k(X)]=\tau[g_k(Y)]=0.
			\end{equation*}
			
		\noindent
		This hypothesis means that $f_i(X) = f_i(L_\vector{x}+R_\vector{x})$ is a polynomial
		in $L_\vector{x},R_\vector{x}$ none of whose terms are words which normally order to $I$, 
		and similarly $g_j(Y) = g_j(L_\vector{y}+R_\vector{y})$ is a polynomial
		in $L_\vector{y},R_\vector{y}$ none of whose terms are words which normally order to $I$.
		Consequently, the alternating product
		
			\begin{equation*}
				f_1(X)g_1(Y) \dots f_k(X)g_k(Y)
			\end{equation*}
			
		\noindent
		is a polynomial in the operators $L_\vector{x},R_\vector{x},L_\vector{y},R_\vector{y}$ whose
		terms are words $W$ of the form 
		
			\begin{equation*}
				W_\vector{x}^1W_\vector{y}^1 \dots W_\vector{x}^kW_\vector{y}^k
			\end{equation*}
			
		\noindent
		with $W_\vector{x}^i$ a word in $L_\vector{x},R_\vector{x}$ which does not normally order to
		$I$ and $W_\vector{y}^j$ a word in $L_\vector{y},R_\vector{y}$ which does not normally
		order to $I$.  Thus the only way that $W$ can have a non-zero vacuum expectation is 
		if we can use the relations $L_\vector{x}R_\vector{y} = B(\vector{y},\vector{x})I$ and $L_\vector{y}R_\vector{x}
		=B(\vector{x},\vector{y})I$ to normally order $W$ as
		
			\begin{equation*}
				B(\vector{x},\vector{y})^mB(\vector{y},\vector{x})^nI
			\end{equation*}
			
		\noindent
		with $m,n$ non-negative integers at least one of which is positive.  
		But, since $\vector{x},\vector{y}$ are $B$-orthogonal, this 
		is the zero element of $\algebra{A}$, which has vacuum expectation zero.
						
		\subsection{Algebraic versus asymptotic models}
		We have constructed algebraic models for a free arcsine pair and a free semicircular pair.  Perhaps
		these should be called examples rather than models, since the term model connotes some
		degree of imprecision or ambiguity and algebra is a subject which allows neither.  
		
		Suppose that $X,Y$ are free random variables living in an abstract non-commutative probability space
		$(\algebra{A},\tau)$.  An approximate model for this pair will consist of a sequence 
		$(\algebra{A}_N,\tau_N)$ of concrete or canonical non-commutative probability spaces
		together with a sequence of pairs $X_N,Y_N$ of random variables from these spaces
		such that $X_N$ models $X$ and $Y_N$ models $Y$, i.e.
		
			\begin{equation*}
				\tau[f(X)]= \lim_{N \rightarrow \infty} \tau[f(X_N)], \quad \tau[g(Y)]= \lim_{N \rightarrow \infty} \tau[g(Y_N)]
			\end{equation*}
			
		\noindent
		for any polynomials $f,g$, and such that free independence holds in the large $N$ limit, i.e.
		
			\begin{equation*}
				\lim_{N \rightarrow \infty} \tau[f_1(X_N)g_1(Y_N) \dots f_k(X_N)g_k(Y_N)]=0
			\end{equation*}
			
		\noindent
		whenever $f_1,g_1,\dots,f_k,g_k$ are polynomials such that
		
			\begin{equation*}
				\lim_{N \rightarrow \infty} \tau_N[f_1(X_N)] = \lim_{N \rightarrow \infty} \tau_N[g_1(Y_N)] = \dots
				=\lim_{N \rightarrow \infty} \tau_N[f_k(X_N)] = \lim_{N \rightarrow \infty} \tau_N[g_k(Y_N)] = 0.
			\end{equation*}
			
		The question of which non-commutative probability spaces are considered concrete or canonical, 
		and could therefore serve as potential models, is subjective and determined by individual experience.  
		Three examples of concrete non-commutative probability spaces are:
		
			\begin{description}
			
				\item[\bf Group probability spaces]
				$(\algebra{A},\tau)$ consists of the group algebra $\algebra{A}=\algebra{A}[\group{G}]$ of a group $\group{G}$,
				and $\tau$ is the coefficient-of-identity expectation.  This non-commutative probability space is
				commutative if and only if $\group{G}$ is abelian.
				
					\smallskip
					
				\item[\bf Classical probability spaces]
				$(\algebra{A},\tau)$ consists of the algebra of complex random variables 
				$\algebra{A}=L^{\infty-}(\Omega,\algebra{F},P)=\bigcap_{p=1}^{\infty} L^p(\Omega,\algebra{F},P)$ 
				defined on a Kolmogorov  triple which have finite absolute moments of all orders, and 
				$\tau$ is the classical expectation $\tau[X]=\E[X]$.  Classical probability spaces are always commutative.
				
					\smallskip
					
				\item[\bf Matrix probability spaces]
				$(\algebra{A},\tau)$ consists of the algebra $\algebra{A}=\algebra{M}_N(\field{C})$ of $N \times N$
				complex matrices $X=[X(ij)]$, and expectation is the normalized trace, $\tau[X]=\tr_N[X]=
				\frac{X(11)+\dots+X(NN)}{N}$.  This non-commutative probability space is commutative if and only if 
				$N=1$.
				
			\end{description}
		
		The first class of model non-commutative probability spaces, group probability spaces, is algebraic 
		and we are trying to move 
		away from algebraic examples.  The second model class,
		classical probability spaces, has
		genuine randomness but is commutative.  The third model class,
		matrix probability spaces, has a parameter $N$ that can be pushed to infinity but has no randomness.
		By combining classical probability spaces and matrix probability spaces we arrive at a class of model
		non-commutative probability spaces which incorporate both randomness and a parameter which 
		can be made large.  Thus we are led to consider random matrices.  
				
		The space of $N \times N$ complex random matrices is the non-commutative probability space
		$(\algebra{A}_N,\tau_N) = (L^{\infty-}(\Omega,\algebra{F},P) \otimes \algebra{M}_N(\field{C}),
		\E \otimes \tr_N)$.  A random variable $X_N$ in this space may be viewed as an $N \times N$
		matrix whose entries $X_N(ij)$ belong to the algebra $L^{\infty-}(\Omega,\algebra{F},P)$.  The expectation
		$\tau_N[X_N]$ is the expected value of the normalized trace:
		
			\begin{equation*}
				\tau_N[X_N] = (\E \otimes \tr_N)[X_N] = \E \bigg{[} 
				\frac{X_N(11)+\dots +X_N(NN)}{N} \bigg{]}.
			\end{equation*}
		
		We have already seen indications of a connection between free probability 
		and random matrices.  The fact that Wigner's semicircle law assumes the role of the Gaussian distribution in free 
		probability signals a connection between these subjects.  Another example is the occurrence 
		of the Marchenko-Pastur distribution in the free version of the Poisson limit theorem --- this distribution
		is well-known in random matrix theory in connection with the asymptotic eigenvalue distribution of Wishart matrices.  
		In Lecture One,
		we were led to free independence when we tried to solve a counting problem associated to 
		graphs drawn in the plane.  The use of random matrices to enumerate planar graphs has been a subject
		of much interest in mathematical physics since the seminal work of Edouard Br\'ezin, Claude Itzykson, 
		Giorgio Parisi and Jean-Bernard Zuber \cite{BIPZ}, which built on insights of Gerardus 't Hooft.
		Then, when we examined the dynamics of the semicircle flow, we found that the free analogue of the 
		heat equation is the complex Burgers equation.  This partial differential equation actually appeared in 
		Voiculescu's work \cite{Voiculescu} before it emerged in random matrix theory \cite{Matytsin} and the discrete 
		analogue of random matrix theory, the dimer model \cite{KO}.

		In the remainder of these notes, we will model a pair of free random variables $X,Y$ living in an
		abstract non-commutative probability space using sequences $X_N,Y_N$ of random matrices
		living in random matrix space.  This is first carried out in the important
		special case where $X,Y$ are semicircular random variables, then adapted to allow $Y$ to have
		arbitrary distribution while $X$ remains semicircular, and finally relaxed to allow $X,Y$
		to have arbitrary specified distributions.		
		The random matrix models of free random variables which we describe below were used by Voiculescu in order to 
		resolve several previously intractable problems in the theory of von Neumann algebras, see \cite{MS,VDN} for more 
		information.  Random matrix models which approximate free random variables in a stronger sense than that
		described here were 
		subsequently used by Uffe Haagerup and Steen Thorbj{\o}rnsen \cite{HT} to resolve another operator algebras
		conjecture, this time concerning the Ext-invariant of the reduced $C^*$-algebra of $\group{F}_2$.  
		An important feature of the connection between
		free probability and random matrices is that it can sometimes be inverted to obtain information about random
		matrices using the free calculus.  For each of the three matrix models constructed we give
		an example of this type.
				
		\subsection{Random matrix model of a free semicircular pair}
		In this subsection we construct a random matrix model for a free semicircular pair $X,Y$.  
		
		In Lecture One, we briefly discussed Wigner matrices.  A real Wigner matrix is a symmetric matrix
		whose entries are centred real random variables which are independent up to the symmetry constraint.
		A complex Wigner matrix is a Hermitian matrix whose entries are centred complex random variables
		which are independent up to the complex symmetry constraint.
		Our matrix model for a free semicircular pair will be built out of complex Wigner matrices
		of a very special type: they will be GUE random matrices.  
		
		To construct a GUE random matrix $X_N$, we start with a Ginibre matrix $Z_N$.  
		Let $(\Omega,\algebra{F},P)$ be a Kolmogorov triple.
		The $N^2$ matrix elements $Z_N(ij) \in L^{\infty-}(\Omega,\algebra{F},P)$ of a Ginibre matrix
		are iid complex Gaussian random variables of mean zero and variance $1/N$.
		Thus $Z_N$ is a random variable in the non-commutative probability space
		$(\algebra{A}_N,\tau_N) = (L^{\infty-}(\Omega,\algebra{F},P) \otimes \algebra{M}_N(\field{C}),
		\E \otimes \tr_N)$.  The symmetrized random matrix $X_N=\frac{1}{2}(Z_N+Z_N^*)$ is again
		a member of random matrix space.  The joint distribution of the 
		eigenvalues of $X_N$ can be explicitly computed, and is given by 
		
			\begin{equation*}
				P(\lambda_N(1) \in I_1,\dots,\lambda_N(N) \in I_N) \propto \int\limits_{I_1} \dots \int\limits_{I_N}
				e^{-N^2 \mathcal{H}(\lambda_1,\dots,\lambda_N)} \mathrm{d}\lambda_1 \dots \mathrm{d}\lambda_N
			\end{equation*}	
			
		\noindent
		for any intervals $I_1,\dots,I_N \subseteq \field{R}$, where $\mathcal{H}$ is the log-gas Hamiltonian \cite{Forrester}
		
			\begin{equation*}
				\mathcal{H}(\lambda_1,\dots,\lambda_N) = \frac{1}{N} \sum_{i=1}^N \frac{\lambda_i^2}{2} - \frac{1}{N^2}
				\sum_{1 \leq i \neq j \leq N} \log |\lambda_i-\lambda_j|.
			\end{equation*}
			
		\noindent
		The random point process on the real line driven by this Hamiltonian is known as the Gaussian Unitary Ensemble,
		 and $X_N$ is termed a GUE random matrix.  GUE random matrices sit at the 
		 nexus of the two principal strains of complex random matrix theory: they are simultaneously Hermitian Wigner matrices
		 and unitarily invariant matrices.  The latter condition means that the distribution of a GUE matrix 
		 in the space of $N \times N$ Hermitian matrices is invariant under conjugation by unitary matrices.
		 The spectral statistics of a GUE random matrix can be computed in gory detail from knowledge
		 of the joint distribution of eigenvalues, and virtually any question can be answered.
		The universality programme in random matrix theory seeks to show that, in the limit 
		$N \rightarrow \infty$ and under mild hypotheses, Hermitian Wigner matrices as well as 
		unitarily invariant Hermitian matrices exhibit the same spectral statistics as GUE matrices.
		 		 
		 Given the central role of the GUE in random matrix theory, it is fitting that our matrix model
		 for a free semicircular pair is built from a pair of independent GUE matrices.  
		 The first step in proving this is to show that a single GUE matrix $X_N$ in random matrix 
		 space $(\algebra{A}_N,\tau_N)$ is an asymptotic model for a single semicircular random variable
		 $X$ living in an abstract non-commutative probability space $(\algebra{A},\tau)$.
		 In other words, we need to prove that
		 
		 	\begin{equation*}
				\lim_{N \rightarrow \infty} \tau_N[X_N^n]=\lim_{N \rightarrow \infty} (\E \otimes \tr_N)[X_N^n] = \begin{cases}
					0, \text{ if $n$ odd}\\
					\Cat_{\frac{n}{2}}, \text{ if $n$ even}
				\end{cases}.
			\end{equation*}
		
		\noindent
		In order to establish this, we will not need
		 access to the eigenvalues of $X_N$.  Rather, we work with the correlation functions of its
		 entries.  
		 
		 Let $X_N = [X_N(ij)]$ be a GUE random matrix.  Mixed moments of the random variables 
		 $X_N(ij)$, i.e. expectations of the form
		 
		 	\begin{equation*}
				\E\bigg{[} \prod_{k=1}^n X_N(i(k)j(k)) \bigg{]}
			\end{equation*}
			
		\noindent
		where $i,j$ are functions $[n] \rightarrow [N]$, are called correlation functions.  
		All correlations may be computed in terms of pair correlations (i.e. covariances)
		
		 	\begin{equation*}
				\E[X_N(ij)\overline{X_N(kl)}] = \E[X_N(ij)X_N(lk)] = \frac{\delta_{ik}\delta_{jl}}{N}
			\end{equation*}
		
		\noindent
		using a convenient combinatorial formula known as Wick's formula.	
		This formula, named for the Italian physicist Gian-Carlo Wick, is yet another manifestation of the
		moment-cumulant/exponential formulas.  It asserts that
		
			\begin{equation*}
				\E \bigg{[} \prod_{k=1}^n X_N(i(k)j(k))) \bigg{]} = \sum_{\pi \in \lattice{P}_2(n)}
				\prod_{\{r,s\} \in \pi} \E[X_N(i(r)j(r)) X_N(i(s)j(s))]
			\end{equation*}
		
		\noindent
		for any integer $n \geq 1$ and functions $i,j:[n] \rightarrow [N]$.  The sum on the right hand side is taken over
		all pair partitions of $[n]$, and the product is over the blocks of $\pi$.  For example,
		
			\begin{equation*}
				\E[X_N(i(1)j(1))X_N(i(2)j(2))X_N(i(3)j(3))] =0
			\end{equation*}
			
		\noindent
		since there are no pairings on three points, whereas
		
			\begin{equation*}
				\begin{split}
					&\E[X_N(i(1)j(1))X_N(i(2)j(2))X_N(i(3)j(3))X_N(i(4)j(4))] \\
					=& \E[X_N(i(1)j(1))X_N(i(2)j(2))]\E[X_N(i(3)j(3))X_N(i(4)j(4))] \\
					+& \E[X_N(i(1)j(1))X_N(i(3)j(3))]\E[X_N(i(2)j(2))X_N(i(4)j(4))] \\
					+& \E[X_N(i(1)j(1))X_N(i(4)j(4))]\E[X_N(i(2)j(2))X_N(i(3)j(3))], 
				\end{split}
			\end{equation*}
			
		\noindent
		corresponding to the three pair partitions $\{1,2\} \sqcup \{3,4\},\{1,3\} \sqcup \{2,4\}, \{1,4\} \sqcup \{2,3\}$
		of $[4]$.
		The Wick formula is a special feature of Gaussian random variables which,
		ultimately, is a consequence of the moment formula
		
			\begin{equation*}
				\E[X^n] = \sum_{\pi \in \lattice{P}_2(n)} 1
			\end{equation*}
			
		\noindent
		for a single standard real Gaussian $X$ which we proved in Lecture One.  A proof of the Wick formula
		may be found in Alexandre Zvonkin's expository article \cite{Zvonkin}.
		
		We now compute the moments of the trace of a GUE matrix $X_N$ using the Wick formula, and then
		take the $N \rightarrow \infty$ limit.  We have
		
			\begin{equation*}
				\begin{split}
					\tau_N[X_N^n] & = \frac{1}{N} \sum_{i:[n] \rightarrow [N]} 
								\E[X_N(i(1)i(2)) X_N(i(2)i(3))) \dots X_N(i(n)i(1))] \\
							&= \frac{1}{N} \sum_{i:[n] \rightarrow [N]} \E \bigg{[} \prod_{k=1}^n X_N(i(k)i\gamma(k)) \bigg{]},
				\end{split}
			\end{equation*}	
			
		\noindent
		where $\gamma = (1\ 2\ \dots \ n)$ is the full forward cycle in the symmetric group $\group{S}(n)$.  Let us 
		apply the Wick formula to each term of this sum, and then use the covariance structure of the matrix
		elements.  We obtain
		
			\begin{equation*}
				\begin{split}
				\E \bigg{[} \prod_{k=1}^n X_N(i(k)i\gamma(k)) \bigg{]} &= \sum_{\pi \in \lattice{P}_2(n)}
							\prod_{\{r,s\} \in \pi} \E[X_N(i(r)i\gamma(r)) X_N(i(s)i\gamma(s))]\\
					&=  N^{-\frac{n}{2}} \sum_{\pi \in \lattice{P}_2(n)}
							\prod_{\{r,s\} \in \pi} \delta_{i(r)i\gamma(s)} \delta_{i(s)i\gamma(r)}.
				\end{split}
			\end{equation*}		
			
		\noindent
		Now, any pair partition of $[n]$ can be viewed as a product of disjoint two-cycles in $\group{S}(n)$.  For
		example, the three pair partitions of $[4]$ enumerated above may be viewed as the fixed
		point free involutions
		
			\begin{equation*}
				(1\ 2)(3\ 4), \ (1\ 3)(2\ 4),\ (1\ 4)(2\ 3)
			\end{equation*}					
			
		\noindent
		in $\group{S}(4)$.  This is a useful shift in perspective because partitions are inert combinatorial objects
		whereas permutations are functions which act on points.  Our computation 
		above may thus be re-written as
		
			\begin{equation*}
				\E \bigg{[} \prod_{k=1}^n X_N(i(k)i\gamma(k)) \bigg{]} =
				  N^{-\frac{n}{2}} \sum_{\pi \in \lattice{P}_2(n)}
							\prod_{k=1}^n \delta_{i(k)i\gamma\pi(k)}.
			\end{equation*}				
			
		\noindent
		Putting this all together and changing order of summation, we obtain
		
			\begin{equation*}
				\begin{split}
					\tau_N[X_N^n] & = N^{1-\frac{n}{2}}  \sum_{i:[n] \rightarrow [N]} 
									\sum_{\pi \in \lattice{P}_2(n)}
									\prod_{k=1}^n \delta_{i(k)i\gamma\pi(k)} \\
								&= N^{1-\frac{n}{2}} \sum_{\pi \in \lattice{P}_2(n)}
									\sum_{i:[n] \rightarrow [N]} \prod_{k=1}^n \delta_{i(k)i\gamma\pi(k)},
				\end{split}
			\end{equation*}			
			
		\noindent
		from which we see that the internal sum is non-zero if and only if the function $i:[n] \rightarrow [N]$
		is constant on the cycles of the permutation $\gamma\pi \in \group{S}(n)$.  In order to build such a function, we must specify 
		one of $N$ possible values to be taken on each cycle.  We thus obtain
		
			\begin{equation*}
				\tau_N[X_N^n] = \sum_{\pi \in \lattice{P}_2(n)} N^{c(\gamma\pi)-1-\frac{n}{2}},
			\end{equation*}
			
		\noindent
		where $c(\sigma)$ denotes the number of cycles in the disjoint cycle decomposition of a permutation
		$\sigma \in \group{S}(n)$.  For example, when $n=3$ we have $\tau_n[X_N^3]=0$ since there are 
		no fixed point free involutions in $\group{S}(3)$.  In order to compute $\tau_N[X_N^4]$, we first
		compute the product of $\gamma$ with all fixed point free involutions in $\group{S}(4)$,
		
			\begin{equation*}
				\begin{split}
					(1\ 2\ 3\ 4)(1\ 2)(3\ 4) &= (1\ 3)(2)(4) \\
					(1\ 2\ 3\ 4)(1\ 3)(2\ 4) &= (1\ 4\ 3\ 2) \\
					(1\ 2\ 3\ 4)(1\ 4)(2\ 3) &= (2\ 4)(1)(3),
				\end{split}
			\end{equation*}
			
		\noindent
		and from this we obtain 
		
			\begin{equation*}
				\tau_N[X_N^4] = 2 + \frac{1}{N^2}.
			\end{equation*}
			
		More generally, $\tau_N[X_N^n]=0$ whenever $n$ is odd since there are no pairings
		on an odd number of points.  When $n=2k$ is even the product
		$\gamma \pi$ has the form
		
			\begin{equation*}
				\gamma\pi=(1\ 2\ \dots\ 2k)(s_1\ t_1)(s_2\ t_2) \dots (s_k\ t_k).
			\end{equation*}
			
\begin{figure}
				\includegraphics{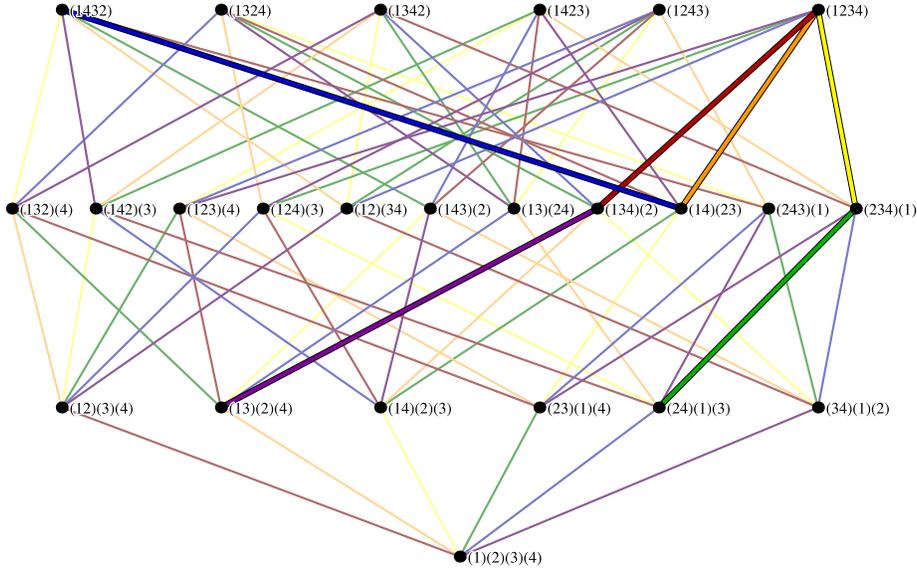}
				\caption{\label{fig:CayleyS4-Path-Highlight}Walks corresponding to the products $\gamma\pi$ in $\group{S}(4)$.}
			\end{figure}

		\noindent
		In this product, each transposition factor $(s_i\ t_i)$ acts either as a ``cut'' or as a ``join'',
		meaning that it may either cut a cycle of $(1\ 2\ \dots\ 2k)(s_1\ t_1) \dots (s_{i-1}\ t_{i-1})$
		in two, or join two disjoint cycles together into one.  More geometrically, we can view the 
		product $\gamma\pi$ as a walk of length $k$ on the (right) Cayley graph of $\group{S}(2k)$; this
		walk is non-backtracking and each step taken augments the distance from the identity permutation by
		$\pm1$, see Figure \ref{fig:CayleyS4-Path-Highlight}. 
		
		 A cut (step towards the identity) occurs when $s_i$ and $t_i$
		reside on the same cycle in the disjoint cycle decomposition of $(1\ 2\ \dots\ 2k)(s_1\ t_1) \dots (s_{i-1}\ t_{i-1})$,
		while a join (step away from the identity) occurs when $s_i$ and $t_i$ are on different cycles.
		In general, the number of cycles in the product will be
		
			\begin{equation*}
				c(\gamma\pi) = 1+\#\text{cuts}-\#\text{joins},
			\end{equation*}
			
		\noindent
		so $c(\gamma\pi)$ is maximal at $c(\gamma\pi)=1+k$ when it is acted on by 
		a sequence of $k$ cut transpositions.  In this case we get a contribution of $N^{1+k-1-k}=N^0$
		to $\tau[X_N^n]$.  In fact, we always have
		
			\begin{equation*}
				\#\text{cuts}-\#\text{joins} = k-2g
			\end{equation*}
			
		\noindent 
		for some non-negative integer $g$, leading to a contribution of the form $N^{-2g}$ and 
		resulting in the formula
		
			\begin{equation*}
				\tau_N[X_N^{2k}] = \sum_{g \geq 0} \frac{\varepsilon_g(2k)}{N^{2g}}
			\end{equation*}
			
		\noindent
		where $\varepsilon_g(2k)$ is the number of products $\gamma\pi$ of the long
		cycle with a fixed point free involution in $\group{S}(2k)$ which terminate at a point
		of the sphere $\partial B(\id,2k-1-2g)$.  We are only interested
		in the first term of this expansion, $\varepsilon_0(2k)$, which counts 
		fixed point free involutions in $\group{S}(2k)$ entirely composed of cuts.
		It is not difficult to see that $(s_1\ t_1) \dots (s_k\ t_k)$ is a sequence of cuts for 
		$\gamma$ if and only if it corresponds to a non-crossing pair partition of $[2k]$, 
		and as we know the number of these is $\Cat_k$. 
		
		We have now shown that
		
			\begin{equation*}
				\lim_{N \rightarrow \infty} \tau_N[X_N^n]=\lim_{N \rightarrow \infty} (\E \otimes \tr_N)[X_N^n] = \begin{cases}
					0, \text{ if $n$ odd}\\
					\Cat_{\frac{n}{2}}, \text{ if $n$ even}
				\end{cases}.
			\end{equation*}
		
		\noindent
		for a GUE matrix $X_N$.  This establishes that
		$X_N$ is an asymptotic random matrix model of a single semicircular random variable
		$X$.  It remains
		to use this fact to construct a sequence of pairs of random matrices which model 
		a pair $X,Y$ of freely independent semicircular random variables.
		
		What should we be looking for?  Let $X^{(1)},X^{(2)}$ be a pair of free semicircular random 
		variables.  Let $e:[n] \rightarrow [2]$ be a function, and apply the free moment-cumulant formula
		to the corresponding mixed moment:
		
			\begin{equation*}
				\begin{split}
				\tau[X^{(e(1))} \dots X^{(e(n))}] &= \sum_{\pi \in \lattice{NC}(n)} \prod_{B \in \pi} \kappa_{|B|}(X^{(e(i))}:i \in B) \\
				&=  \sum_{\pi \in \lattice{NC}_2(n)} \prod_{\{r,s\} \in \pi} \delta_{e(r)e(s)}.
				\end{split}
			\end{equation*}
		
		\noindent
		This reduction occurs because $X^{(1)},X^{(2)}$ are free, so that all mixed free cumulants in these variables vanish.
		Moreover, these variables are semicircular so only order two pure cumulants survive.  We can think of the function
		$e$ as a bicolouring of $[n]$.  The formula for mixed moments of a semicircular pair then becomes
		
			\begin{equation*}
				\tau[X^{(e(1))} \dots X^{(e(n))}] =  \sum_{\pi \in \lattice{NC}_2^{(e)}(n)} 1,
			\end{equation*}
			
		\noindent
		where $\pi \in \lattice{NC}_2^{(e)}(n)$ is the set
		of non-crossing pair partitions of $[n]$ which pair elements of the same colour.  This is very much like the 
		Wick formula for Gaussian expectations, but with Gaussians replaced by semicirculars and summation 
		restricted to non-crossing pairings.
		We need to realize this structure in the combinatorics of GUE random matrices.  
				
		This construction goes as follows.  Let $Z_N^{(e)}(ij),\ 1 \leq e \leq 2,\ 1 \leq i,j \leq N$
		be a collection of $2N^2$ iid centred complex Gaussian random variables of variance $1/N$.  
		Form the corresponding Ginibre matrices $Z_N^{(1)}=[Z_N^{(1)}(ij)],\ Z_N^{(2)}=[Z_N^{(2)}(ij)]$ and 
		GUE matrices $X_N^{(1)} = \frac{1}{2}(Z_N^{(1)}+(Z_N^{(1)})^*),\ X_N^{(2)} = \frac{1}{2}(Z_N^{(2)}+(Z_N^{(2)})^*)$.
		The resulting covariance structure of matrix elements is 
		
			\begin{equation*}
				\E[X_N^{(p)}(ij)\overline{X_N^{(q)}(kl)}] = \E[X_N^{(p)}(ij)X_N^{(q)}(lk)] = \frac{\delta_{ik}\delta_{jl}\delta_{pq}}{N}. 
			\end{equation*}
			
		\noindent
		We can prove that $X_N^{(1)},X_N^{(2)}$ are asymptotically free by showing that
				
			\begin{equation*}
				\lim_{N \rightarrow \infty} \tau_N[X_N^{(e(1))} \dots X_N^{(e(n))}]  = |\lattice{NC}_2^{(e)}(n)|,
			\end{equation*} 
			
		\noindent
		and this can in turn be proved using the Wick formula and the above covariance structure.
		Computations almost exactly like those appearing in the one-matrix case lead to the formula
		
			\begin{equation*}
				\tau_N[X_N^{(e(1))} \dots X_N^{(e(n))}] = \sum_{\pi \in \lattice{P}_2^{(e)}(n)} N^{c(\gamma\pi)-1-\frac{n}{2}},
			\end{equation*}
			
		\noindent
		with the summation being taken over the set $\lattice{P}_2^{(e)}(n)$ of pairings on 
		$[n]$ which respect the colouring $e:[n] \rightarrow [2]$.  Arguing as above, each such 
		pairing makes a contribution of the form $N^{-2g}$ for some $g \geq 0$, and 
		those which make contributions on the leading order $N^0$ correspond to sequences
		of cut transpositions for the full forward cycle $\pi$, which we know come from non-crossing pairings.
		So in the limit $N \rightarrow \infty$ this expectation converges to $|\lattice{NC}_2^{(e)}(n)|$,
		as required.
			
		\subsection{Random matrix model of a free pair with one semicircle}
		In the previous subsection we modelled a free pair of semicircular random variables 
		$X,Y$ living in an abstract non-commutative probability space $(\algebra{A},\tau)$ 
		using a sequence of independent GUE random matrices $X_N,Y_N$ living in random matrix
		space $(\algebra{A}_N,\tau_N)$.  
		
		It is reasonable to wonder whether we have not overlooked 
		the possibility of modelling $X,Y$ in a simpler way, 
		namely using deterministic matrices.  Indeed, we have
		
			\begin{equation*}
				\tau[X^n] = \int\limits_{\field{R}} t^n \mu_X(\mathrm{d}t)
			\end{equation*}
		
		\noindent
		with
		
			\begin{equation*}
				\mu_X(\mathrm{d}t) = \frac{1}{2\pi} \sqrt{4-t^2} \mathrm{d}t
			\end{equation*}
		
		\noindent	
		the Wigner semicircle measure, and this fact leads to a
		deterministic matrix model for $X$.
		For each $N \geq 1$, define the $N^{\text{th}}$ classical 
		locations $L_N(1) < L_N(2) < \dots <L_N(N)$ of  $\mu_X$ implicitly by
		
			\begin{equation*}
				\int\limits_{-2}^{L_N(i)} \mu_X(\mathrm{d}t) = \frac{i}{N}.
			\end{equation*}
			
		\noindent
		That is, we start at $t=-2$ and integrate along the semicircle until a mass of $i/N$ is
		achieved, at which time we mark off the corresponding location $L_N(i)$ on the 
		$t$-axis.  The measure $\mu_N$ which places mass $1/N$
		at each of the $N^{\text{th}}$ classical locations converges weakly to $\mu_X$
		as $N \rightarrow \infty$.  Consequently, the diagonal matrix $X_N$ with entries $X_N(ij) = \delta_{ij}L_N(i)$
		is a random variable in deterministic matrix space 
		$(\algebra{M}_N(\field{C}),\tr_N)$ which models $X$,
		
			\begin{equation*}
				\lim_{N \rightarrow \infty} \tr_N[X_N^n] = \tau[X^n].
			\end{equation*}
		
		\noindent
		Since $X$ and $Y$ are equidistributed,
		putting $Y_N:=X_N$ we have that $X_N$ models $X$ and 
		$Y_N$ models $Y$.  However, $X_N$ and $Y_N$ are not asymptotically free.
		Indeed, asymptotic freeness of $X_N$ and $Y_N$ would imply that
		
			\begin{equation*}
				\lim_{N \rightarrow \infty} \tr_N[X_NY_N] = \lim_{N \rightarrow \infty} \tr_N[X_N]
				\lim_{N \rightarrow \infty}\tr_N[X_N]=0,
			\end{equation*}
			
		\noindent
		but instead we have
		
			\begin{equation*}
				\tr_N[X_NY_N] = \frac{L_N(1)^2 + \dots +L_N(N)^2}{N},
			\end{equation*}
			
		\noindent
		the mean squared classical locations of the Wigner measure,
		which is strictly positive and increasing in $N$.  Thus while $X_N$ and $Y_N$ model 
		$X$ and $Y$ respectively, they cannot model the free relation between them.
		However, this does not preclude the possibility that a pair of free random variables can
		be modelled by one random and one deterministic matrix.  
		
		Let $X$ and $Y$
		be a pair of free random variables with $X$ semicircular, and $Y$ of arbitrary 
		distribution.  
		Let $X_N$ be a sequence of GUE matrices modelling $X$, and 
		suppose that $Y_N$ is a sequence of deterministic matrices
		modelling $Y$,
		
			\begin{equation*}
				\lim_{N \rightarrow \infty} \tr_N[Y_N^n] = \tau[Y^n].
			\end{equation*}
			
		\noindent
		$X_N$ lives in random matrix space 
		$(\algebra{A}_N,\tau_N)=(L^{\infty-}(\Omega,\algebra{F},P) \otimes \algebra{M}_N(\field{C}), \E \otimes \tr_N)$
		while $Y_N$ lives in deterministic matrix space $(\algebra{M}_N(\field{C}),\tr_N)$, so a priori
		it is meaningless to speak of the potential asymptotic free independence of $X_N$ and $Y_N$.
		However, we
		may think of 
		a deterministic matrix as a random matrix whose entries are constant random variables
		in $L^{\infty-}(\Omega,\algebra{F},P)$.  This corresponds to an embedding of deterministic 
		matrix space in random matrix space satisfying $\tau_N|_{\algebra{M}_N(\field{C})} 
		= (\E \otimes \tr_N)|_{\algebra{M}_N(\field{C})}= \tr_N$.
		From this point of view, $Y_N$ is a random matrix 
		model of $Y$ and we can consider the possibility that $X_N,Y_N \in \algebra{A}_N$ are
		asymptotically free with respect to $\tau_N$.
		We now show that this is indeed the case.
		
		As in the previous subsection, we proceed by identifying the combinatorial structure
		governing the target pair $X,Y$ and then looking for this same structure in the 
		$N \rightarrow \infty$ asymptotics of $X_N,Y_N$.  Our target is a pair of free random 
		variables with $X$ semicircular and $Y$ arbitrary.  Understanding their joint distribution 
		means understanding the collection of mixed moments
		
			\begin{equation*}
				\tau[X^{p(1)}Y^{q(1)} \dots X^{p(n)}Y^{q(n)}],
			\end{equation*}
			
		\noindent
		with $n \geq 1$ and $p,q:[n] \rightarrow \{0,1,2,\dots\}$.  This amounts to understanding mixed moments
		of the form
		
			\begin{equation*}
				\tau[XY^{q(1)} \dots XY^{q(n)}],
			\end{equation*}
			
		\noindent
		since we can artificially insert copies of $Y^0=1_\algebra{A}$ to break up 
		powers of $X$ greater than one.  We can expand this expectation using the free 
		moment-cumulant formula and simplify the resulting expression using the fact that 
		mixed cumulants in free random variables vanish.  Further simplification results 
		from the fact that, since $X$ is semicircular, its only non-vanishing pure cumulant
		is $\kappa_2(X)=1$.  This leads to a formula for $\tau[XY^{q(1)} \dots XY^{q(n)}]$
		which is straightforward but whose statement requires some notions which 
		we have not covered (in particular, the complement of a non-crossing partition,
		see \cite{NS}).  However, in the case where $\tau$ is a tracial expecation,
		meaning that $\tau[AB]=\tau[BA]$, the formula in question can be stated 
		more simply as 
		
			\begin{equation*}
				\tau[XY^{q(1)} \dots XY^{q(n)}] = \sum_{\pi \in \lattice{NC}_2(n)}
				\tau_{\pi\gamma}[Y^{q(1)},\dots,Y^{q(n)}].
			\end{equation*}
			
		\noindent 
		Here, as in the last subsection, we think of a pair partition $\pi \in \lattice{P}_2(n)$ as a product of 
		disjoint two-cycles in the symmetric group $\group{S}(n)$, and $\gamma$ is the full forward
		cycle $(1\ 2\ \dots\ n)$.  Given a permutation $\sigma \in \group{S(n)}$, the expression 
		$\tau_\sigma[A_1,\dots,A_N]$ is defined to be the product of $\tau$ extended over the 
		cycles of $\sigma$.  For example,
		
			\begin{equation*}
				\tau_{(1\ 6\ 2)(4\ 5)(3)}[A_1,A_2,A_3,A_4,A_5,A_6] = \tau[A_1A_6A_2]\tau[A_4A_5]\tau[A_3].
			\end{equation*}
			
		\noindent
		This definition is kosher since $\tau$ is  tracial.
		We now have our proof strategy: we will prove that $X_N,Y_N$ are 
		asymptotically free by showing that
		
			\begin{equation*}
				\lim_{N \rightarrow \infty} \tau_N[X_NY_N^{q(1)} \dots X_NY_N^{q(n)}] =
				\sum_{\pi \in \lattice{NC}_2(n)}
				\tau_{\pi\gamma}[Y^{q(1)},\dots,Y^{q(n)}].
			\end{equation*}
			
		The computation proceeds much as in the last section --- we expand everything in sight
		and apply the Wick formula.  We have
		
			\begin{equation*}
				\begin{split}
				&\tau_N[X_NY_N^{q(1)} \dots X_NY_N^{q(n)}] \\ 
				&= \frac{1}{N} \sum_a 
				\E[X_N(a(1)a(2))Y_N^{q(1)}(a(2)a(3)) \dots X_N(a(2n-1)a(2n))Y_N^{q(n)}(a(2n)a(1))],
				\end{split}
			\end{equation*}
			
		\noindent
		the summation being over all functions $a:[2n] \rightarrow [N]$.  Let us reparameterize each 
		term of the sum with $i,j:[n] \rightarrow [N]$ defined by
		
			\begin{equation*}
				(a(1),a(2),\dots,a(2n-1),a(2n)) = (i(1),j(1),\dots,i(n),j(n)).
			\end{equation*}
			
		\noindent
		Our computation so far becomes
		
			\begin{equation*}
				\tau_N[X_NY_N^{q(1)} \dots X_NY_N^{q(n)}] = \frac{1}{N} \sum_{i,j} 
				\E\bigg{[} \prod_{k=1}^n X_N(i(k)j(k)) \bigg{]} \prod_{k=1}^n Y_N^{q(k)}(j(k)i\gamma(k)).
			\end{equation*}
			
		\noindent
		Applying the Wick formula, the calculation evolves as follows:
		
			\begin{equation*}
				\begin{split}
					&\tau_N[X_NY_N^{q(1)} \dots X_NY_N^{q(n)}] \\
					&= \frac{1}{N} \sum_{i,j} \sum_{\pi \in \lattice{P}_2(n)} \prod_{\{r,s\} \in \pi} 
						\E[X_N(i(r)j(r))X_N(i(s)j(s))] \prod_{k=1}^n Y_N^{q(k)}(j(k)i\gamma(k)) \\
					&= N^{-1-\frac{n}{2}} \sum_{i,j} \sum_{\pi \in \lattice{P}_2(n)} 
						\prod_{k=1}^n \delta_{i(k)j\pi(k)} Y_N^{q(k)}(j(k)i\gamma(k)) \\
					&= N^{-1-\frac{n}{2}} \sum_{\pi \in \lattice{P}_2(n)} \sum_j
						\prod_{k=1}^n Y_N^{q(k)}(j(k)j\pi\gamma(k)) \\
					&= N^{-1-\frac{n}{2}} \sum_{\pi \in \lattice{P}_2(n)} \Tr_{\pi\gamma}[Y_N^{q(1)},\dots,Y_N^{q(n)}] \\
					&= \sum_{\pi \in \lattice{P}_2(n)} N^{c(\pi\gamma)-1-\frac{n}{2}} \tr_{\pi\gamma}[Y_N^{q(1)},\dots,Y_N^{q(n)}].
				\end{split}
			\end{equation*}
			
		\noindent
		As in the previous subsection, the dominant contributions to this sum are of order $N^0$ and come from 
		those pair partitions $\pi \in \lattice{P}_2(n)$ for which $c(\pi\gamma)$ is maximal, and these are 
		the non-crossing pairings.  Hence we obtain
		
			\begin{equation*}
				\lim_{N \rightarrow \infty} \tau_N[X_NY_N^{q(1)} \dots X_NY_N^{q(n)}] =
				\sum_{\pi \in \lattice{NC}_2(n)}
				\tau_{\pi\gamma}[Y^{q(1)},\dots,Y^{q(n)}],
			\end{equation*}
			
		\noindent
		as required.
								
		\subsection{Random matrix model of an arbitrary free pair}
		In the last section we saw that a pair of free random variables can be modelled
		by one random and one deterministic matrix provided that at least one of the target variables 
		is semicircular.  In this case, the semicircular target is modelled by
		a sequence of GUE random matrices.
				
		In this section we show that any pair of free random variables can be modelled by 
		one random and one deterministic matrix, provided each target variable can 
		be individually modelled by a sequence of deterministic matrices.  
		The idea is to randomly rotate one of the deterministic matrix models so as
		to create the free relation.
		
		Let $X,Y$ be a pair of free random variables living in an abstract
		non-commutative probability space $(\algebra{A},\tau)$.  We make no assumption on 
		their moments.  What we assume is the existence of a pair of deterministic matrix
		models
		
			\begin{equation*}
				\tau[X^n]=\lim_{N \rightarrow \infty} \tr_N[X_N^n] , \quad
				 \tau[Y^n]=\lim_{N \rightarrow \infty} \tr_N[Y_N^n].
			\end{equation*}
			
		\noindent
		If $X,Y$ happen to have distributions $\mu_X,\mu_Y$ 
		which are compactly supported probability measures on $\field{R}$,
		then such models can always be constructed.  In particular, this
		will be the case if $X,Y$ are bounded self-adjoint random variables living
		in a $*$-probability space.
		
		As in the previous subsection,
		we view $X_N,Y_N$ as random matrices with constant entries so 
		that they reside in random matrix space $(\algebra{A}_N,\tau_N)$, with 
		the $\E$ part of $\tau_N=\E \otimes \tr_N$ acting trivially.  As we saw above,
		there is no guarantee that $X_N,Y_N$ are asymptotically free.  On the other hand, we also 
		saw that special pairs of free random variables can be modelled by 
		one random and one deterministic matrix.  Therefore it is reasonable to hope that making
		$X_N$ genuinely random might lead to asymptotic freeness.  We have to randomize
		$X_N$ in such a way that its moments will be preserved.
		This can be achieved via conjugation by a unitary random
		matrix $U_N \in \algebra{A}_N$,
				
			\begin{equation*}
				X_N \mapsto U_NX_NU_N^*.
			\end{equation*} 
			
		\noindent
		The deterministic matrix $X_N$ and its randomized version $U_NX_NU_N^*$
		have the same moments since
		
			\begin{equation*}
				\begin{split}
				\tau_N[(U_NX_NU_N^*)^n] &= (\E \otimes \tr_N)[(U_NX_NU_N^*)^n] \\ 
				&= (\E \otimes\tr_N)[U_NX_N^nU_N^*] \\
				&= (\E \otimes\tr_N)[U_N^*U_NX_N^n] \\
				&= (\E \otimes \tr_N)[X_N^n] \\
				&= \tau_N[X_N^n].
				\end{split}
			\end{equation*}
			
		\noindent
		Consequently, the sequence $U_NX_NU_N^*$ is a random
		matrix model for $X$.  
				
		We aim to prove that $U_NX_NU_N^*$ and $Y_N$ are asymptotically free.
		Since we are making no assumptions on the limiting 
		variables $X,Y$, we cannot verify this by looking for special structure in 
		the limiting mixed moments of $U_NX_NU_N^*$ and $Y_N$, as we did above.
		Instead, we must verify asymptotic freeness directly, using the definition:
		
			\begin{equation*}
				\lim_{N \rightarrow \infty} \tau_N[f_1(U_NX_NU_N^*)g_1(Y_N) \dots 
				f_n(U_NX_NU_N^*)g_n(Y_N)] = 0
			\end{equation*}
			
		\noindent
		whenever $f_1,g_1,\dots,f_n,g_n$ are polynomials such that
		
			\begin{equation*}
				\lim_{N \rightarrow \infty} \tau_N[f_1(U_NX_NU_N^*)]=\lim_{N \rightarrow \infty} \tau_n[g_1(Y_N)]=
				\dots 
				= \lim_{N \rightarrow \infty} \tau_N[f_n(U_NX_NU_N^*)]
				=\lim_{N \rightarrow \infty} \tau_n[g_n(Y_N)]=0.
			\end{equation*}
			
		\noindent
		Though the brute force verification of this criterion may seem an impossible task,
		we will see that it can be accomplished for a well-chosen sequence of 
		unitary random matrices $U_N$.  Let us advance as far as possible before 
		specifying $U_N$ precisely.
		
		As an initial reduction, note the identity
		
			\begin{equation*}
				\begin{split}
				&\tau_N[f_1(U_NX_NU_N^*)g_1(Y_N) \dots f_n(U_NX_NU_N^*)g_n(Y_N)] \\
				&=\tau_N[U_Nf_1(X_N)U_N^*g_1(Y_N) \dots U_Nf_n(X_N)U_N^*g_n(Y_N)].
				\end{split}
			\end{equation*}
			
		\noindent
		Since the $f_i$'s and $g_j$'s are polynomials and $\tau_N$ is linear, the right hand side
		of this equation may be expanded as a sum of monomial expectations,
		
			\begin{equation*}
				\begin{split}
				&\tau_N[U_Nf_1(X_N)U_N^*g_1(Y_N) \dots U_Nf_n(X_N)U_N^*g_n(Y_N)] \\
				= &\sum_{p,q} c(pq) \tau_N[U_NX_N^{p(1)}U_N^*Y_N^{q(1)} \dots U_NX_N^{p(n)}U_N^*Y_N^{q(n)}]
				\end{split}
			\end{equation*}
			
		\noindent
		weighted by some scalar coefficients $c(pq)$, the sum being over 
		functions $p:[n] \rightarrow \{0,\dots,\max \deg f_i\},q:[n] \rightarrow \{0,\dots,\max \deg g_j\}$.
		Each monomial expectation can in turn be expanded as 
		
			\begin{equation*}
				\begin{split}
				&\tau_N[U_NX_N^{p(1)}U_N^*Y_N^{q(1)} \dots U_NX_N^{p(n)}U_N^*Y_N^{q(n)}] \\
				=&\frac{1}{N} \sum_{a} 
				\E[U_N(a(1)a(2))X_N^{p(1)}(a(2)a(3)) \dots U_N^*(a(4n-1)a(4n))Y_N^{q(n)}(a(4n)a(1))] \\
				=&\frac{1}{N} \sum_{a} 
				\E[U_N(a(1)a(2))X_N^{p(1)}(a(2)a(3)) \dots \overline{U}_N(a(4n)a(4n-1))Y_N^{q(n)}(a(4n)a(1))].
				\end{split}
			\end{equation*}
			
		\noindent
		Let us reparameterize the summation index $a:[4n] \rightarrow [N]$ by a quadruple of functions
		$i,j,i',j':[n] \rightarrow [N]$ according to 
		
			\begin{equation*}
				\begin{split}
				&(a(1),a(2),a(3),a(4),\dots,a(4n-3),a(4n-2),a(4n-1),a(4n)) \\
				=&(i(1),j(1),j'(1),i'(1),\dots,i(n),j(n),j'(n),i'(n)).
				\end{split}
			\end{equation*}
			
		\noindent
		Our monomial expectations then take the more streamlined form
		
			\begin{equation*}
				\begin{split}
				&\tau_N[U_NX_N^{p(1)}U_N^*Y_N^{q(1)} \dots U_NX_N^{p(n)}U_N^*Y_N^{q(n)}] \\
				=&\frac{1}{N} \sum_{i,j,i',j'} \E\bigg{[} \prod_{k=1}^n U_N(i(k)j(k)) \overline{U}_N(i'(k)j'(k))\bigg{]}
				\prod_{k=1}^n X_N^{p(k)}(j(k)j'(k))Y_N^{q(k)}(i'(k)i\gamma(k)),
				\end{split}
			\end{equation*}
			
		\noindent
		where as always $\gamma=(1\ 2\ \dots \ n)$ is the full forward cycle in the symmetric group 
		$\group{S}(n)$.  In order to 
		go any further with this calculation, we must deal with the correlation functions
		
			\begin{equation*}
				\E\bigg{[} \prod_{k=1}^n U_N(i(k)j(k)) \overline{U}_N(i'(k)j'(k))\bigg{]}.
			\end{equation*}
			
		\noindent
		of the matrix elements of $U_N$.  We would like to have an analogue of the 
		Wick formula which will enable us to address these correlation functions.  
		A formula of this type is 
		known for random matrices sampled from the Haar probability measure on 
		the unitary group $\group{U}(N)$.
		
		Haar-distributed unitary matrices are the 
		second most important class of random matrices after GUE matrices.  Like GUE
		matrices, they can be constructively obtained from Ginibre matrices.  Let 
		$\tilde{Z}_N=\sqrt{N}Z_N$ be an $N \times N$ random matrix whose entries $\tilde{Z}_N(ij)$ 
		are iid complex Gaussian random variables of mean zero and variance one.
		This is a renormalized version of the Ginibre matrix which we previously used to
		construct a GUE random matrix.  The Ginibre matrix 
		$\tilde{Z}_N$ is almost surely non-singular.  Applying the Gram-Schmidt 
		orthonormalization procedure to the columns of $\tilde{Z}_N$, we obtain 
		a random unitary matrix $U_N$ whose distribution in the unitary group 
		$\group{U}(N)$ is given by the Haar probability measure.  The entries $U_N(ij)$ are bounded 
		random variables, so $U_N$ is a non-commutative random 
		variable living in random matrix space $(\algebra{A}_N,\tau_N)$.  
		The eigenvalues $\lambda_N(1)=e^{\mathbf{i}\theta_N(1)},\dots,\lambda_N(N)=e^{\mathbf{i}\theta_N(N)}$,
		$0 \leq \theta_N(1) \leq \dots \leq \theta_N(N) \leq 2\pi$ of $U_N$ form a
		random point process on the unit circle with joint distribution
		
			\begin{equation*}
				P(\theta_N(1) \in I_1,\dots,\theta_N(N)\in I_N) \propto \int\limits_{I_1} \dots \int\limits_{I_N}
				e^{-N^2 \mathcal{H}(\theta_1,\dots,\theta_N)} \mathrm{d}\theta_1 \dots \mathrm{d}\theta_N
			\end{equation*}	
			
		\noindent
		for any intervals $I_1,\dots,I_N \subseteq [0,2\pi]$, where $\mathcal{H}$ is the log-gas Hamiltonian \cite{Forrester}
		
			\begin{equation*}
				\mathcal{H}(\theta_1,\dots,\theta_N) = - \frac{1}{N^2}
				\sum_{1 \leq i \neq j \leq N} \log |e^{\mathbf{i}\theta_i}-e^{\mathbf{i}\theta_j}|.
			\end{equation*}
			
		\noindent
		The random point process on the unit circle driven by this Hamiltonian is known as the Circular Unitary Ensemble,
		 and $U_N$ is termed a CUE random matrix.  As with GUE random matrices, almost any 
		 question about the spectrum of CUE random matrices can be answered using this explicit formula,
		 see e.g. \cite{Diaconis} for a survey of many interesting results.
		 
		 We are not interested in the eigenvalues of CUE matrices, but rather in the correlation functions
		 of their matrix elements.  These can be handled using a Wick-type formula 
		 known as the Weingarten formula, after the American physicist Donald H. Weingarten\footnote{Further information
		 regarding Weingarten and his colleagues in the first Fermilab theory group may be found at
		 \textsf{http://bama.ua.edu/$\sim$lclavell/Weston/}}.  Like
		 the Wick formula, the Weingarten formula is a combinatorial rule which reduces the computation 
		 of general correlation functions to the computation of a special class of correlations.  
		 Unfortunately,
		 the Weingarten formula is more complicated than the Wick formula.
		 It reads:
		 
		 	\begin{equation*}
				\E\bigg{[} \prod_{k=1}^n U_N(i(k)j(k)) \overline{U}_N(i'(k)j'(k))\bigg{]} = \sum_{\rho,\sigma \in \group{S}(n)}
				\delta_{i\sigma,i'}\delta_{j\rho,j'} \E\bigg{[} \prod_{k=1}^n U_N(kk) \overline{U}_N(k\rho^{-1}\sigma(k))\bigg{]}.
			\end{equation*}
		
		\noindent
		Note that his formula only makes sense 
		when $N \geq n$, and
		instead of a sum over fixed point free involutions we are faced with
		a double sum over all of $\group{S}(n)$.  Worse still, the Weingarten formula 
		does not reduce our problem to the computation of pair correlators, but only to the 
		computation of arbitrary permutation correlators
		
			\begin{equation*}
				\E\bigg{[} \prod_{k=1}^n U_N(kk) \overline{U}_N(k\pi(k))\bigg{]}, \quad \pi \in \group{S}(n),
			\end{equation*}	 
			
		\noindent
		and these have a rather complicated structure.  Their computation is the subject of a large
		literature both in physics and mathematics, a unified treatment of which may be
		found in \cite{CMN}.  We delay dealing with these averages for the moment and press
		on in our calculation.
				
		We return to the expression
		
			\begin{equation*}
				\begin{split}
				&\tau_N[U_NX_N^{p(1)}U_N^*Y_N^{q(1)} \dots U_NX_N^{p(n)}U_N^*Y_N^{q(n)}] \\
				=&\frac{1}{N} \sum_{i,j,i',j'} \E\bigg{[} \prod_{k=1}^n U_N(i(k)j(k)) \overline{U}_N(i'(k)j'(k))\bigg{]}
				\prod_{k=1}^n X_N^{p(k)}(j(k)j'(k))Y_N^{q(k)}(i'(k)i\gamma(k)),
				\end{split}
			\end{equation*}
			
		\noindent
		and apply the Weingarten formula.  The calculation evolves as follows:
		
			\begin{equation*}
				\begin{split}
				&\tau_N[U_NX_N^{p(1)}U_N^*Y_N^{q(1)} \dots U_NX_N^{p(n)}U_N^*Y_N^{q(n)}] \\
				=&\frac{1}{N} \sum_{i,j,i',j'} \sum_{\rho,\sigma \in \group{S}(n)} \delta_{i\sigma,i'}\delta_{j\rho,j'}
				 \E\bigg{[} \prod_{k=1}^n U_N(kk) \overline{U}_N(k\rho^{-1}\sigma(k))\bigg{]}
				\prod_{k=1}^n X_N^{p(k)}(j(k)j'(k))Y_N^{q(k)}(i'(k)i\gamma(k)) \\
				=&\frac{1}{N}\sum_{\rho,\sigma \in \group{S}(n)}  \E\bigg{[} \prod_{k=1}^n U_N(kk) \overline{U}_N(k\rho^{-1}\sigma(k))\bigg{]}
				\sum_{i',j} \prod_{k=1}^n X_N^{p(k)}(j(k)j\rho(k))Y_N^{q(k)}(i'(k)i\sigma^{-1}\gamma(k)) \\
				=& \frac{1}{N}\sum_{\rho,\sigma \in \group{S}(n)}  \E\bigg{[} \prod_{k=1}^n U_N(kk) \overline{U}_N(k\rho^{-1}\sigma(k))\bigg{]}
				\Tr_\rho(X_N^{p(1)},\dots,X_N^{p(n)}) \Tr_{\sigma^{-1}\gamma}(Y_N^{p(1)},\dots,Y_N^{p(n)}) \\
				=&\sum_{\rho,\sigma \in \group{S}(n)}  \E\bigg{[} \prod_{k=1}^n U_N(kk) \overline{U}_N(k\rho^{-1}\sigma(k))\bigg{]}
				N^{c(\rho)+c(\sigma^{-1}\gamma)-1}\tr_\rho(X_N^{p(1)},\dots,X_N^{p(n)}) \tr_{\sigma^{-1}\gamma}(Y_N^{p(1)},\dots,Y_N^{p(n)}) .
				\end{split}
			\end{equation*}
			
		\noindent
		At this point we are forced to deal with the permutation correlators $\E[\prod U_N(kk) \overline{U}_N(k\pi(k))]$.
		Perhaps the most appealing presentation of these expectations is as a power series in 
		$N^{-1}$.  It may be shown \cite{Novak:BCP} that
		
			\begin{equation*}
				\E\bigg{[} \prod_{k=1}^n U_N(kk) \overline{U}_N(k\pi(k))\bigg{]} = 
				\frac{1}{N^n}\sum_{r=0}^{\infty} (-1)^r \frac{c_{n,r}(\pi)}{N^r},
			\end{equation*}

		\noindent
		for any $\pi \in \group{S}(n)$, where the coefficient $c_{n,r}(\pi)$ equals the number of 
		factorizations
		
			\begin{equation*}
				\pi = (s_1\ t_1) \dots (s_r\ t_r)
			\end{equation*}
			
		\noindent
		of $\pi$ into $r$ transpositions $(s_i\ t_i) \in \group{S}(n),\ s_i<t_i$, which
		have the property that 
		
			\begin{equation*}
				t_1 \leq \dots \leq t_r.
			\end{equation*}
			
		\noindent
		This series is absolutely convergent for $N \geq n$, but divergent for $N<n$.
		This will not trouble us since we are looking for $N \rightarrow \infty$
		asymptotics with $n$ fixed.  Indeed, let $|\pi|=n-c(\pi)$ denote the distance 
		from the identity permutation to $\pi$ in the Cayley graph of $\group{S}(n)$.
		Then, since any permutation is either even or odd, we have
		
			\begin{equation*}
				\begin{split}
				\E\bigg{[} \prod_{k=1}^n U_N(kk) \overline{U}_N(k\pi(k))\bigg{]} &= 
				\frac{1}{N^n}\sum_{r=0}^{\infty} (-1)^r \frac{c_{n,r}(\pi)}{N^r} \\
				&=\frac{(-1)^{|\pi|}}{N^{n+|\pi|}}\sum_{g=0}^{\infty}\frac{c_{n,|\pi|+2g}(\pi)}{N^{2g}} \\
				&= \frac{a(\pi)}{N^{n+|\pi|}} + O\bigg{(}\frac{1}{N^{n+|\pi|+2}} \bigg{)},
				\end{split}
			\end{equation*}
		
		\noindent
		where $a(\pi)=(-1)^{|\pi|} c_{n,|\pi|}(\pi)$ is the leading asymptotics.  
		We may now continue our calculation:
		
			\begin{equation*}
				\begin{split}
				&\tau_N[U_NX_N^{p(1)}U_N^*Y_N^{q(1)} \dots U_NX_N^{p(n)}U_N^*Y_N^{q(n)}] \\
				=&\sum_{\rho,\sigma \in \group{S}(n)}  \bigg{(} \frac{a(\rho^{-1}\sigma))}{N^{n+|\rho^{-1}\sigma|}} + 
					O\bigg{(}\frac{1}{N^{n+|\rho^{-1}\sigma|+2}} \bigg{)} \bigg{)}
				N^{c(\rho)+c(\sigma^{-1}\gamma)-1}\tr_\rho(X_N^{p(1)},\dots,X_N^{p(n)}) \tr_{\sigma^{-1}\gamma}(Y_N^{p(1)},\dots,Y_N^{p(n)}) \\
				=&\sum_{\rho,\sigma \in \group{S}(n)}  \bigg{(} a(\rho^{-1}\sigma)+ 
					O\bigg{(}\frac{1}{N^2} \bigg{)} \bigg{)}
				N^{|\gamma|-|\rho|-|\rho^{-1}\sigma|-|\sigma^{-1}\gamma|}
				\tr_\rho(X_N^{p(1)},\dots,X_N^{p(n)}) \tr_{\sigma^{-1}\gamma}(Y_N^{p(1)},\dots,Y_N^{p(n)}).
				\end{split}
			\end{equation*}
			
\begin{figure}
				\includegraphics{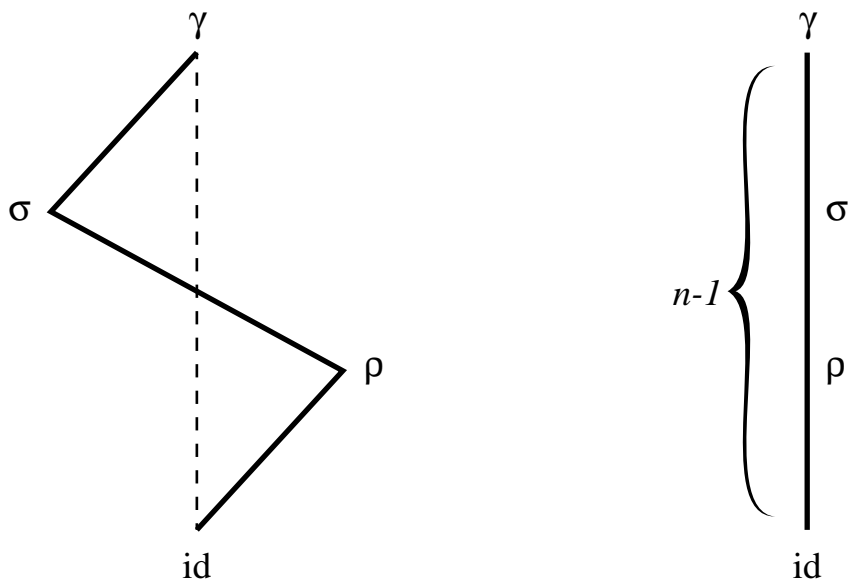}
				\caption{\label{fig:three-leg-b}Only geodesic paths survive in the large $N$ limit}
			\end{figure}

		Putting everything together, we have shown that
		
			\begin{equation*}	
				\begin{split}
				&\tau_N[U_Nf_1(X_N)U_N^*g_1(Y_N) \dots U_Nf_n(X_N)U_N^*g_n(Y_N)] \\
				=&\sum_{\rho,\sigma \in \group{S}(n)}  \bigg{(} a(\rho^{-1}\sigma)+ 
					O\bigg{(}\frac{1}{N^2} \bigg{)} \bigg{)}
				N^{|\gamma|-|\rho|-|\rho^{-1}\sigma|-|\sigma^{-1}\gamma|}
				\tr_\rho(f_1(X_N),\dots,f_n(X_N)) \tr_{\sigma^{-1}\gamma}(g_1(Y_N),\dots,g_n(Y_N)),
				\end{split}
			\end{equation*}
			
		\noindent
		and it remains to show that the $N \rightarrow \infty$ limit of this complicated expression is zero.
		To this end, consider the order $|\gamma|-|\rho|-|\rho^{-1}\sigma|-|\sigma^{-1}\gamma|$ of the
		$\rho,\sigma$ term in this sum.  The positive part, $|\gamma|=n-1$, is simply the length of any 
		geodesic joining the identity permutation to $\gamma$ in the Cayley graph of $\group{S}(n)$.
		The negative part,
		$-|\rho|-|\rho^{-1}\sigma|-|\sigma^{-1}\gamma|$, is the length of a walk from the identity to
		$\gamma$ made up of three legs: a geodesic from $\id$ to $\rho$, followed by a geodesic from
		$\rho$ to $\sigma$, followed by a geodesic from $\sigma$ to $\gamma$.  Thus the order of
		the $\rho,\sigma$ term is at most $N^0$, and this occurs precisely when $\rho$ and $\sigma$
		lie on a geodesic from $\id$ to $\gamma$, see Figure \ref{fig:three-leg-b}.  Thus
		
			\begin{equation*}	
				\begin{split}
				&\lim_{N \rightarrow \infty} \tau_N[U_Nf_1(X_N)U_N^*g_1(Y_N) \dots U_Nf_n(X_N)U_N^*g_n(Y_N)] \\
				=&\sum_{|\rho|+|\rho^{-1}\sigma|+|\sigma^{-1}\gamma|=|\gamma|} a(\rho^{-1}\sigma)
				\tau_\rho(f_1(X),\dots,f_n(X)) \tau_{\sigma^{-1}\gamma}(g_1(Y),\dots,g_n(Y)).
				\end{split}
			\end{equation*}
			
		\noindent
		Since 
		
			\begin{equation*}
				\tau[f_1(X)]=\tau[g_1(Y)] = \dots = \tau[f_n(X)]=\tau[g_n(Y)]=0,
			\end{equation*}
			
		\noindent
		in order to show that the sum on the right has all terms equal to zero it suffices to show
		that the condition $|\rho|+|\rho^{-1}\sigma|+|\sigma^{-1}\gamma|=|\gamma|$ forces either
		$\rho$ or $\sigma^{-1}\gamma$ to have a fixed point.  This is because $\tau_{\rho}$ and
		$\tau_{\sigma^{-1}\gamma}$ are products determined by the cycle 
		structure of the indexing permutation.  Since $\rho,\sigma$ lie on a
		geodesic $\id \rightarrow \gamma$, we have $|\rho|+|\sigma^{-1}\gamma| \leq |\gamma|=n-1$,
		so that one of $\rho$ or $\sigma^{-1}\gamma$ is a product of at most
		$(n-1)/2$ transpositions.  In the extremal case, all of these transpositions are joins, leading 
		to a permutation consisting of an $(n-1)$-cycle and a fixed point.

		\subsection{$\text{GUE}+\text{GUE}$}
		Imagine that we had been enumeratively lazy in our construction of the GUE 
		matrix model of a free semicircular pair, and had only shown that two iid
		GUE matrices $X_N^{(1)},X_N^{(2)}$ are asymptotically free without 
		determining their individual limiting distributions.  We could then appeal to 
		the free central limit theorem to obtain that the limit distribution of the random 
		matrix 
		
			\begin{equation*}
				S_N = \frac{X_N^{(1)}+ \dots + X_N^{(n)}}{\sqrt{N}},
			\end{equation*}
			
		\noindent
		where the $X_N^{(i)}$'s are iid GUE samples, is standard semicircular.  On the 
		other hand, since the matrix elements of the $X_N^{(i)}$'s are independent
		Guassians whose variances add, we see that the rescaled 
		sum $S_N$ is itself an $N \times N$ GUE random 
		matrix for each finite $N$.  Thus we recover Wigner's semicircle law (for GUE 
		matrices) from the free central limit theorem.
						
		\subsection{$\text{GUE}+\text{deterministic}$}
		Let $X_N$ be an $N \times N$ GUE random matrix.
		Let $Y_N$ be an $N \times N$ deterministic Hermitian matrix whose spectral 
		measure $\nu_N$ converges weakly to a compactly supported probability 
		measure $\nu$.  Let $\sigma$ be the limit distribution of the random matrix
		$X_N+Y_N$.  Since $X_N,Y_N$ are asymptotically free, we have
		
			\begin{equation*}
				\sigma = \mu \boxplus \nu,
			\end{equation*}
			
		\noindent
		where $\mu$ is the Wigner semicircle.
				
		\subsection{$\text{randomly rotated}+\text{diagonal}$}
		Consider the $2N \times 2N$ diagonal matrix
		
			\begin{equation*}
				D_{2N} = \begin{bmatrix}
					1 & {} & {} & {} & {}\\
					{} & -1 & {} & {} & {}\\
					{} & {} & \ddots & {} & {}\\
					{} & {} & {} & 1 & {} \\
					{} & {} & {} & {} & -1
				\end{bmatrix}
			\end{equation*}
			
		\noindent
		whose diagonal entries are the first $2N$ terms of an alternating sequence of $\pm1$'s,
		all other entries being zero.  Let $U_{2N}$ be a $2N \times 2N$ CUE random matrix,
		and consider the random Hermitian matrix
		
			\begin{equation*}
				A_{2N} = U_{2N}D_{2N}U_{2N}^* + D_{2N}.
			\end{equation*}
			
		\noindent
		Let 
		$\mu_{2N}$ denote the spectral measure of $A_{2N}$.  We claim that $\mu_{2N}$ 
		converges weakly to the arcsine distribution 
		
			\begin{equation*}
				\mu(\mathrm{d}t) = \frac{1}{\pi\sqrt{4-t^2}}\mathrm{d}t, \quad t \in [-2,2],
			\end{equation*}
			
		\noindent
		as $N \rightarrow \infty$.
		
		Proof:  Set $X_{2N}=U_{2N}D_{2N}U_{2N}^*$ and $Y_{2N}=D_{2N}$.  Then 
		$X_N,Y_N$ is a random matrix model for a pair of free random variables $X,Y$
		each of which has the $\pm1$-Bernoulli distribution
		
			\begin{equation*}
				\frac{1}{2}\delta_{-1} + \frac{1}{2}\delta_{+1}.
			\end{equation*}
			
		\noindent
		Thus the limit distribution of their sum is 
		
			\begin{equation*}
				\text{Bernoulli} \boxplus \text{Bernoulli} = \text{Arcsine}.
			\end{equation*}
						
\bibliographystyle{amsplain}

\end{document}